\pdfoutput=1

\documentclass[3p]{elsarticle}

\usepackage[english]{babel}
\usepackage{moreverb}
\usepackage{subfigure}
\usepackage{float}
\usepackage{graphicx}
\usepackage{epstopdf}
\usepackage{hyperref} 

\usepackage{times}
\usepackage{pifont,latexsym,ifthen,theorem,rotating,calc}
\usepackage{amsfonts,amssymb,amsbsy,amsmath}

\makeatletter
\def\ps@pprintTitle{%
  \let\@oddhead\@empty
  \let\@evenhead\@empty
  \def\@oddfoot{\reset@font\hfil\thepage\hfil}
  \let\@evenfoot\@oddfoot
}
\makeatother

\newcommand\BibTeX{{\rmfamily B\kern-.05em \textsc{i\kern-.025em b}\kern-.08em
T\kern-.1667em\lower.7ex\hbox{E}\kern-.125emX}}

\newcommand{\vek}[1]{\mathchoice{\displaystyle\boldsymbol#1}
{\textstyle\boldsymbol#1}{\scriptstyle\boldsymbol#1}
{\scriptscriptstyle\boldsymbol#1}}

\begin{document}

\begin{frontmatter}
\title{Conformal higher-order remeshing schemes\\ for implicitly defined
interface problems}

\author{S. Omerovi{\'c} and T.P. Fries}
\address{Institute of Structural Analysis, Graz University of Technology,
Lessingstr. 25/II, 8055 Graz, Austria}

\begin{abstract}
A new higher-order accurate method is proposed that combines the advantages
of the classical $p$-version of the FEM on body-fitted meshes with
embedded domain methods. A background mesh composed by higher-order
Lagrange elements is used. Boundaries and interfaces are described
implicitly by the level set method and are within elements. In the
elements cut by the boundaries or interfaces, an automatic decomposition
into higher-order accurate sub-elements is realized. Therefore, the
zero level sets are detected and meshed in a first step which is called
reconstruction. Then, based on the topological situation in the cut
element, higher-order sub-elements are mapped to the two sides of
the boundary or interface. The quality of the reconstruction and the
mapping largely determines the properties of the resulting, automatically
generated conforming mesh. It is found that optimal convergence rates
are possible although the resulting sub-elements are not always well-shaped.
\end{abstract}


\end{frontmatter}

\section{Introduction\label{sec:Introduction}}

Numerical simulations in all fields of computational mechanics often
involve complex-shaped, possibly curved domains and/or interfaces
with arbitrary geometry. Examples and applications are numerous and
include e.g.~in solid mechanics, the modelling of composite materials
with complex microstructure \cite{Moes_2003a}, material inhomogeneities
\cite{Sukumar_2001c}, solidification \cite{Chessa_2002a} and topology
optimization \cite{Belytschko_2003b,Wang_2004a}. Further applications
can be found in fluid mechanics and include interfaces in multi-phase
flows \cite{Glowinski_1999a} or problems in multi-physics. Such applications
feature non-smooth solutions and involve differential equations with
discontinuities in the coefficients \cite{Babuska_1970a} or even
different PDEs across the interfaces. It is well-known that the accuracy
of FEM approximations is severely degraded if meshes are used with
elements that do not conform to the boundaries and interfaces.

Therefore, over the years a lot of attention has been drawn to the
development of techniques to handle interface problems. Several different
approaches have been proposed and in terms of mesh representation,
the methods can roughly be divided into two different categories.
The first category consists of \emph{body fitted finite element methods},
where all interfaces and boundaries align with element edges. The
second category involves the large spectrum of \emph{embedding domain
methods} or \emph{fictitious domain methods} that use a non-conforming
(background) mesh and enforce boundary and interface conditions by
tailored approaches, e.g.~Lagrange multipliers, penalty methods or
Nitsche's method. 

In body fitted finite element methods, the computational meshes are
usually generated in two different ways. The first approach is to
create a mesh with straight edges and then transform this linear discretization
into higher-order elements, see \cite{Sherwin_2002a,Shephard_2005a,Luo_2004a,Lu_2014a}.
The other possibility is to construct a higher-order mesh directly
by creating curved boundaries and using them to discretize the interior
of the domain, see \cite{Sherwin_2002a}. Because the interface can
be arbitrarily complex, curved or even in motion, the generation of
interface-conforming meshes can be very costly from a computational
point of view. Even worse, a large amount of human interaction may
be required to generate suitable meshes, in particular with higher-order
elements. However, the outstanding advantage is that once a conforming
mesh is created, basically any finite element solver based on the
isoparametric concept can process it. No specialized treatment for
the boundary conditions is needed and standard Gaussian quadrature
rules are used for the numerical integration of the weak form. All
these advantages are also present in the approach proposed herein.

The second category, i.e.~embedding domain methods, often uses a
structured background mesh as a starting point, which does not account
for the element boundaries. The location of the interfaces may, for
example, be described via the level set method, see \cite{Osher_2003a,Osher_2001a,Belytschko_2003c}.
Methods characterized as embedding domain methods are e.g.~the immersed
boundary method \cite{Peskin_1972a,Cottet_2004a}, cut finite element
method \cite{Burman_2010a,Burman_2012a,Burman_2014a,Hansbo_2002a},
finite cell method \cite{Parvizian_2007a,Duester_2008a,Schillinger_2012a,Schillinger_2014a},
unfitted finite element methods \cite{Bastian_2009a,Lehrenfeld_2016a}
and the parametric finite element method \cite{Frei_2014a}. Even
the extended finite element method (XFEM) can be cast in the realm
of embedding domain methods for some specific applications, see e.g
\cite{Fries_2009b}, or in a higher order context \cite{Byfut_2009b,Legrain_2013a}.
In all these methods, the boundaries or interfaces are within elements
so that the application of boundary and interface conditions as well
as the numerical integration in cut elements is quite involved. Furthermore,
the conditioning of the system of equations may be degraded due to
largely different supports of the involved shape functions.

Herein, based on the method introduced in \cite{Noble_2010a} and
\cite{Kramer_2014a} in a linear, two-dimensional context, we propose
a higher-order accurate method which is located somewhere in-between
the two aforementioned methodologies. The method uses a background
mesh composed by higher-order elements; the mesh does not have to
be Cartesian nor structured. Elements that are cut by the implicitly
defined boundaries and interfaces are automatically decomposed into
conforming sub-elements. This decomposition is based on a two-step
procedure where, firstly, the zero levels of level set functions are
identified and meshed. Then, secondly, tailored mappings are introduced
which define conforming, higher-order sub-elements on the two sides
of the interface. Due to the similarities to the low-order method
proposed in \cite{Noble_2010a}, the introduced approach is labeled
a higher-order conformal decomposition finite element method (CDFEM).
Similar approaches that follow from mesh subdivision and mesh adaptation
are sometimes labelled cut-cell method are found in the context of
higher-order discontinuous Galerkin methods (DG) in \cite{Fidkowski_2007a}. 

The organization of this paper is as follows. Section \ref{sec:Overview}
describes some model problems and gives a short overview of the proposed
method. Section \ref{sec: Reconstruction} shows how to identify and
discretize the approximated zero level set by means of higher-order
one-dimensional interface elements. This is further referred to as
reconstruction. The employed element-wise subdivision technique is
discussed in detail in Section \ref{sec: Remeshing} and valid topologies
for triangular and quadrilateral elements are defined. Furthermore
coordinate mappings are discussed that define the location of the
internal nodes of the generated higher-order sub-elements. Section
\ref{sec:Numerical results} contains numerical examples illustrating
the accuracy and the stability of the proposed method. All results
are in excellent agreement with the analytical solutions and exhibit
optimal higher-order convergence rates. We conclude in Section \ref{sec:Conclusions}
and give an overview of open topics as well as an outlook on future
work.

\section{Proposed technique in a nutshell\label{sec:Overview}}

Let $\Omega$ be a domain embedded in the Euclidean space $\mathbb{R}^{2}$.
The domain is bounded by the external interface $\Gamma_{ext}$. Inside
$\Omega$, there may be internal interfaces $\Gamma_{int}$ separating,
for example, regions with different material properties or the two
sides of a crack path. The interfaces $\Gamma_{int}$ and $\Gamma_{ext}$
are sharp and possibly curved.

\begin{figure}[H]
\centering
\subfigure[Conforming mesh]{\includegraphics[width=6.9cm]{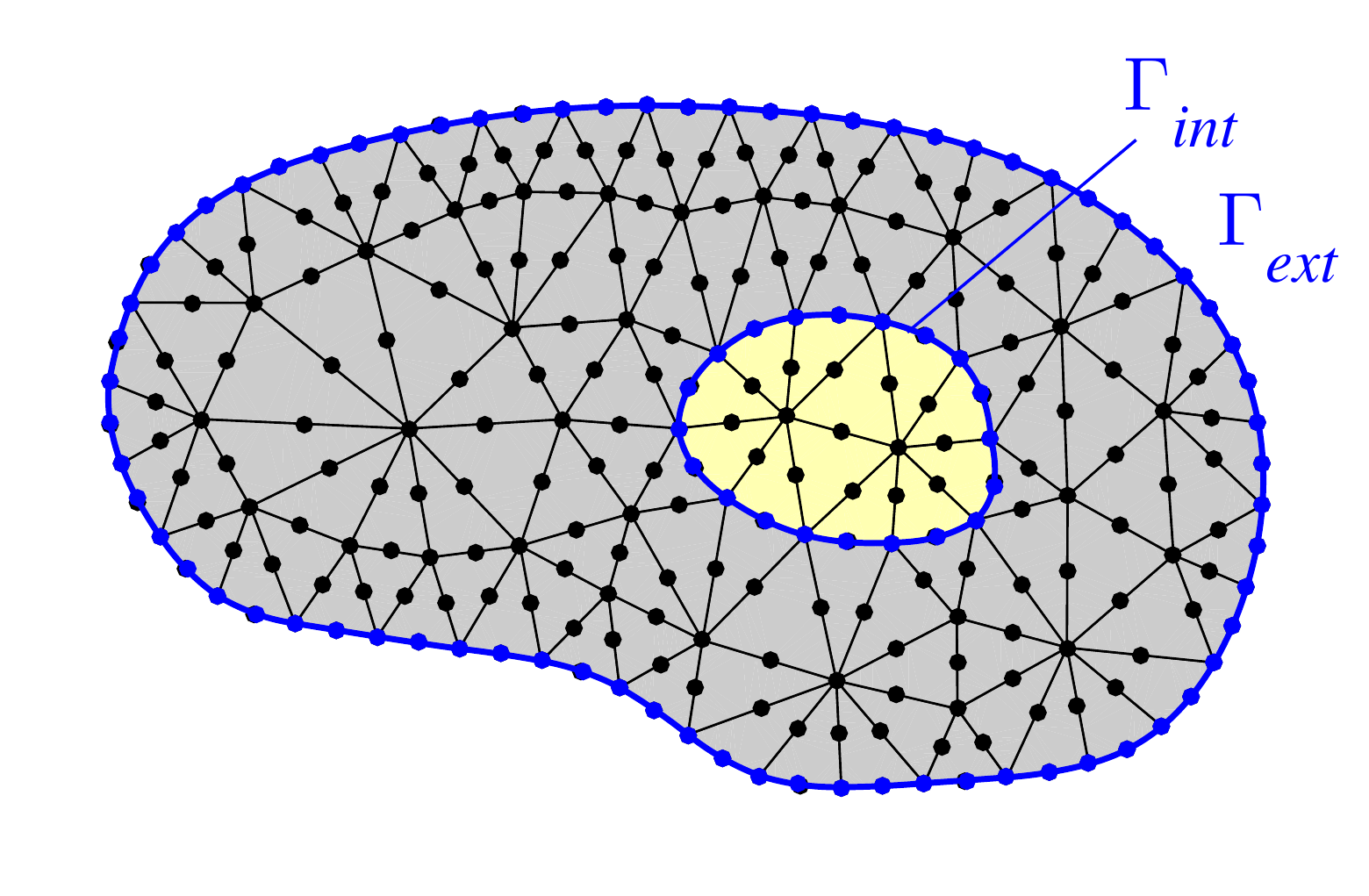}}$\quad$\subfigure[Automatically decomposed background mesh]{\includegraphics[width=6.9cm]{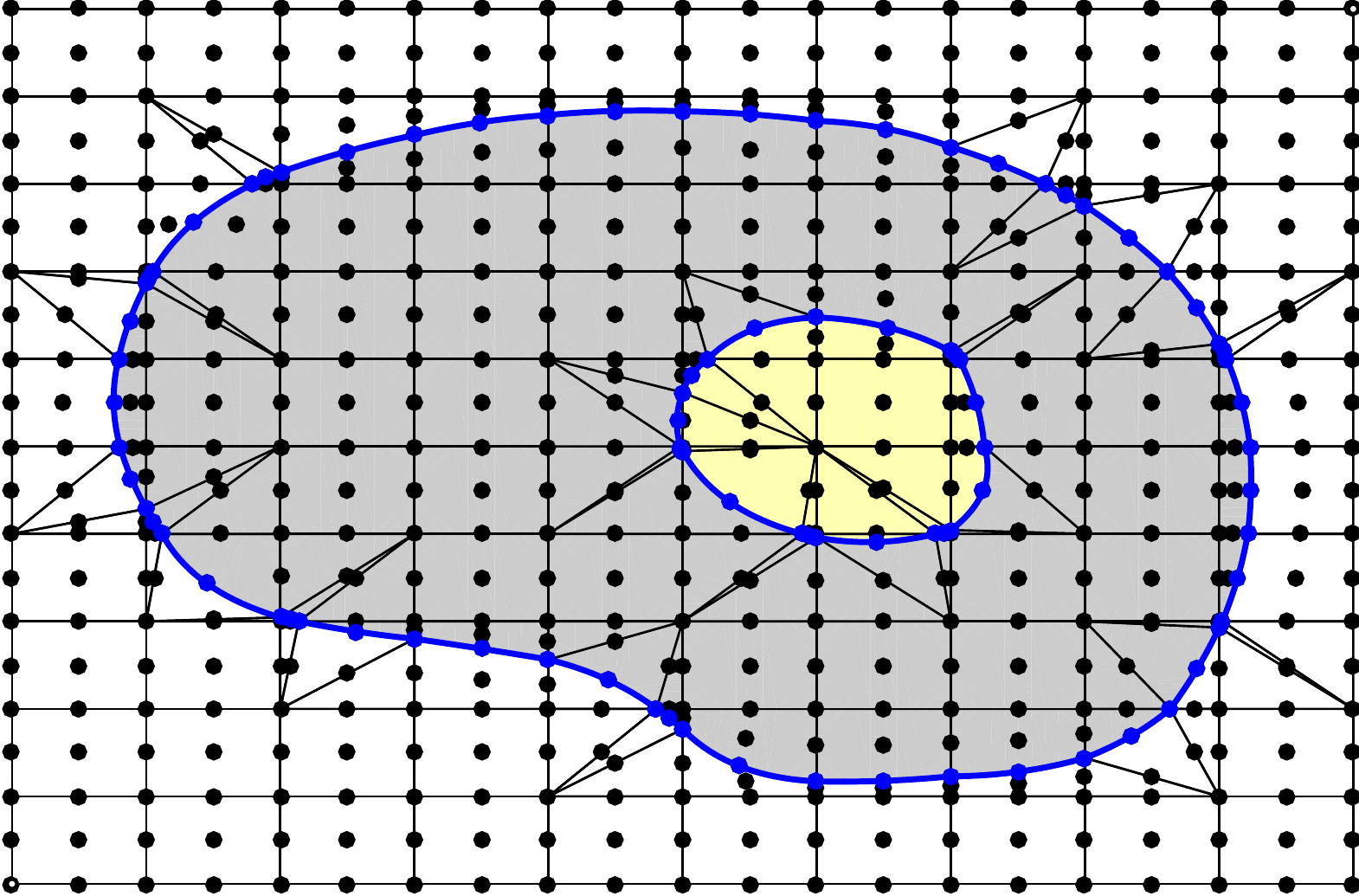}}

\caption{(a) Classical conforming mesh and (b) the mesh automatically generated
from decomposing the background mesh. Both meshes feature higher-order
(here: quadratic) Lagrange elements. \label{fig:Chap1MeshedDomains}}
\end{figure}

For the solution of a BVP in such a domain with the classical FEM,
it is crucial to construct a conforming mesh in order to achieve optimal
convergence properties. Thus, the elements must consider the position
of the interfaces, i.e.~element edges must align there. Furthermore,
no hanging nodes are typically accepted. A mesh for some schematic
domain is shown in Fig.~\ref{fig:Chap1MeshedDomains}(a). It is not
surprising that the generation of the mesh and the maintenance for
possibly evolving interfaces are not trivial tasks, especially not
with higher-order elements..

Therefore, the aim is to use simple higher-order meshes, further on
referred as background meshes, which do not have to account for any
of the (external and internal) interfaces. These are composed by standard
Lagrange elements which may be triangles or quadrilaterals. These
background meshes are simple to construct, but the proposed approach
is not necessarily restricted to Cartesian nor structured grids. We
propose a procedure to automatically decompose elements in the background
mesh that are cut by an arbitrary interface into isoparametric, higher-order
sub-elements. This leads to meshes as shown in Fig.~\ref{fig:Chap1MeshedDomains}(b).
It will be seen later that approximations with optimal accuracy are
possible on these meshes, although the element shapes are often quite
awkward. 

Herein, it is assumed that the position of the interfaces are defined
by the level set method. That is, a scalar function $\phi(\vek x)$
is introduced whose zero level set defines the position of the interfaces.
In the present work we restrict ourselves to the case of \emph{one}
level set function which however enables the definition of several
unconnected interfaces. The level set values are prescribed at the
nodes of the higher-order background mesh and interpolated by the
finite element shape functions in between. The discretized level set
function is further on referred to as $\phi^{h}$.

\begin{figure}[H]
\centering\includegraphics[scale=0.83]{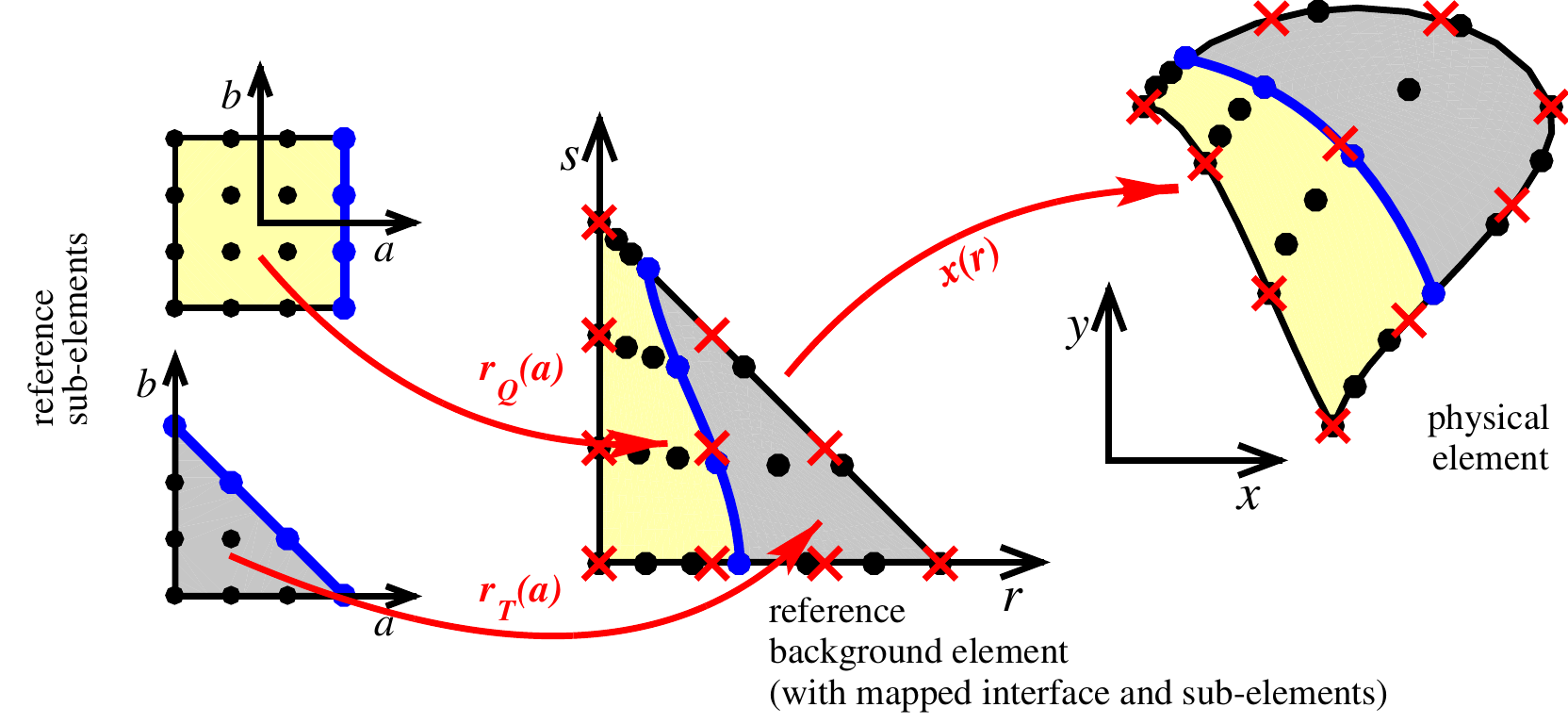}

\caption{\label{fig:Chap1Overview}This sketch demonstrates the two mappings
(coordinate transformations) that are needed for the proposed remeshing
scheme. The mapping $\vek x\left(\vek r\right)$ is defined by the
background element. The mappings $\vek r_{Q}\left(\vek a\right)$
and $\vek r_{T}\left(\vek a\right)$ are both based on the interface
element. (black dots: nodes of the decomposed sub-elements; red crosses:
nodes of the background element; blue dots: nodes of the interface
element) }
\end{figure}

The procedure for the automatic decomposition of cut elements is sketched
as follows: Following the work of \cite{Fries_2015a}, the zero level
set of $\phi^{h}$ is first detected and meshed by higher-order interface
elements. This is done in the background reference elements as depicted
by the blue line in Fig.~\ref{fig:Chap1Overview} and is called \emph{reconstruction}.
Here, we moderately modify the procedure to achieve the goals of this
work. The second step examines the decomposition into higher-order
sub-elements which align with the reconstructed level set. Therefore,
customized mappings $\vek r_{Q}(\vek a)$ and $\vek r_{T}(\vek a)$
for quadrilaterals and triangles, respectively, are introduced from
the reference elements to the sub-elements in the cut reference background
element. Finally, the reference background element is mapped to the
physical domain using the isoparametric mapping, denoted as $\vek x(\vek r)$.

\section{Reconstruction of the zero level sets \label{sec: Reconstruction}}

The aim of this section is to approximate the zero level set of the
level set function by means of higher-order finite elements, called
interface elements herein. As this paper is restricted to two dimensions,
the discretized zero level set is composed by \emph{line} elements.
The topic of reconstructing zero level sets is discussed in more detail
in \cite{Fries_2015a}, including the extension to the three-dimensional
case.

The strategy has a local character, that is, the reconstruction of
the zero isoline advances element by element. The choice of Lagrange
elements enforces continuity between the elements. The first step
is to transfer the level set data, given in the physical domain, to
the corresponding reference domain, as the reconstruction is executed
in the reference element. This is particularly simple for Lagrange
elements as the nodal values $\phi(\vek x_{n})$ are taken as nodal
values $\phi(\vek r_{i})$ in the reference element. Thereby, the
indices $n$ and $i$ denote a specific node in global and local element
numbering, respectively. Based on this information, elements that
are cut by the zero level set are detected. Once the cut elements
are identified, it must be ensured that the level set data is valid
in the sense that a unique subdivision of the element is possible.
For topologically valid cut elements, a scalar root finding algorithm
is provided, giving any desired number of points on the interface.
These points serve as nodes of the higher-order interface elements
approximating the zero level set.

\subsection{Detection of cut elements}

An element is considered as cut, when the data at a finite number
of points indicates that the approximated zero level set function
$\phi^{h}$ crosses the element domain. To decide whether the element
is cut, a sample grid with a user defined resolution is used in the
parent element, see Fig.~\ref{fig:RootSearchEdges}. The resolution
of this auxiliary grid should depend on the number of element nodes.
For the calculations performed in this work three points between element
nodes in each dimension were chosen. Using the element functions and
the nodal values $\phi^{h}(\vek r_{i}),i=\{1,..n\}$, the level set
is interpolated at all grid points. The set of grid points assigned
to an element is defined as $I_{\mathrm{grid}}$. The check for cut
elements is then done by an inspection of the sign of the grid level
set values. The considered element is cut if
\begin{equation}
\min_{i\in I_{\mathrm{grid}}}(\phi_{i})\,\cdot\,\max_{j\in I_{\mathrm{grid}}}(\phi_{j})<0.\label{eq:LevelSetCheck}
\end{equation}
It is important to note that nodal values are not sufficient for the
check due to the non-positiveness of the Lagrange shape functions.
Hence, a sample grid with suitable resolution is crucial for the check.

\subsection{Valid level set data}

Arbitrary level set data in a higher-order element may result in very
complex patterns of zero level sets and the implied regions where
$\phi(\vek x)$ is positive or negative, respectively. Therefore,
it is crucial to impose restrictions on the complexity of the level
set data. Herein, we call the level set data in a cut element according
to Eq.~(\ref{eq:LevelSetCheck}) \emph{valid} if (i) the element
boundary is cut \emph{two} times and (ii) the zero level set divides
the element domain $\Omega_{\vek r}^{h}$ into \emph{two} sub-domains.
However, even valid level set data has to meet some restrictions regarding
the curvature of the zero level set, which is discussed in more detail
in Section~\ref{sec: Remeshing}. 

\begin{figure}
\centering

\subfigure[]{\includegraphics[width=3.5cm]{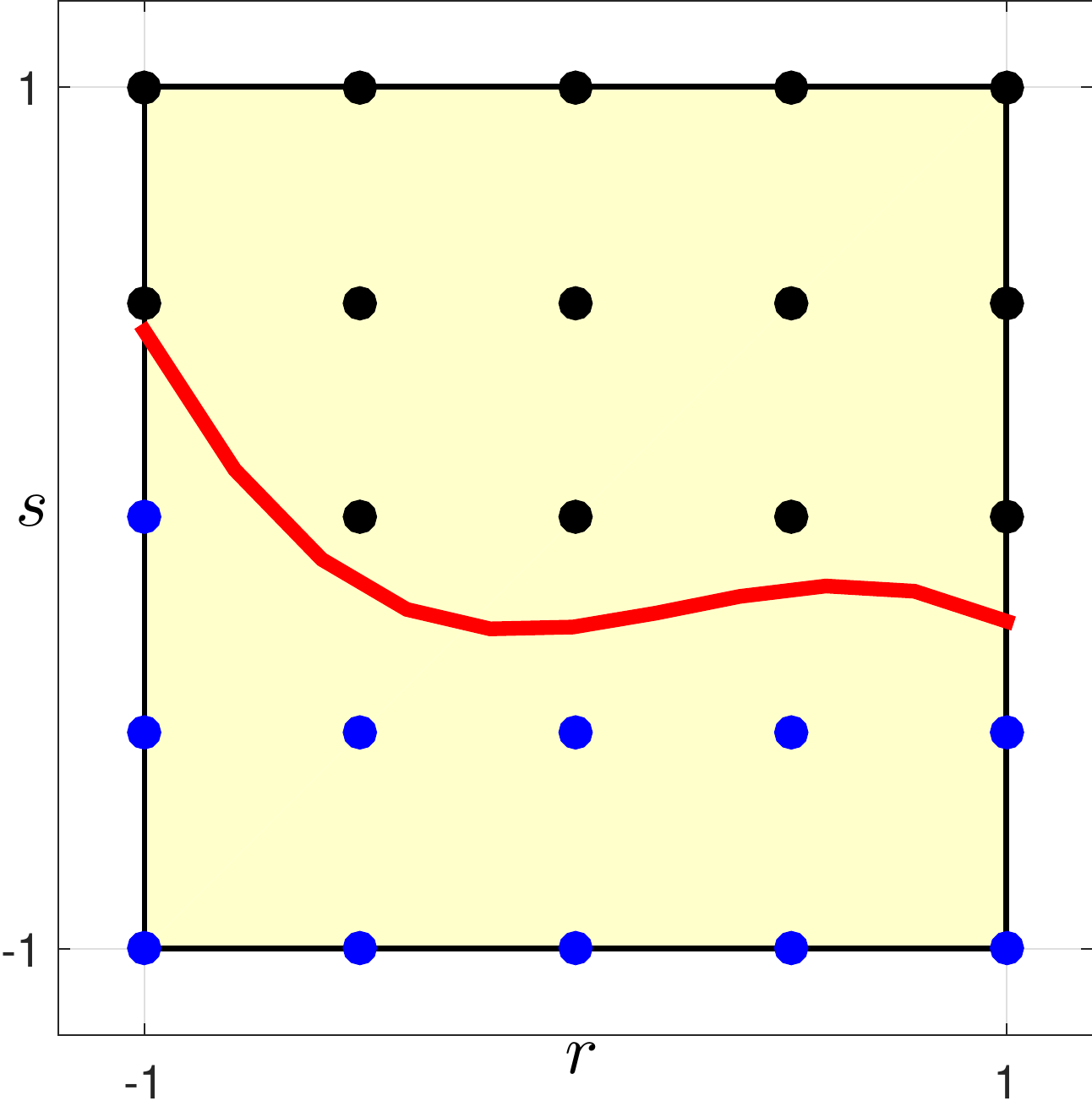}}\subfigure[]{\includegraphics[width=3.5cm]{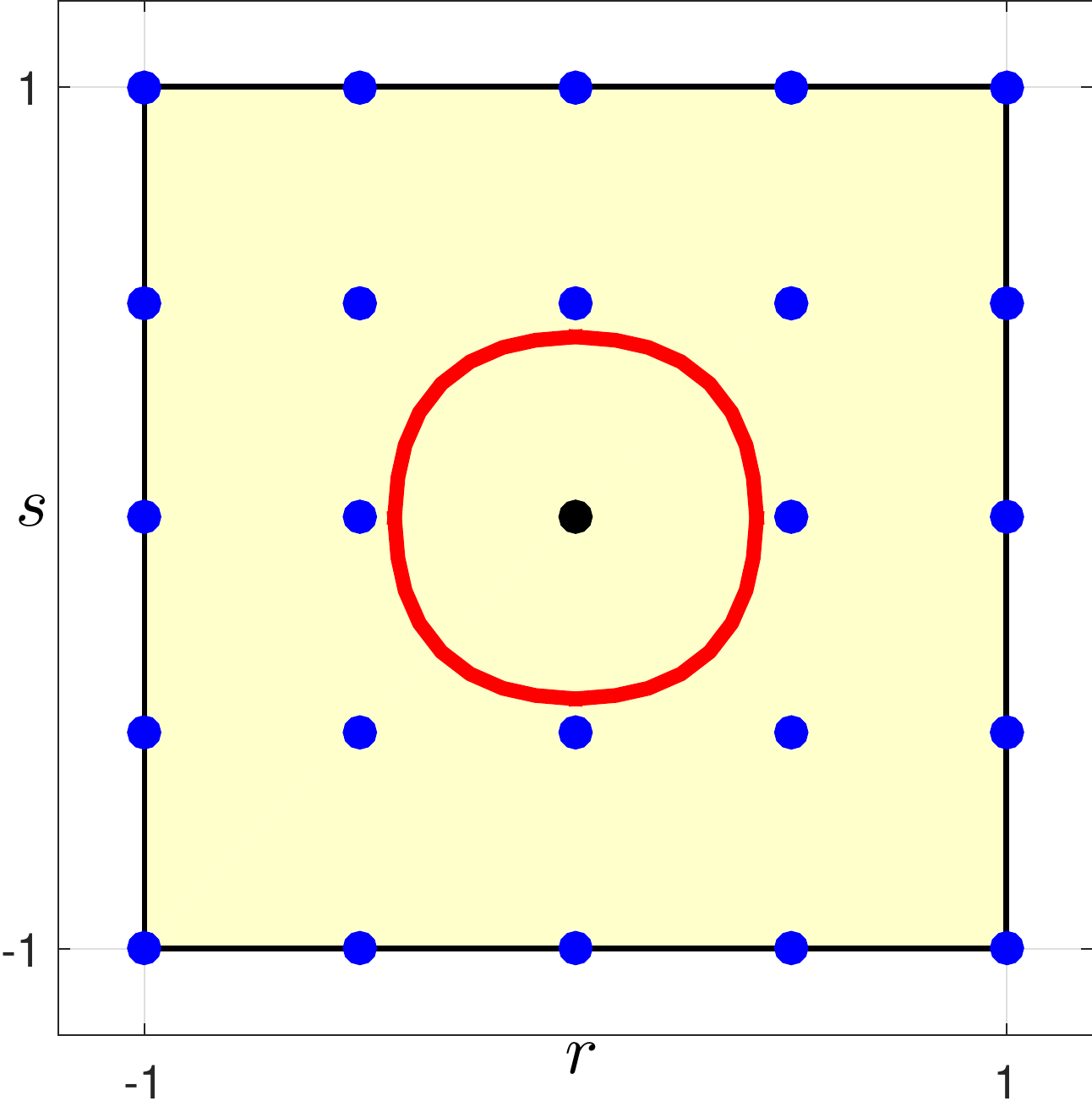}}\subfigure[]{\includegraphics[width=3.5cm]{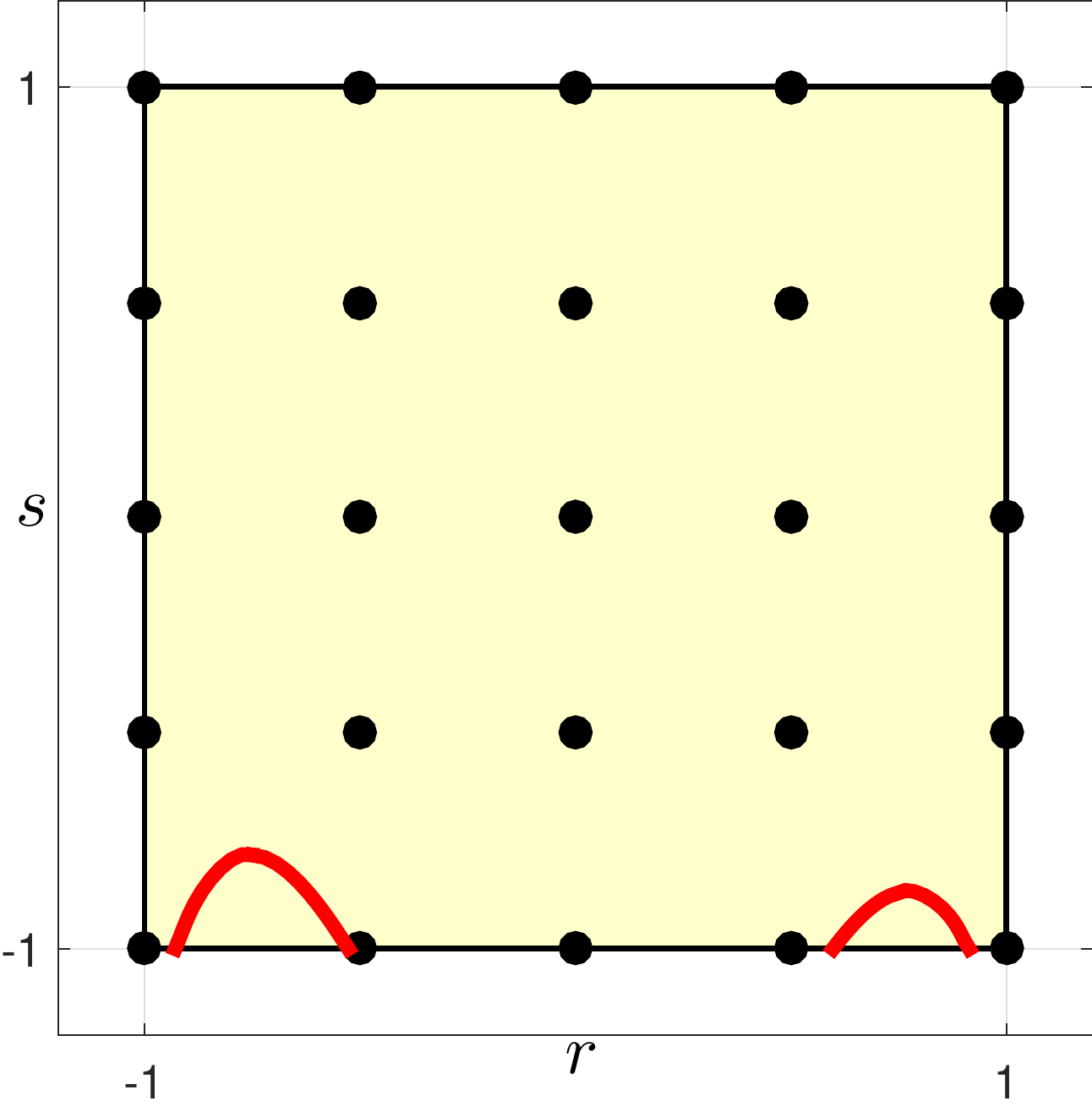}}\subfigure[]{\includegraphics[width=3.5cm]{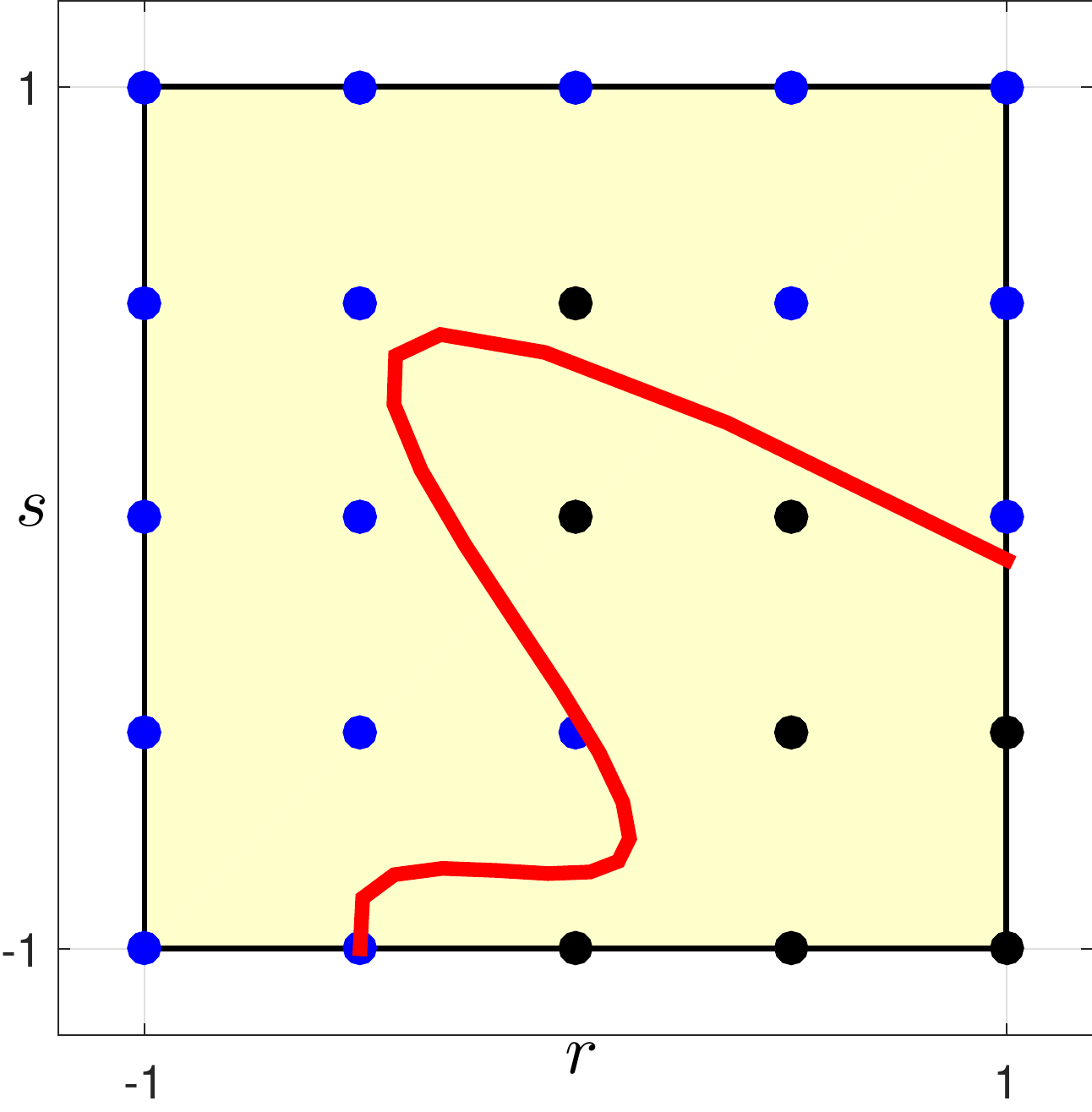}}

\caption{\label{fig:ValidityIsoLines} Zero level sets in reference elements
are shown according to (a) valid level set data, (b) and (c) invalid
data. (d) results from valid level set data, however the reconstruction
is likely to fail.}
\end{figure}

Fig.~\ref{fig:ValidityIsoLines}(a) shows valid level set data subdividing
the reference element into two sub-domains. Inclusions as seen in
Fig.~\ref{fig:ValidityIsoLines}(b) violate the conditions, hence,
the level set data is not valid. A situation violating the second
condition is shown in Fig.~\ref{fig:ValidityIsoLines}(c). Here the
level set data changes more than two times over the element edges,
producing more than two separate domains. If there is a high curvature
locally, see Fig.~\ref{fig:ValidityIsoLines}(d), the level set data
is valid according to the conditions above, however, the reconstruction
may not be successful anyway (which is automatically detected). If
the element is regarded as cut but features invalid level set data,
we note that a (locally) refined background mesh always produces valid
level set data that meets the conditions from above. 

\subsection{Reconstruction of zero level set by Newton-Raphson procedure}

Assuming that the element is cut and the level set data is valid,
the task is now to locate element nodes of the interface element with
a user-defined order on the zero level set. The process starts with
finding two roots on the element boundary. Then, a number of points
$\vek P_{i}$, where $i$ depends on the desired element order, are
evenly distributed between the two roots on the boundary. These points
serve as the starting guesses for the iterative root-finding inside
the element.

\subsubsection{Root search on the element boundary}

\begin{figure}
\centering
\subfigure[Triangular element]{\includegraphics[width=5cm]{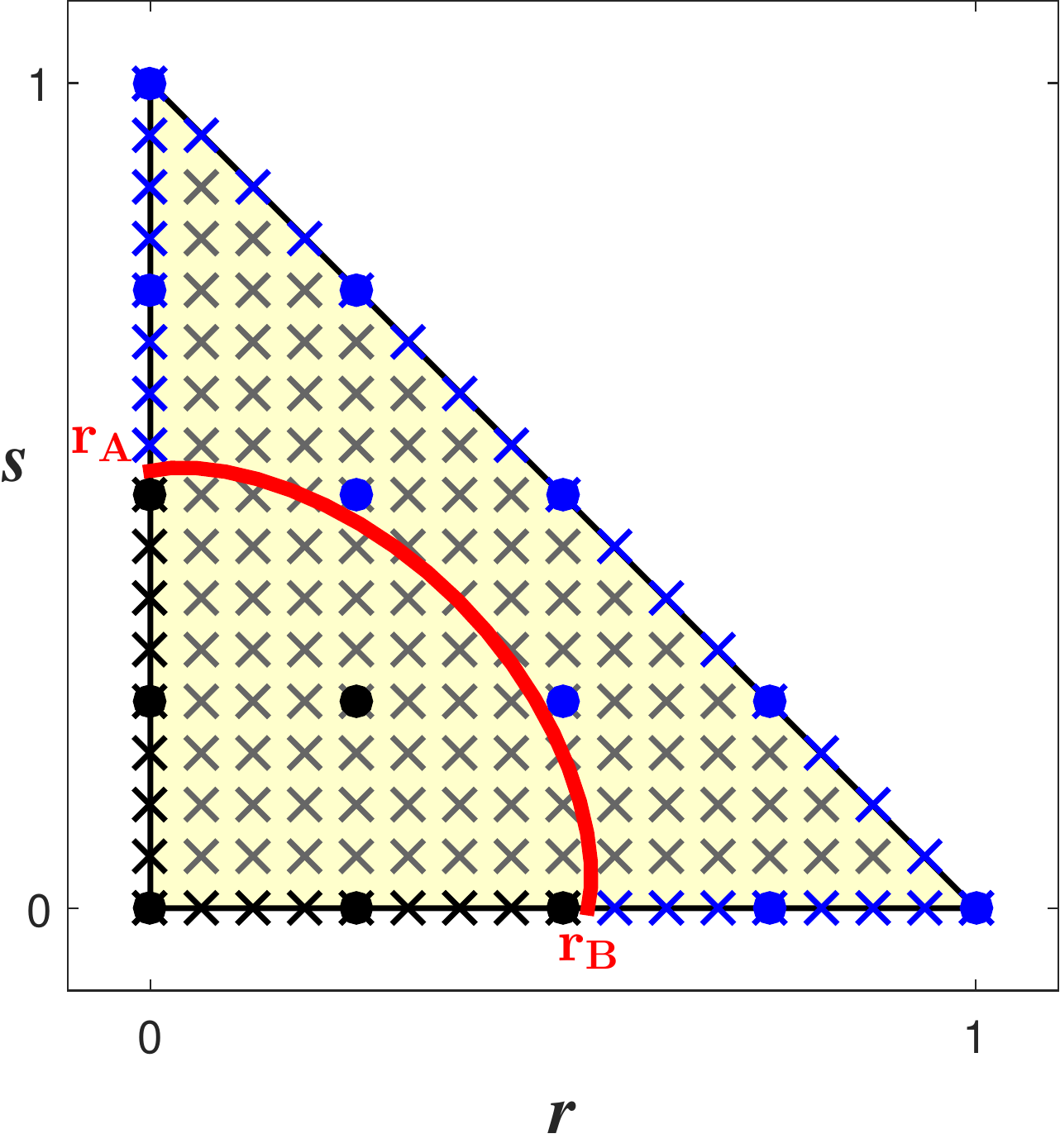}}$\qquad$\subfigure[Quadrilateral element]{\includegraphics[width=5cm]{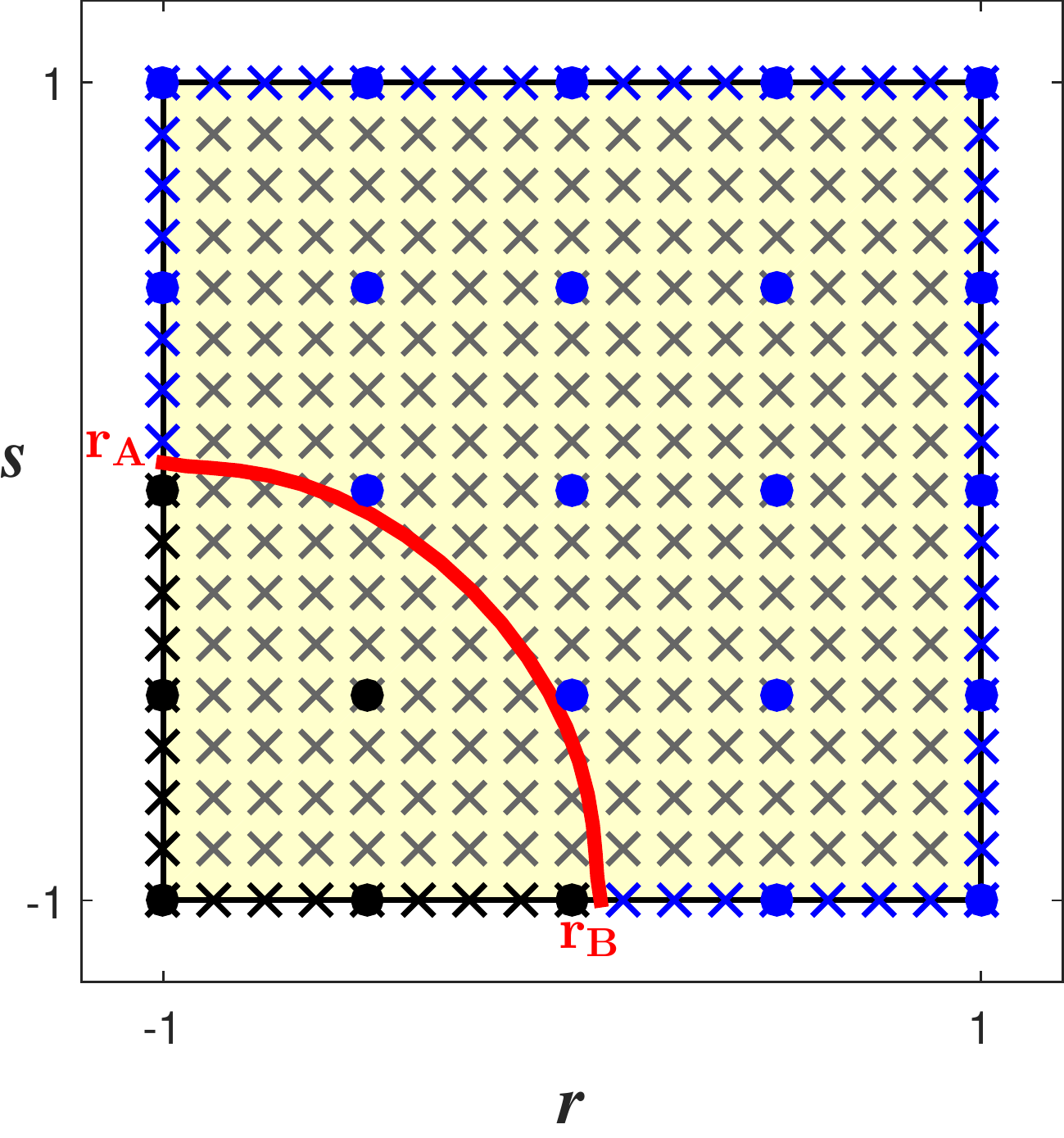}}

\caption{\label{fig:RootSearchEdges}Detection of the intersections of the
zero level set with the element boundary based on a sample grid. (blue
and black dots: element nodes with positive and negative level-set
values, respectively. The crosses build the sample grid) }
\end{figure}

The procedure starts with finding the two intersections of the zero
level set with the element boundary. Therefore, the sample grid from
above is considered on the boundary, see Fig.~\ref{fig:RootSearchEdges}.
Based on the signs of the level set data in two adjacent grid points,
intervals containing zeros of $\phi^{h}$ can be located immediately.
The zeros are labeled as $\vek r_{A}$ and $\vek r_{B}$ and are found
by applying an iterative root finding algorithm along the edge. The
middle of the interval is chosen as a starting guess for the iteration
and the algorithm is given by
\begin{equation}
\vek r^{i+1}=\vek r^{i}-\frac{\phi^{h}\left(\vek r\right)}{\nabla_{\vek r}\phi^{h}(\vek r^{j})}.\label{eq:NREdge}
\end{equation}
For all considered topological cases, a cut element \emph{edge} has
either one or two roots, depending on the topological situation. In
this section we restrict ourselves to the case of an edge having only
one root. As elaborated further in Section \ref{sec: Remeshing} this
leads to a topological case labeled as \emph{local topology}, which
is the starting point of the subdivision scheme. The configuration
with two roots on one element edge, designated as \emph{global topology,}
is also described in Section \ref{sec: Remeshing}. In this case the
approach is different, but for the reconstruction of the zero level
set reduces to the same procedure as described here. 

\subsubsection{Root search inside an element}

\begin{figure}
\centering
\subfigure[]{\includegraphics[width=5cm]{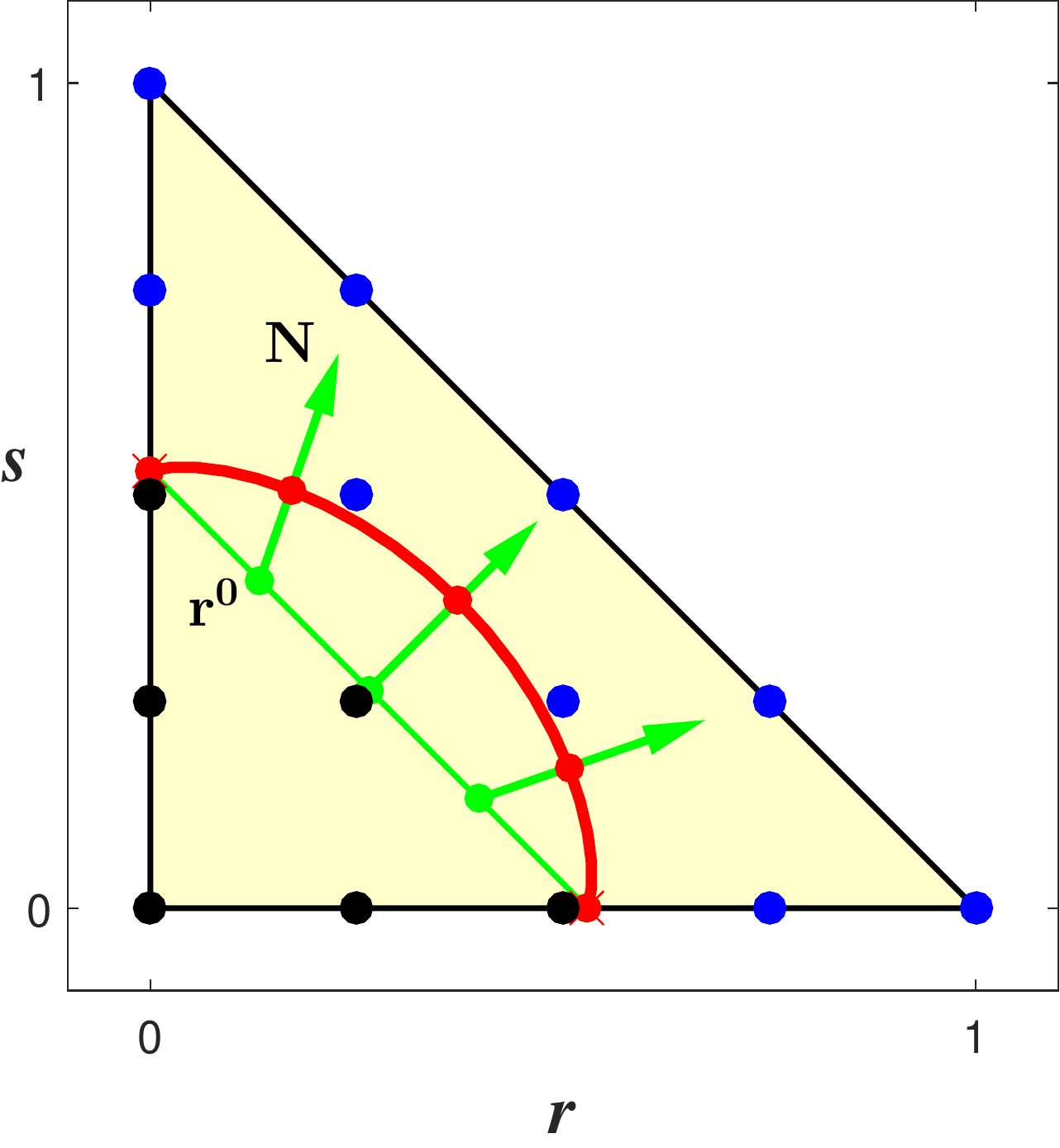}}$\qquad$\subfigure[]{\includegraphics[width=5cm]{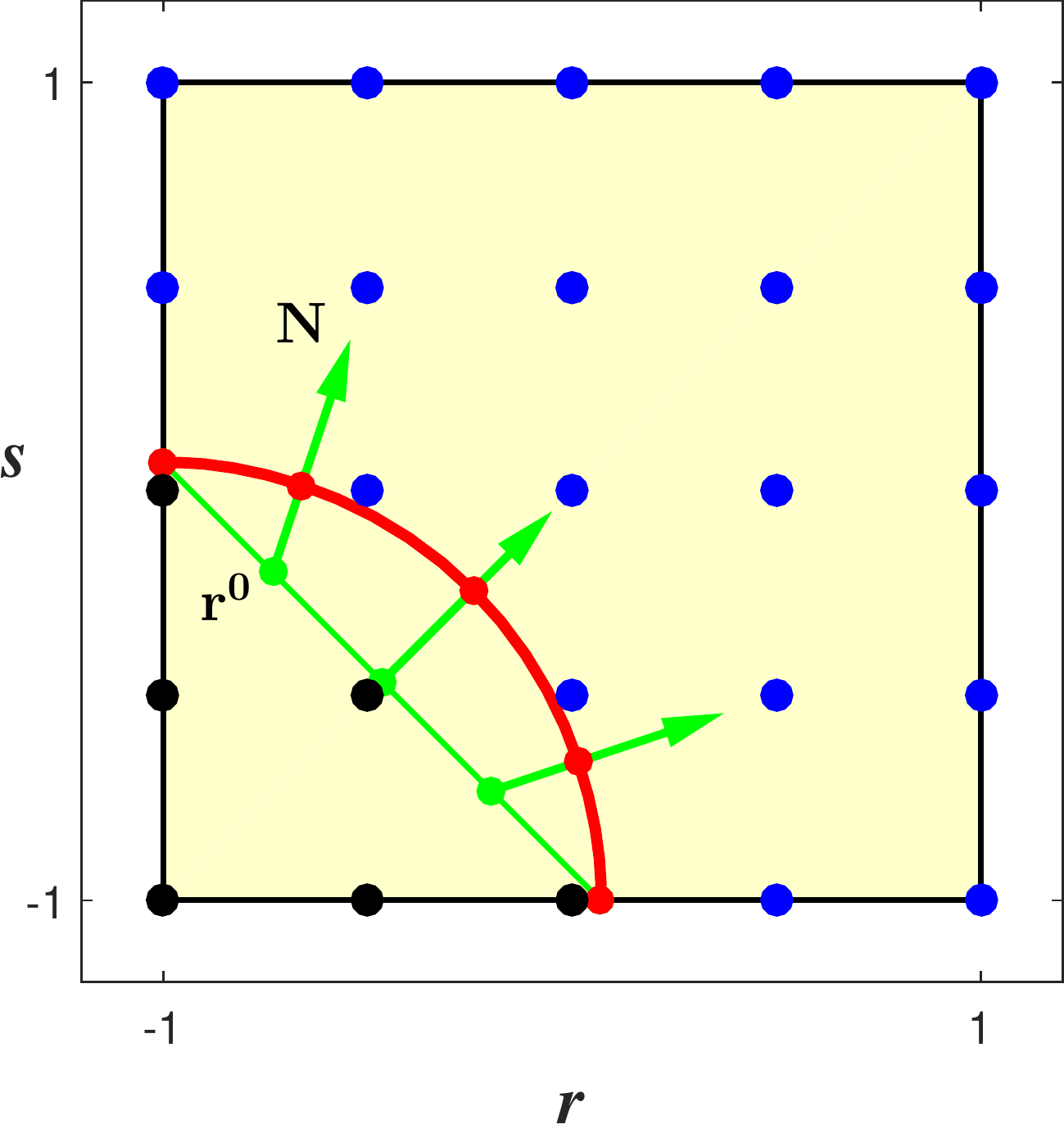}}

\caption{Search for the roots of $\phi^{h}$ inside the reference element,
using a search path specified by the starting point $\vek P=\vek r^{0}$
and the search direction $\vek N$ for a triangular (a) and a quadrilateral
(b) element.\label{fig:RootSearchElem}}
\end{figure}

For a detailed discussion concerning the root search inside an element,
with different choices of start points and search directions, we
refer to \cite{Fries_2015a}. Once the root search along the boundary
is done, the two roots are connected by a straight line. On this line,
a number of nodes $\vek P_{i}$, associated to the desired order of
the interface element are distributed equally spaced. This gives a
rough approximation of the zero isoline, and is further referred to
as ``linear reconstruction'', depicted by the green line in Fig.~\ref{fig:RootSearchElem}.
The root search is then executed along a straight path, going through
the start point $\vek P=\vek r^{0}$ with a direction $\vek N$ given
by the gradient of the level set function at $\vek P$, i.e.~$\vek N=\nabla_{\vek r}\phi^{h}(\vek P)$.
Applying a Newton-Raphson approach along $\vek N$, similar to Eq.~(\ref{eq:NREdge})
yields
\begin{equation}
\vek r^{i+1}=\vek r^{i}-\frac{\phi^{h}\left(\vek r^{i}\right)}{\nabla_{\vek r}\phi^{h}(\vek r^{j})\cdot\vek N}\cdot\vek N\label{eq:NR3}
\end{equation}
After the residuum $|\vek r^{i+1}-\vek r^{i}|$ of Eq.~(\ref{eq:NR3})
falls below a user-specified convergence tolerance (in the performed
calculations this tolerance is set as $\approx10^{-12}$) the newly
found coordinates of the roots are set as nodal coordinates of the
higher-order interface element, see the red line in Fig.~\ref{fig:RootSearchElem}.
This interface element, further on referred to as $\Gamma_{\vek r}^{h}$,
will finally serve as a new (curved) edge for higher order triangular
or quadrilateral elements on the two sides of the interface.

\clearpage

\section{Conform remeshing of a background element\label{sec: Remeshing}}

At this point we assume that a general, higher-order cut background
element $\Omega_{\vek r}^{h}$ is given and the zero level set is
successfully approximated by the reconstructed interface element $\Gamma_{\vek r}^{h}$,
see Section \ref{sec: Reconstruction}. The aim is to: firstly, to
identify all topological cases leading to low-order sub-cells with
possibly one higher-order side coinciding with $\Gamma_{\vek r}^{h}$
and, secondly, to find a map for the element nodes of higher-order
elements into these sub-cells, which is not trivial as seen below.
It is noted that the resulting decomposition features triangular and
quadrilateral elements.

\begin{figure}[H]
\centering

\includegraphics[width=13cm]{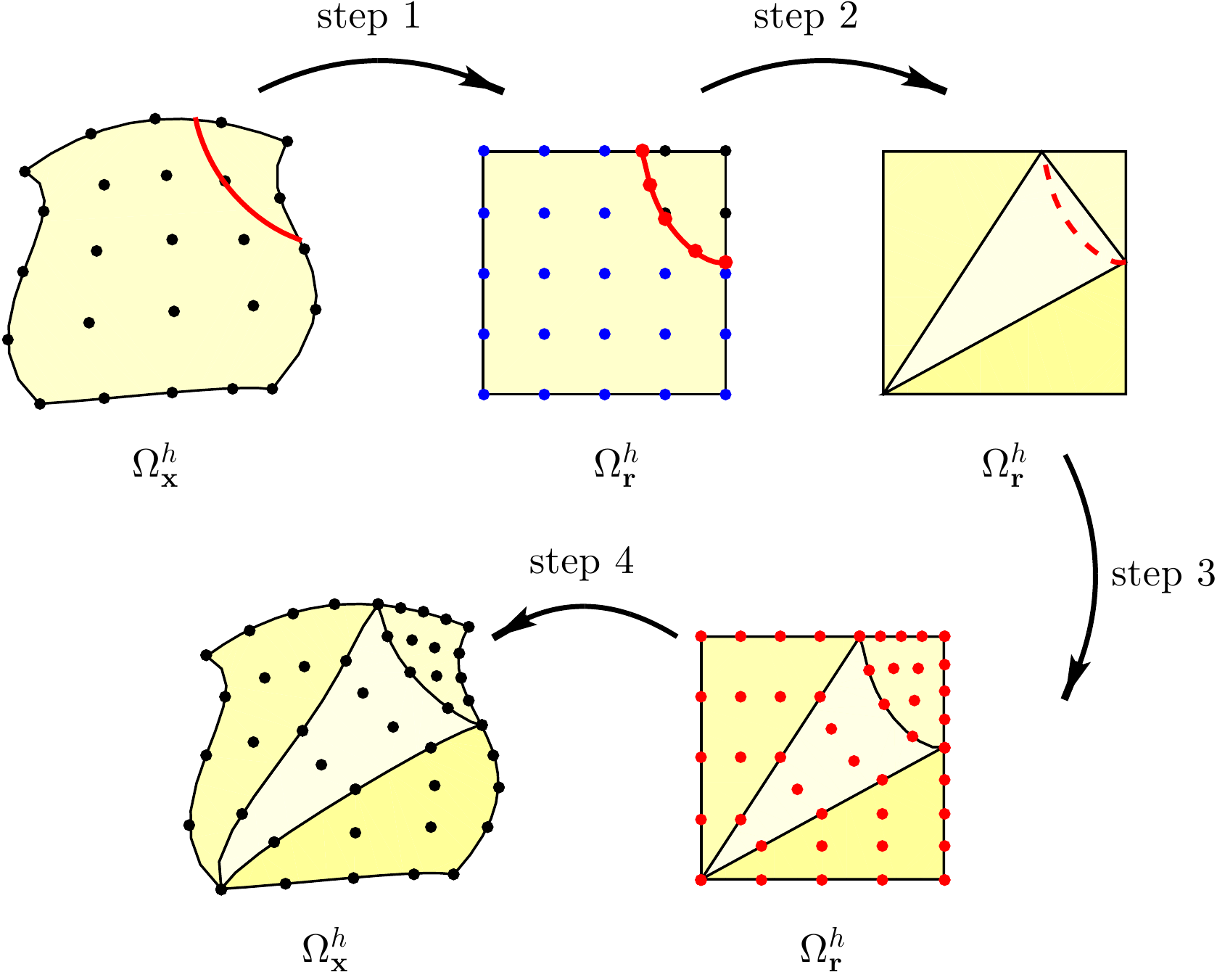}\caption{\label{fig:VisAllMappings}Different steps involved in the remeshing
strategy.}
\end{figure}

There are several mappings and transformations involved and Fig.~\ref{fig:VisAllMappings}
shows the complete process in one diagram. Emanating from the physical
element $\Omega_{\vek x}^{h}$, in \textbf{step\,1} the level set
values are transferred to the parent domain $\Omega_{\vek r}^{h}$,
where the reconstruction of the zero level set is realized. In \textbf{step\,2}
the subdivision of $\Omega_{\vek r}^{h}$ takes place. The element
is subdivided into a few predefined topologies depending on the cut
situation. The only missing link for a working computational mesh
is a node distribution within the element patch that preserves optimal
accuracy for higher-order elements. This is done by means of a customized
mapping procedure, sketched in \textbf{step\,3} of Fig.~\ref{fig:VisAllMappings},
wherefore the approach for triangular and quadrilateral elements differentiates
significantly. Once the nodal coordinates are determined, a conforming
higher-order element patch is obtained. At last, in \textbf{step\,4}
the data is projected back to the physical domain via an isoparametric
mapping.

\begin{figure}
\centering
\subfigure[Linear triangular element]{\includegraphics[width=5cm]{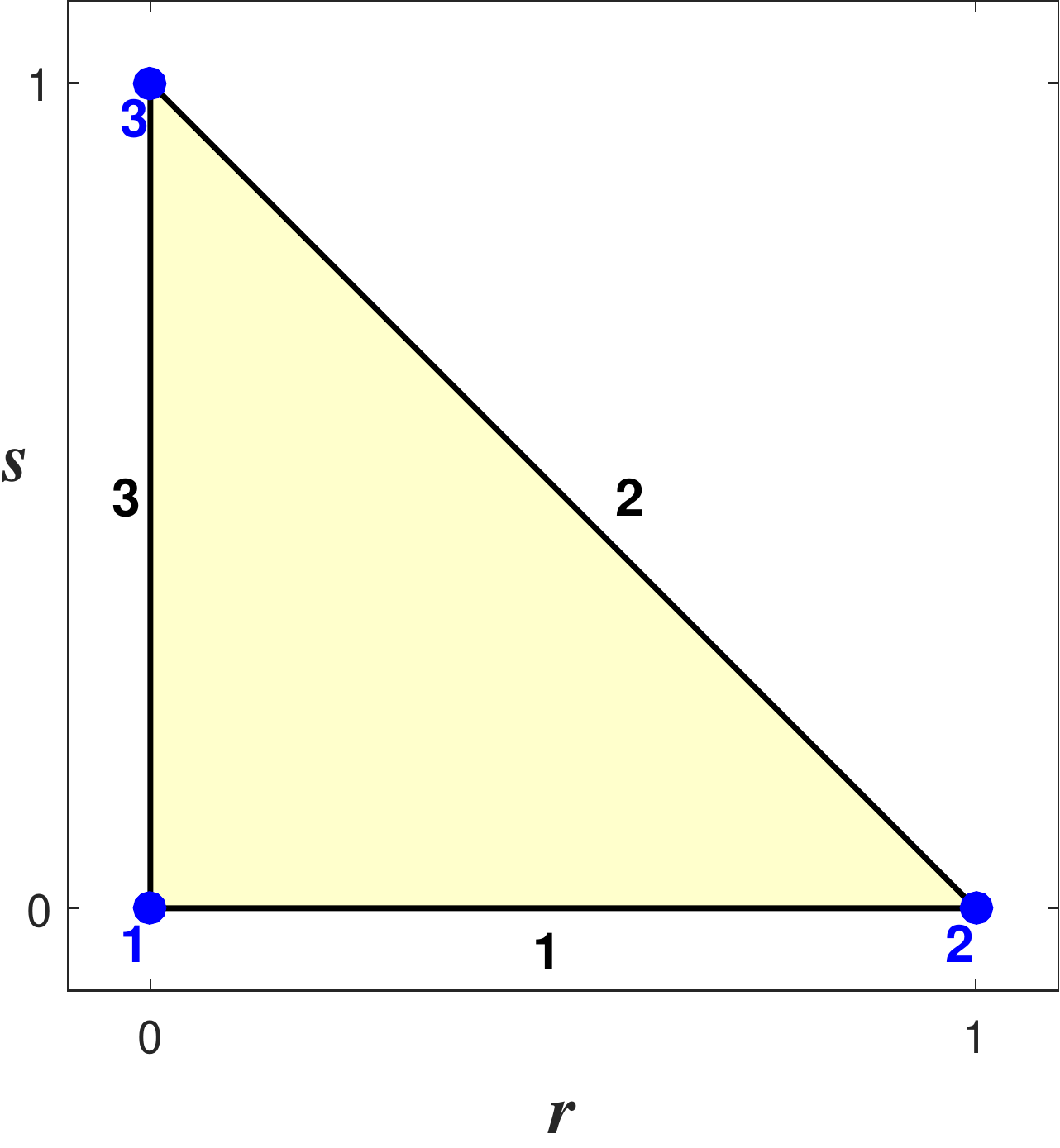}}$\qquad$\subfigure[Bi-linear quadrilateral element]{\includegraphics[width=5cm]{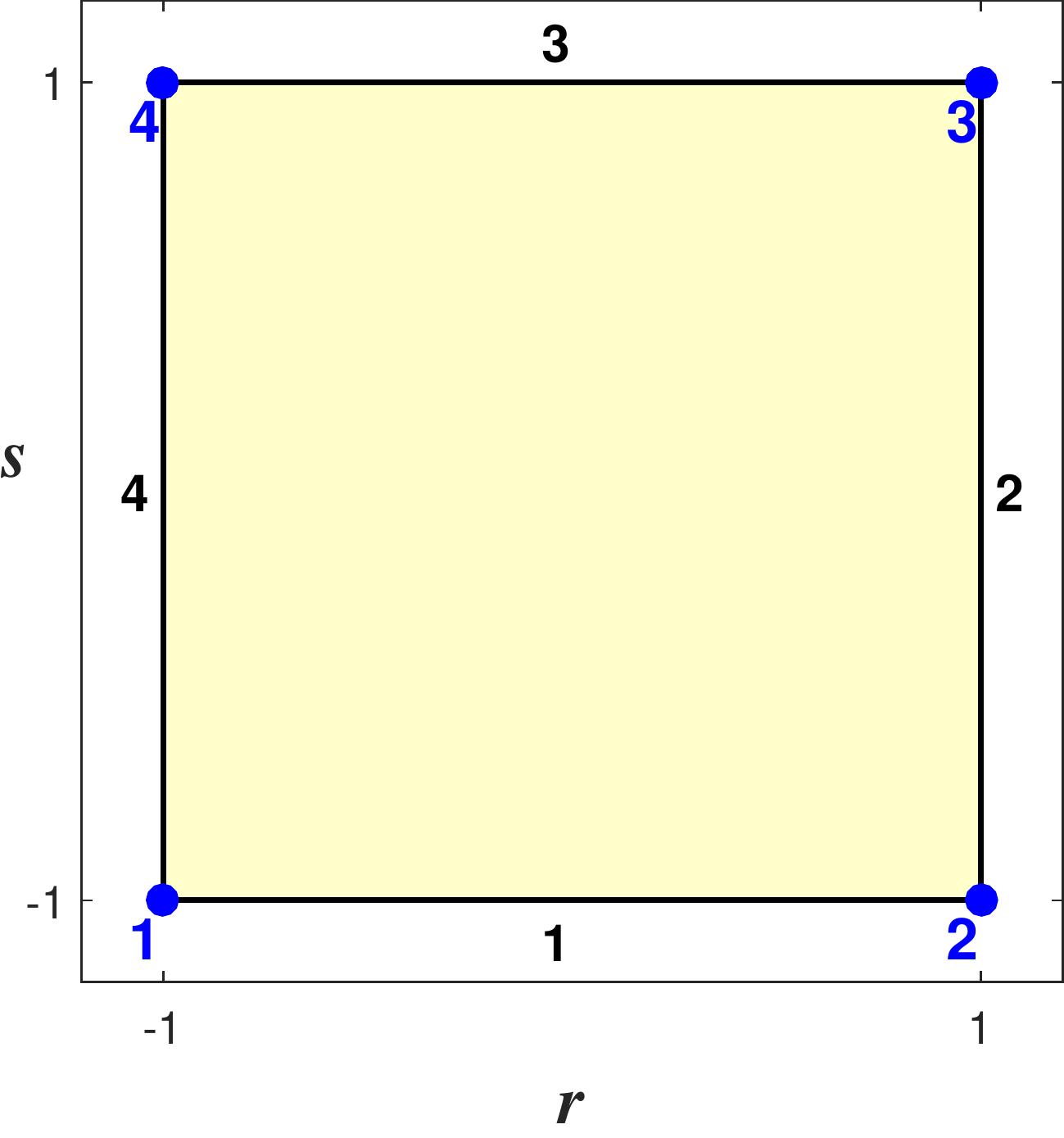}}

\caption{\label{fig:Chap4EdgeNodeNumbering}Edge and node numbering for (bi-)linear
Lagrange elements.}
\end{figure}

Regarding the nomenclature of the following topology classifications,
some abbreviations are introduced. The first letter describes the
type of element. This can be either \textbf{$\mathbf{T}$} for \textbf{t}riangular
or \textbf{$\mathbf{Q}$} for \textbf{q}uadrilateral elements followed
by a combination of subscripts and/or a superscripts. The index of
the subscripts gives the information which \emph{edges} are cut, whereas
the superscripts describes which \emph{nodes} are cut. So \emph{$\mathbf{Q}_{IJ}$}
denotes a quadrilateral, with a zero level set which cuts edge $I\le4$
and edge $J\le4$. Another example is $\mathbf{T}{}_{I}^{J}$, describing
a triangular element where one edge $I\le3$ and one node $J\le3$
are cut. The edge and node numbering are evident from Fig.~\ref{fig:Chap4EdgeNodeNumbering}.
Note that only the vertex nodes of a higher-order element are relevant
for the subdivision scheme.

The introduced notation has the advantage that it indicates directly
the number of possible topologies. As mentioned in Section \ref{sec: Reconstruction},
a valid cut element has exactly two intersections of the zero level
set along the whole element boundary. As a result, all possible cut
situations for triangles are described by $\mathbf{T}{}_{I}^{J}$,
$\mathbf{T}{}_{IJ}$ and $\mathbf{T}{}^{IJ}$, with $I\le3$ and $J\le3$.
The same is valid for quadrilaterals, $\mathbf{Q}{}_{I}^{J}$, $\mathbf{Q}{}_{IJ}$
and $\mathbf{Q}{}^{IJ}$, with $I\le4$ and $J\le4$. It is worth
recalling that a symmetry in the indices holds, e.g.~$\mathbf{T}{}_{IJ}=\mathbf{T}{}_{JI}$
which changes only the orientation of the normal vector but not the
topological situation itself. Consequently, the permutations of $I$
and $J$ for subscripts and superscripts cover all possible topological
situations. The only invalid situation is \emph{one} node that is
cut twice by the level set. In the proposed notation that would result
in $\mathbf{T}{}^{II}$ and $\mathbf{Q}{}^{JJ}$with $I\le3$ and
$J\le4$.

\subsection{Subdivision strategy for triangles}

A configuration is defined as ``local'' if the element can be subdivided
without influencing neighbour elements. For triangular elements there
are two \emph{local} configurations, that affect only the considered
element. Another three \emph{non-local} configurations are present,
influencing the subdivision of \emph{one} neighbour element. The non-local
scheme works successively which means that the (non-local) subdivision
scheme produces a local configuration. Once this is done, the local
subdivision algorithm yields the desired element patches where some
element nodes coincide with the reconstructed interface element.

\subsubsection{Local subdivisions}

\begin{figure}
\centering
\subfigure[$\mathbf{T}_{13}$]{\includegraphics[width=5cm]{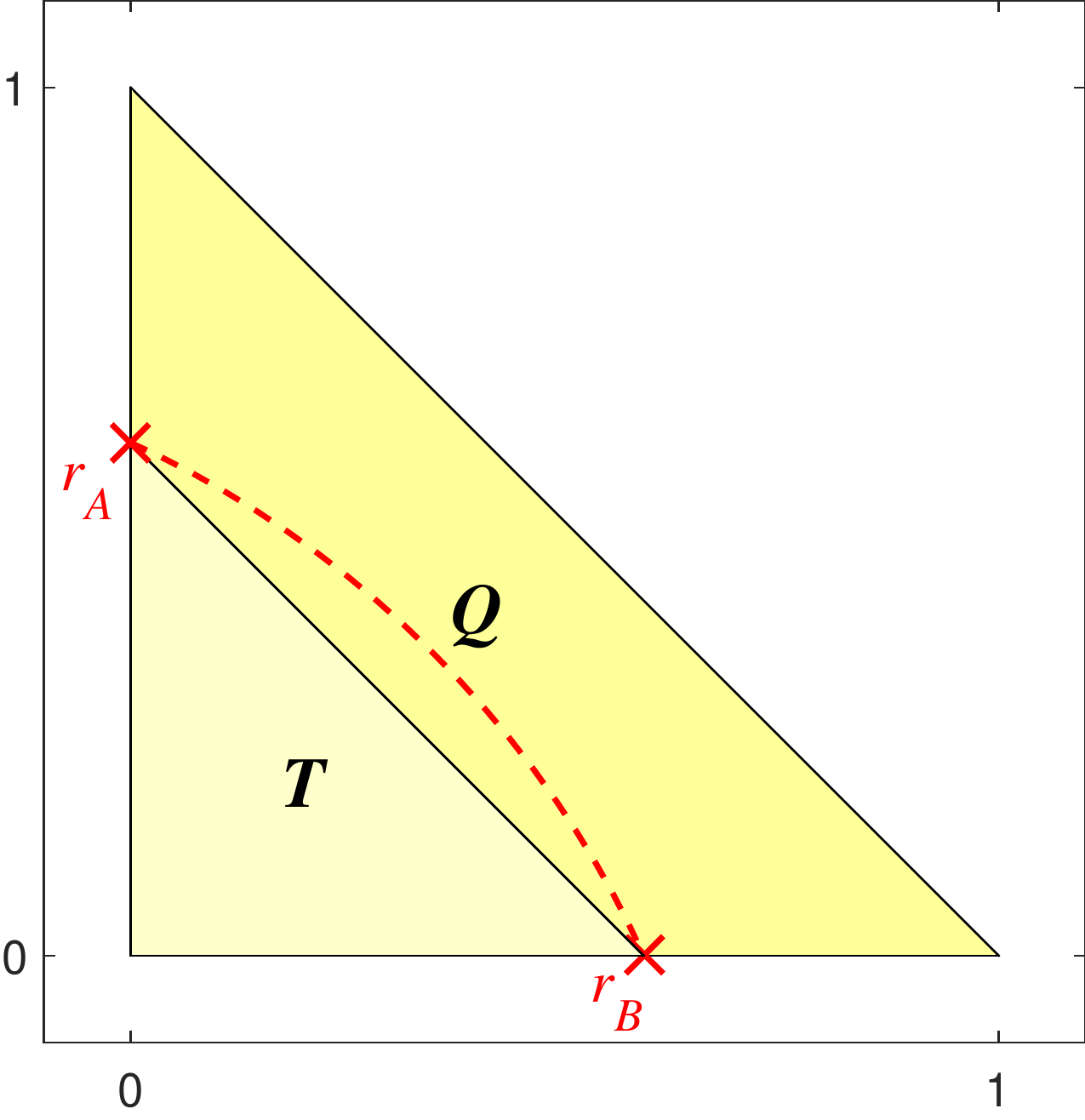}}$\qquad$\subfigure[$\mathbf{T}^{1}_{2}$]{\includegraphics[width=5cm]{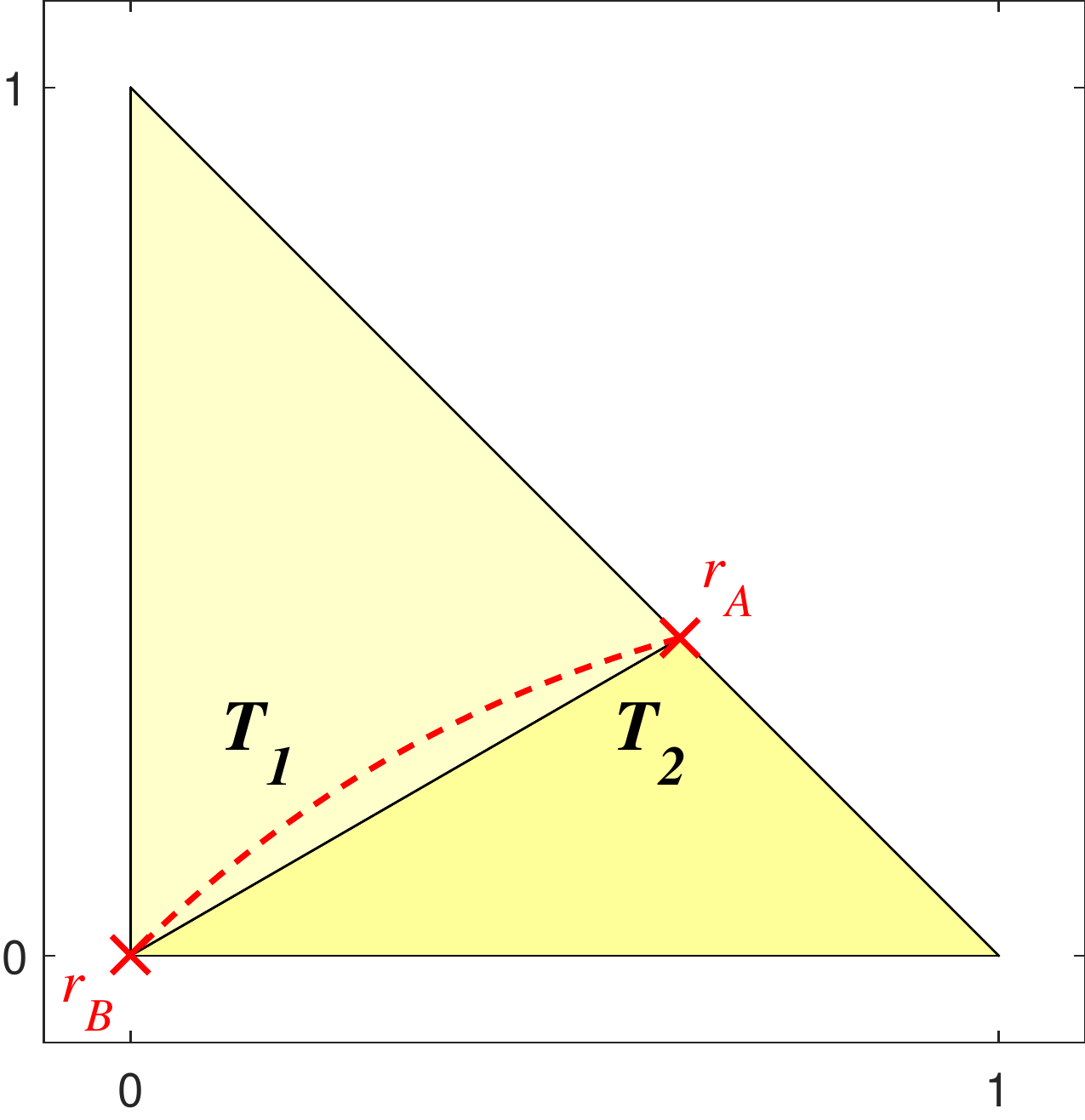}}

\caption{\label{fig:SubDivTrLocal}Local subdivisions for triangles including
the zero level set $\phi^{h}$ (a) when two different edges are cut
and (b) a node and the opposite edge are cut.}
\end{figure}

The subdivisions, which do not affect the neighbour elements or do
not need any information from the neighbour elements to produce a
unique subdivision, are referred to as local subdivisions. The configuration
in Fig.~\ref{fig:SubDivTrLocal}(a) is the first local topology.
It consists of a triangular element divided into a triangle and a
quadrilateral by two cut edges. The situation in the figure is labelled
as $\mathbf{T}_{13}$, but all situations fulfilling $\mathbf{T}_{IJ}$,
$I\le3,J\le3$ and $I\neq J$ are also included in this topological
situation. The other standard topology occurring for triangular elements
is given in Fig.~\ref{fig:SubDivTrLocal}(b) as $\mathbf{T}_{2}^{1}$,
but again the basic topology covers in total the combinations $\{\mathbf{T}_{2}^{1},\mathbf{T}_{3}^{2},\mathbf{T}_{1}^{3}\}$.
In this case a vertex and the opposite edge node are cut, see Fig.~\ref{fig:SubDivTrLocal}(b).
For this topological case, the element patch is easily subdivided
into two triangular elements. This case is particularly important
when distinguished points occur e.g.~a sudden change in the geometry
(e.g.~corners) or boundary conditions. This is considered in detail
in a follow-up paper.

\subsubsection{Non-local subdivisions}

\begin{figure}
\centering
\subfigure[$\mathbf{T}_{11}$]{\includegraphics[width=4.5cm]{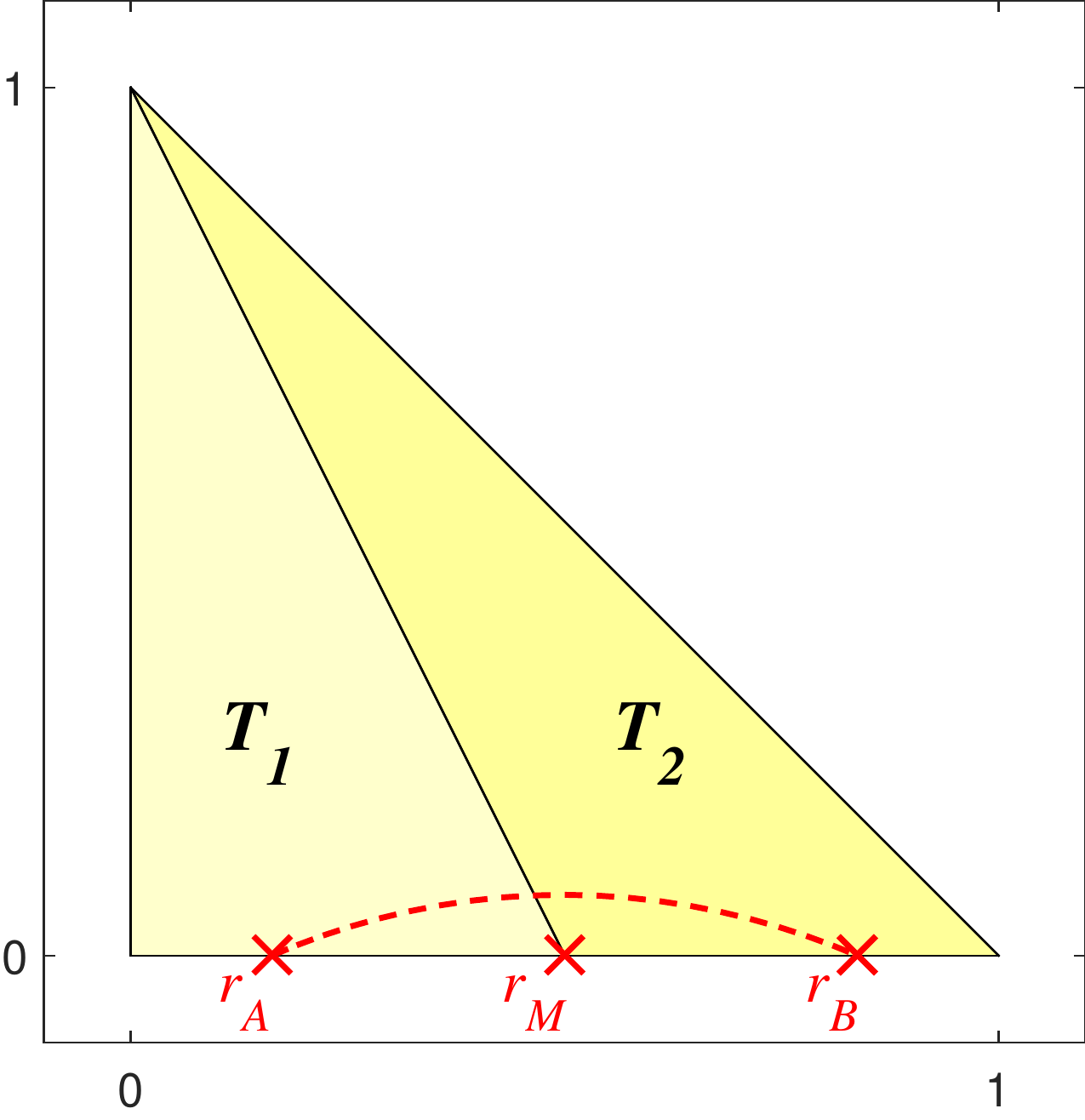}}$\quad$\subfigure[$\mathbf{T}^{1}_{1}$]{\includegraphics[width=4.5cm]{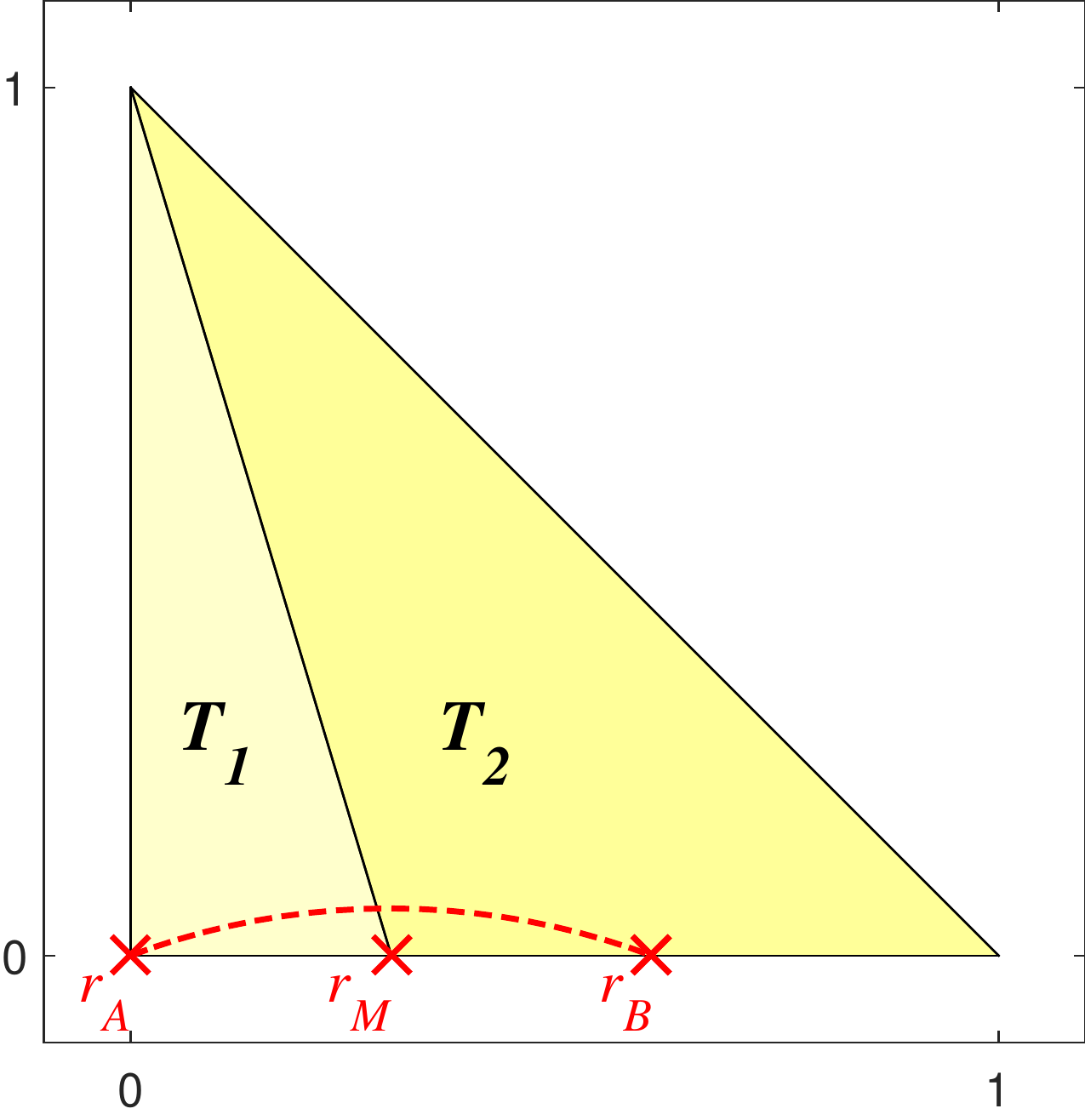}}$\quad$\subfigure[$\mathbf{T}^{12}$]{\includegraphics[width=4.5cm]{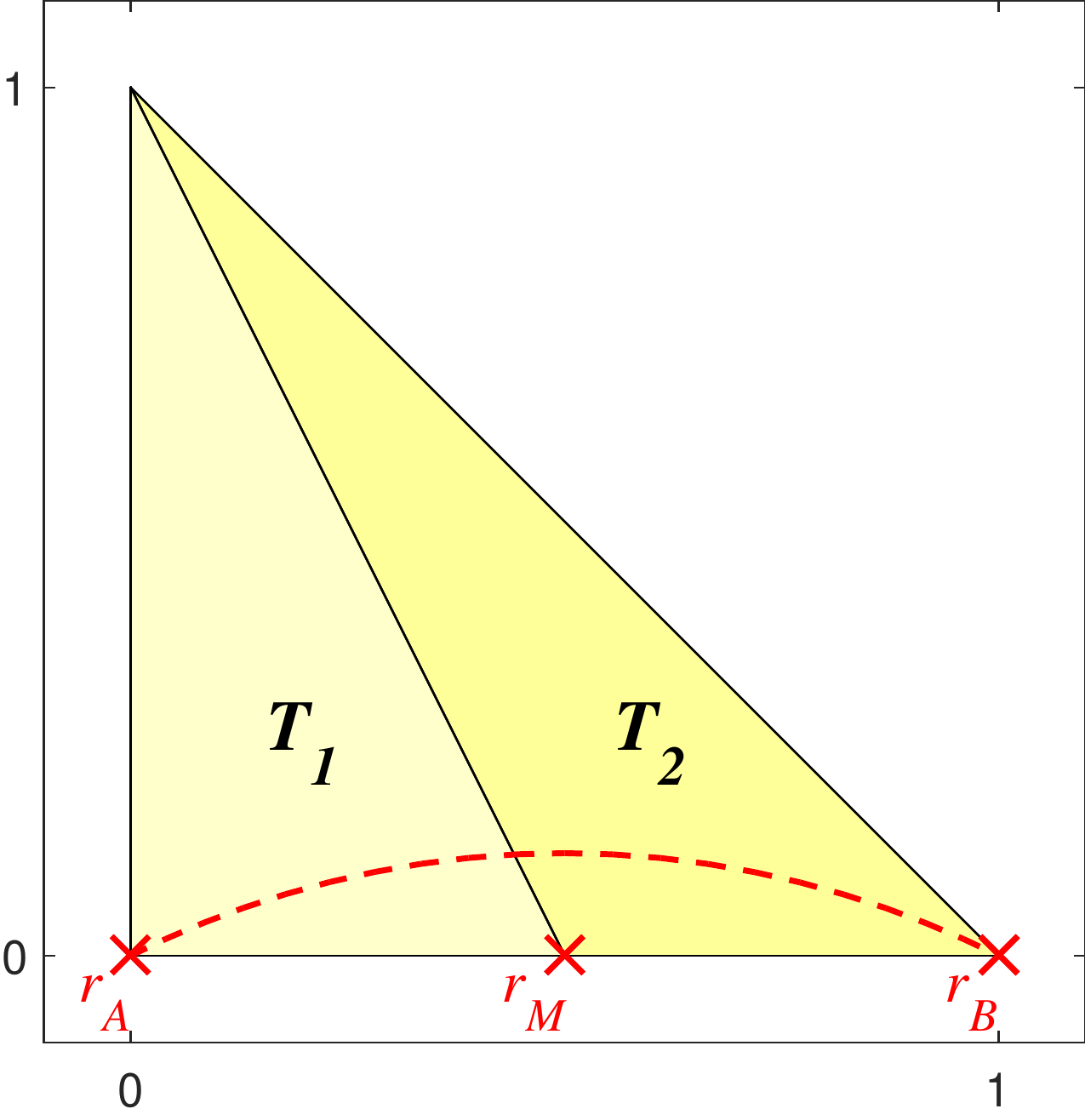}}

\caption{\label{fig:SubDivTrNonLocal}Non-local subdivisions of a triangle
as a consequence of (a) one edge cut twice, (b) one cut node with
a neighbouring cut edge and (c) two cut nodes. }
\end{figure}

Some topologies are considered that need a specialized treatment as
they influence the neighbour elements. These cases occur in connection
with a glancing intersection of the level set with one element edge.
Often the problem could also be resolved by a refinement of the background
mesh. However, from a computational point of view, it is still preferable
to proceed as follows. We go through all elements of the mesh and
find the intersections of the level set function with the element
boundaries. For edges that are cut twice, the adjacent elements are
subdivided as in Fig.~\ref{fig:SubDivTrNonLocal}(a). The element
is subdivided in the middle of the interval defined by the two roots
on the edge, $\vek r_{M}=\frac{1}{2}(\vek r_{A}+\vek r_{B})$. Collectively,
the topological combinations $\{\mathbf{T}_{11},\mathbf{T}_{22},\mathbf{T}_{33}\}$
are covered by this subdivision.

The two other cases are depicted in Fig.~\ref{fig:SubDivTrNonLocal}(b)
where a node is cut together with a neighbouring edge, described by
$\{\mathbf{T}_{1}^{1},\mathbf{T}_{1}^{2},\mathbf{T}_{2}^{2},\mathbf{T}_{2}^{3},\mathbf{T}_{3}^{3},\mathbf{T}_{3}^{1}\}$,
and Fig.~\ref{fig:SubDivTrNonLocal}(c) where two nodes are cut,
described by $\{\mathbf{T}^{11},\mathbf{T}^{22},\mathbf{T}^{33}\}$.
The subdivision strategy is similar as before: Both cases can be subdivided
into two triangles, with the point $\vek r_{M}$ located in the middle
of the intersections of the zero level sets with the element boundary. 

Once the non-local subdivisions are executed, all topological cases
can be treated by the local subdivision strategy described before. 

\subsection{Subdivision strategy for quadrilaterals}

Starting with a mesh composed by quadrilateral elements, one approach
would be to convert the cut elements into triangles and apply the
scheme for triangles from above. However, this would obviously lead
to (slightly) different zero level sets compared to those implied
by the interpolation functions in a quadrilateral element. Therefore,
we rather develop schemes that work without prior conversion. A total
of four \emph{local} subdivisions are introduced, depending on the
cut situation. As before, it is also necessary to consider cases where
the neighbour elements are influenced, which are again named \emph{non-local}
subdivisions.

\subsubsection{Local subdivisions}

\begin{figure}
\centering
\subfigure[$\mathbf{Q}_{14}$]{\includegraphics[width=4.5cm]{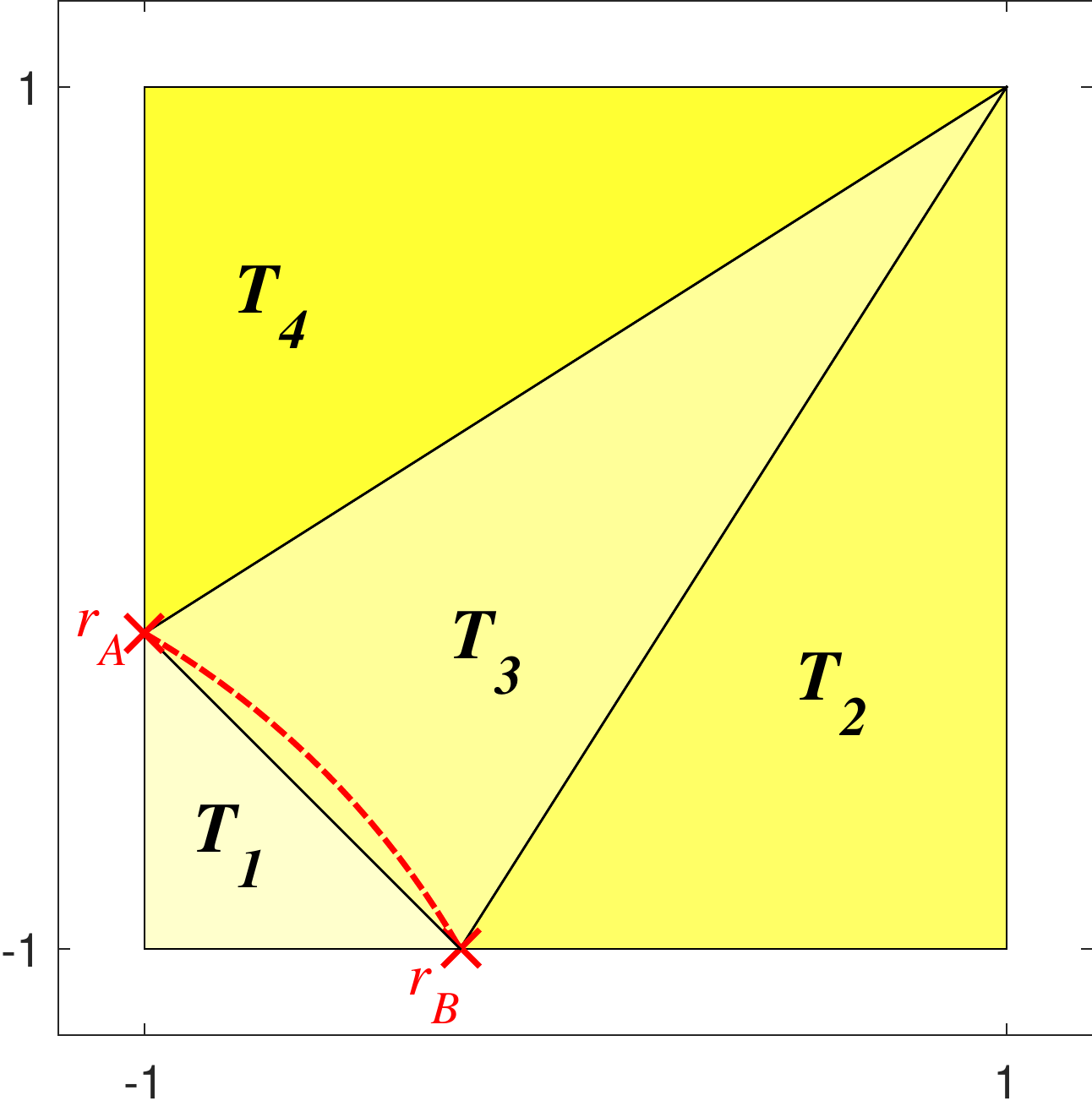}}\qquad\subfigure[$\mathbf{Q}_{24}$]{\includegraphics[width=4.5cm]{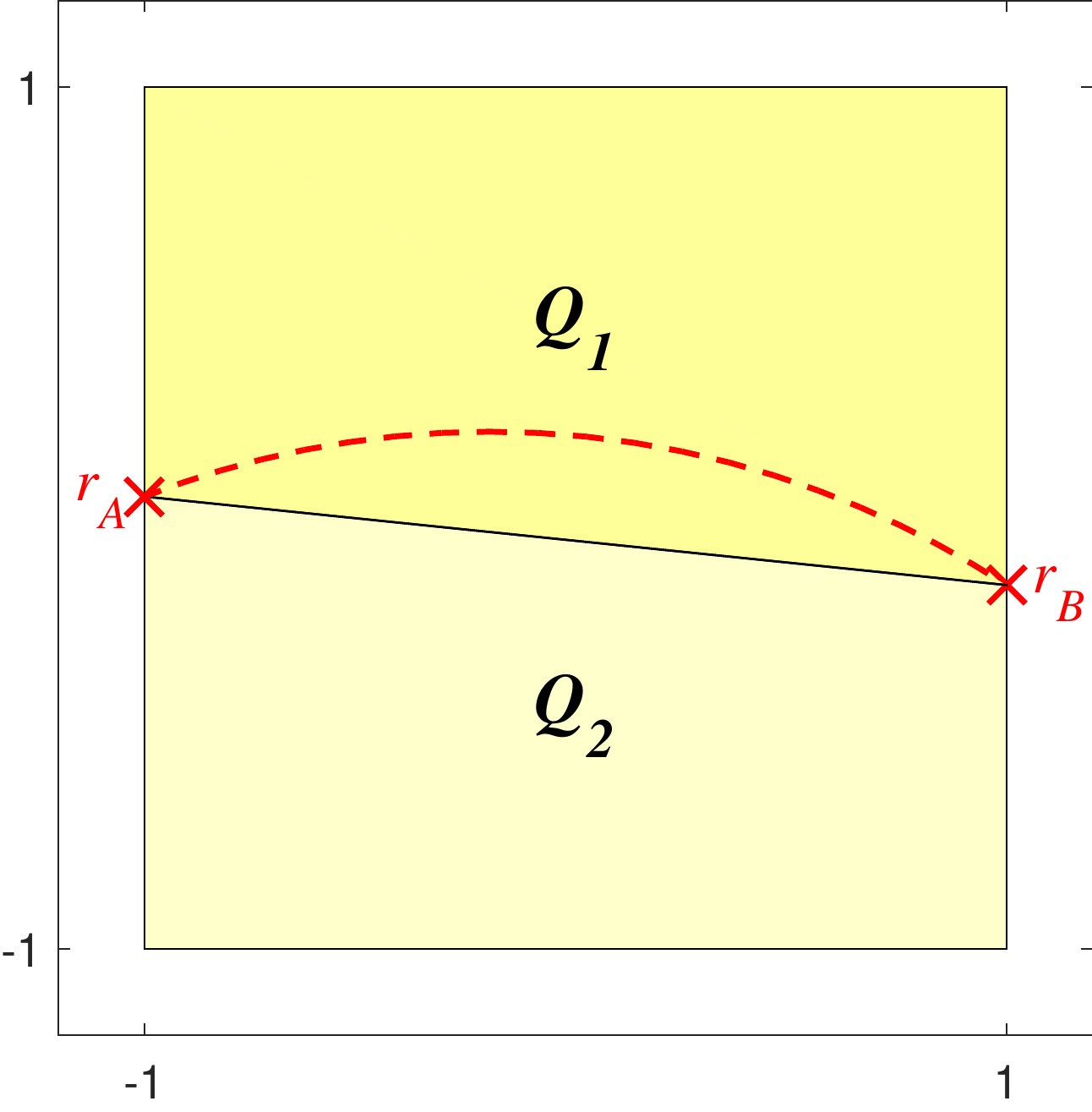}}\\\subfigure[$\mathbf{Q}_{2}^{1}$]{\includegraphics[width=4.5cm]{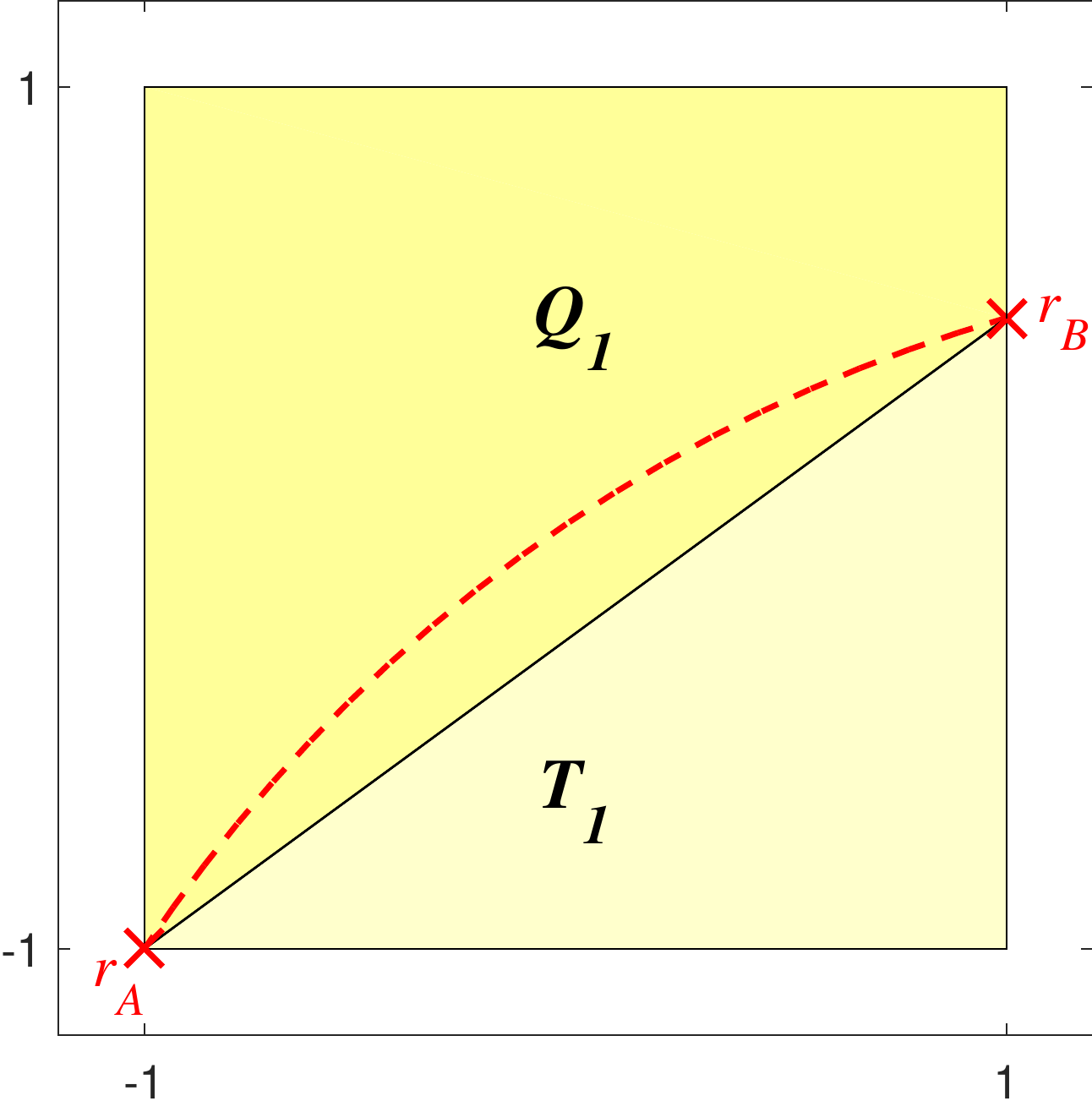}}\qquad\subfigure[$\mathbf{Q}^{13}$]{\includegraphics[width=4.5cm]{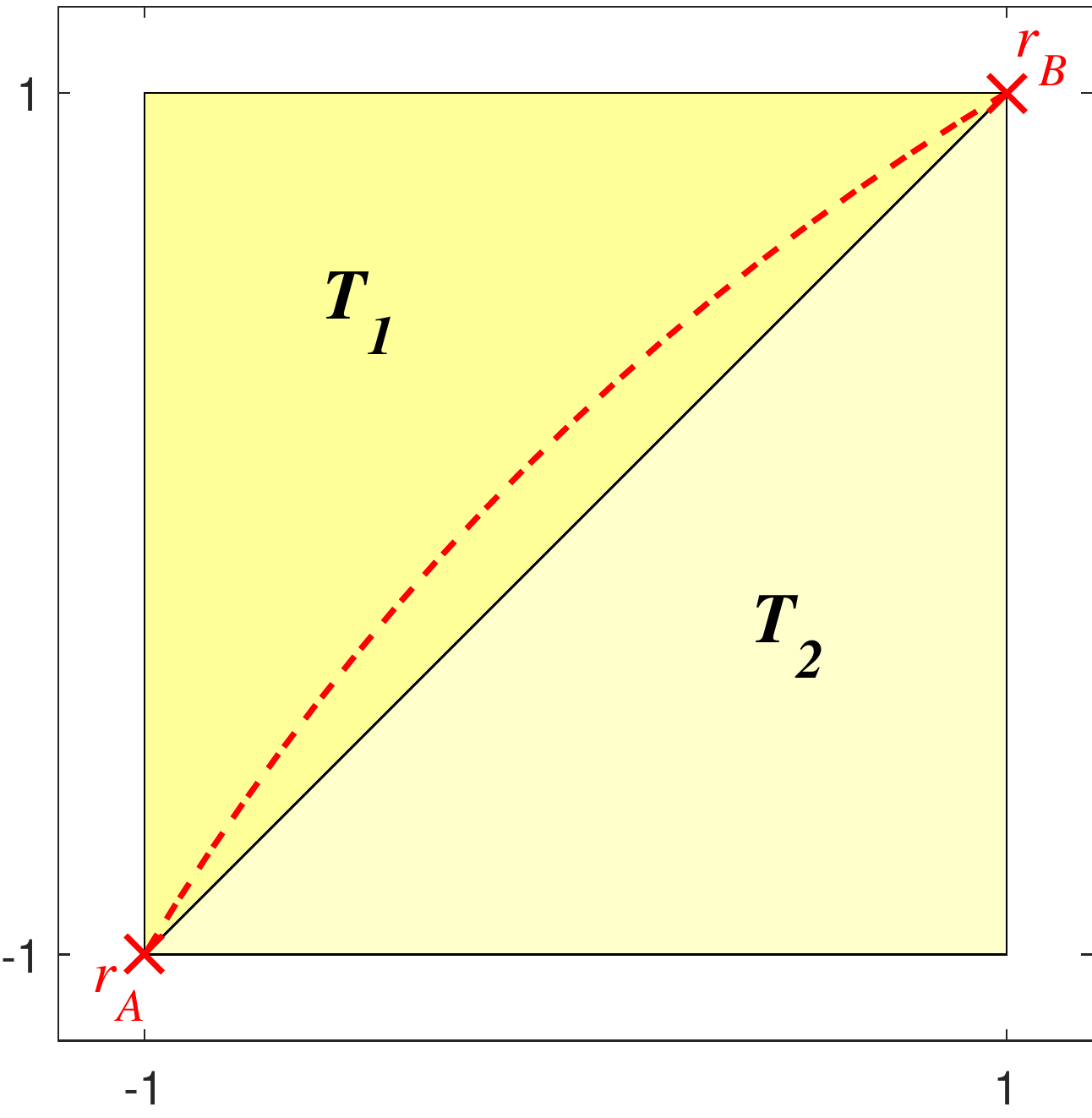}}

\caption{\label{fig:SubDivQuadLocal}All variants of local subdivision schemes
for a quadrilateral elements: (a) features two neighbouring cut edges,
(b) two opposite cute edges, (c) a cut vertex and a non-neighbour
cut edge, and (d) two cut opposite vertices.}
\end{figure}

All situations featuring two cut \emph{adjacent} edges are topologically
equivalent, see Fig.~\ref{fig:SubDivQuadLocal}(a) for $\mathbf{Q}_{14}$.
The subdivision for this topological case produces four triangular
elements. All other topological cases $\{\mathbf{Q}_{12},\mathbf{Q}_{23},\mathbf{Q}_{34}\}$
are also covered with the proposed subdivision. Another case is when
two \emph{neighbouring} edges are cut as depicted in Fig.~\ref{fig:SubDivQuadLocal}(b)
for $\mathbf{Q}_{24}$. This subdivision creates two quadrilaterals.
Again, also the rotated situation $\mathbf{Q}_{13}$, is covered by
the subdivision. The first topology considering also nodal cuts is
depicted in Fig.~\ref{fig:SubDivQuadLocal}(c). In that case, a node
and a non-adjacent edge are cut and subdivide the original quadrilateral
into a quadrilateral and a triangle. The covered situations are $\{\mathbf{Q}_{2}^{1},\mathbf{Q}_{3}^{1},\mathbf{Q}_{3}^{2},\mathbf{Q}_{4}^{2},\mathbf{Q}_{4}^{3},\mathbf{Q}_{1}^{3},\mathbf{Q}_{1}^{4},\mathbf{Q}_{2}^{4}\}$.
The last local subdivision is depicted in Fig.~\ref{fig:SubDivQuadLocal}(c).
This topology features two cut opposite nodes, $\{\mathbf{Q}^{13},\mathbf{Q}^{24}\}$,
leading to two triangles. 

\subsubsection{Non-local subdivisions}

\begin{figure}
\centering
\subfigure[$\mathbf{Q}_{11}$]{\includegraphics[width=4.5cm]{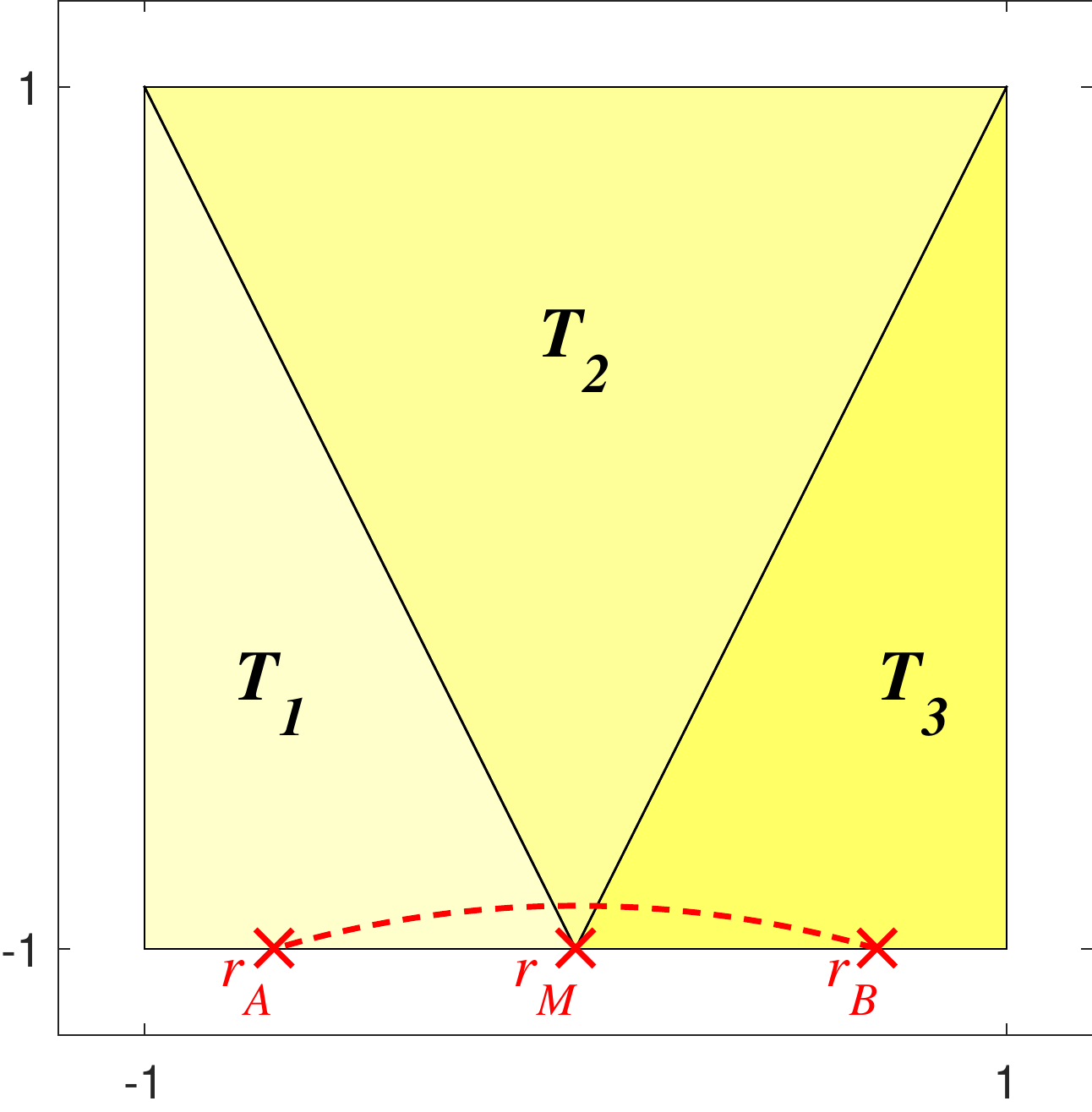}}$\quad$\subfigure[$\mathbf{Q}_{1}^{1}$]{\includegraphics[width=4.5cm]{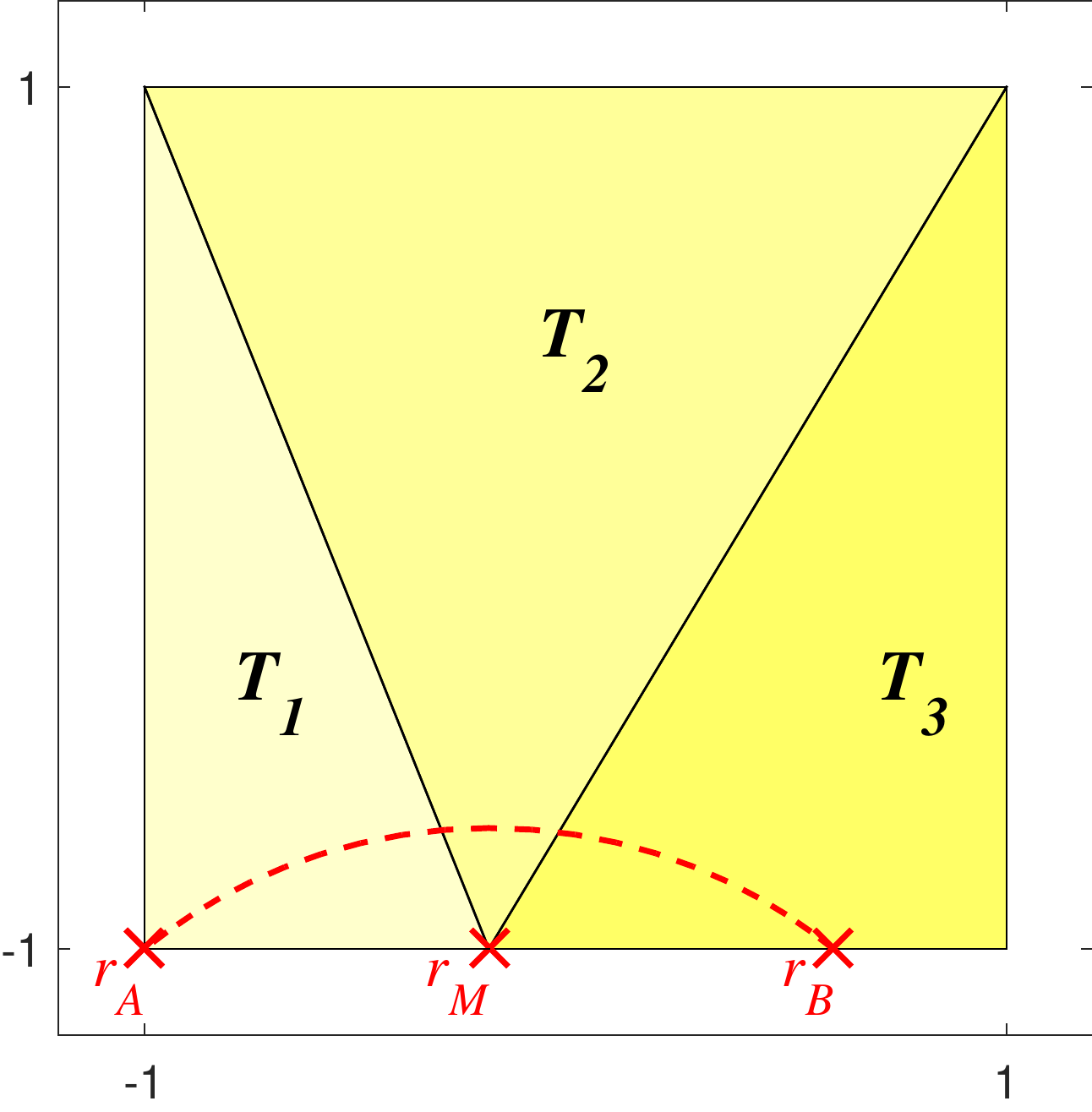}}$\quad$\subfigure[$\mathbf{Q}^{12}$]{\includegraphics[width=4.5cm]{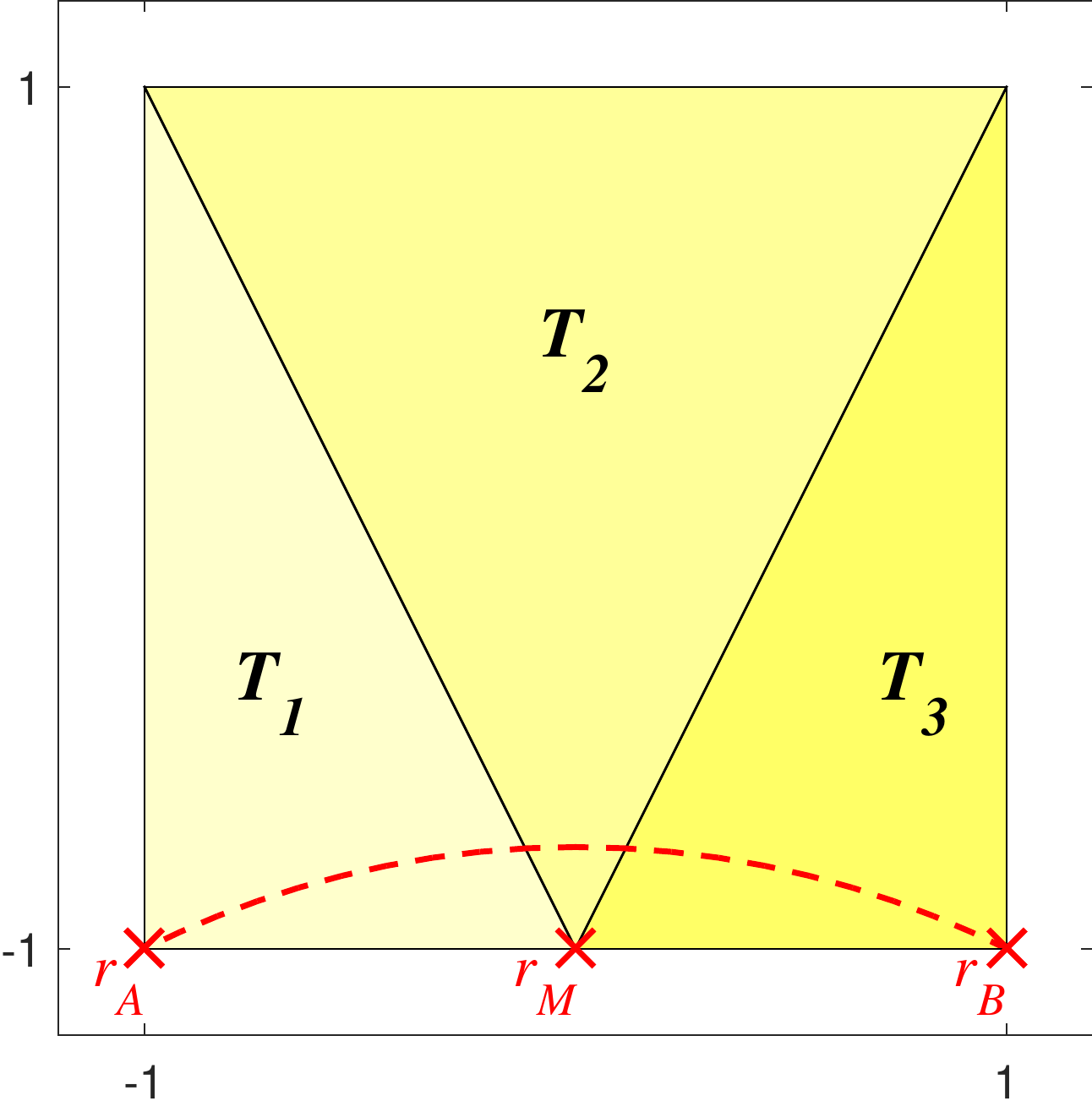}}

\caption{\label{fig:SubDivQuadNonLocal}Non-local subdivisions of a quadrilateral
with (a) one edge cut twice, (b) a cut node next to a cut edge, (c)
two cut neighbouring nodes.}
\end{figure}

To deal with glancing intersections, as depicted in Fig.~\ref{fig:SubDivQuadNonLocal}(a),
it is necessary to subdivide the elements adjacent to the considered
edge. In total, as for the triangular elements, three different topological
cases arise. All situations have in common that one edge is cut twice.
The first case, see Fig.~\ref{fig:SubDivQuadNonLocal}(a), includes
the topologies $\{\mathbf{Q}_{11},\mathbf{Q}_{22},\mathbf{Q}_{33}\}$.
The subdivision of this topological case is done by three triangles,
whereby the triangles share the point $\vek r_{M}=\frac{1}{2}(\vek r_{A}+\vek r_{B})$.
All other topological cases are similar to the former one when placing
one root, as in Fig.~\ref{fig:SubDivQuadNonLocal}(b), or both roots
on the vertex nodes, see Fig.~\ref{fig:SubDivQuadNonLocal}(c). The
covered situations are then $\{\mathbf{Q}_{1}^{1},\mathbf{Q}_{1}^{2},\mathbf{Q}_{2}^{2},\mathbf{Q}_{2}^{3},\mathbf{Q}_{3}^{3},\mathbf{Q}_{3}^{4},\mathbf{Q}_{4}^{4},\mathbf{Q}_{4}^{1}\}$
for one cut vertex and $\{\mathbf{Q}^{11},\mathbf{Q}^{22},\mathbf{Q}^{33}\}$
for two cut vertices. 

\subsection{Coordinate mappings}

The task is now to construct a map for all points from a higher-order
reference element to the previously defined low-order sub-cells with
potentially one higher-order side. The mapping is constructed so that
the following properties are fulfilled: (i) the element resulting
from the coordinate transformation yields optimal interpolation properties
and (ii) one element edge aligns with the (curved) higher-order interface
element $\Gamma_{\vek r}^{h}$.

\begin{figure}
\centering\includegraphics[scale=0.83]{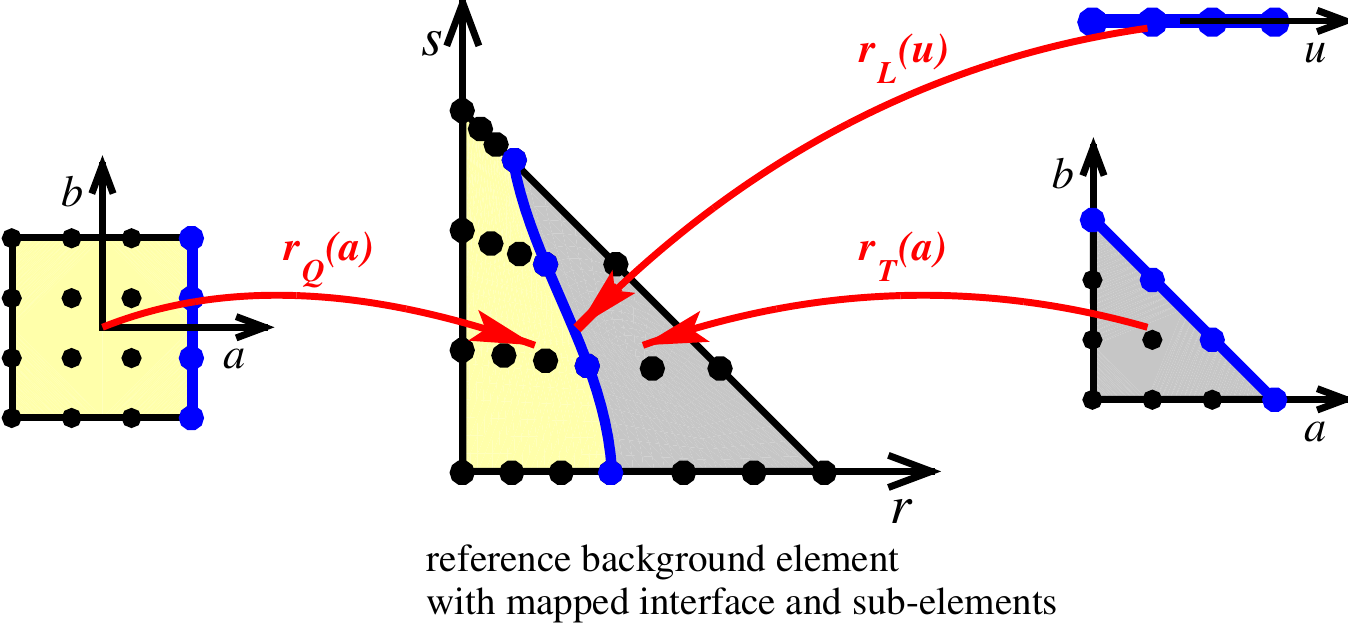}

\caption{\label{fig:Chap4MapBlendOverview}A higher-order background element
is cut by the reconstructed interface element $\Gamma_{\vek r}^{h}$
and subdivided into low-order sub-cells defining topology $\mathbf{T}_{12}$.
The mapping from the reconstructed element to the coordinates $\vek r$
is given by $\vek r_{L}(u)$. A set of higher-order reference elements
in the coordinate system $\vek a$ is introduced. The coordinate transformations
$\vek r_{T}(\vek a)$ and $\vek r_{Q}(\vek a)$ map the higher-order
reference sub-elements to the coordinates $\vek r$ such that one
special side (in our case always edge 2) matches exactly $\Gamma_{\vek r}^{h}$
.}
\end{figure}

As seen from Fig.~\ref{fig:Chap4MapBlendOverview}, several element
types and corresponding coordinate systems are involved: Central to
the whole procedure are the low-order sub-cells (with one higher-order
side) described in the coordinate system $\vek r$, as defined above.
They define the outer boundary of the higher-order elements resulting
from the automatic conformal decomposition. Nevertheless, the mapping
for inner points in the sub-cells is not unique, so that inner nodes
of the resulting higher-order elements may be placed differently which
has important implications on the convergence properties. The higher-order
interface element representing the discretized zero level set $\phi^{h}$
is initially defined in the reference coordinate system $u$ and is
mapped to the reference element $\Omega_{\vek r}$ by $\vek r(u)$.
Finally, there are higher-order reference sub-elements defined in
the coordinate system $\vek a$. These sub-elements are used to define
the initial location of the nodes and the mapping to $\Omega_{\vek r}$
is denoted as $\vek r(\vek a)$.

It is useful to introduce the following definitions related to nodal
sets and corresponding shape functions of the higher-order sub-elements
in the reference coordinate system $\vek a$. The set of the corner
nodes for each element is defined by $I_{\mathrm{lin}}^{\Omega}$.
These are three nodes for a triangular and four for a quadrilateral
element. From these nodes, either one node (see Fig.~\ref{fig:SubDivQuadLocal}(a)
topology $\mathbf{Q}_{14}$, elements $\mathbf{T}_{2}$ and $\mathbf{T}_{4}$)
or two nodes (see Fig.~\ref{fig:SubDivQuadLocal}(a) for topology
$\mathbf{Q}_{14}$, elements $\mathbf{T}_{1}$ and $\mathbf{T}_{3}$)
belonging to $I_{\mathrm{side}}^{\Omega}\subset I_{\mathrm{lin}}^{\Omega}$
are on the linear reconstruction when mapped to $\vek r$. For simplicity,
we choose a node numbering in the sub-elements such that edge 2 is
the one that after the mapping to $\vek r$ will coincide exactly
with the interface element. The set of nodes on edge 2 is called $I_{\mathrm{ho}}^{\Gamma}$,
and the corner nodes are $I_{\mathrm{lin}}^{\Gamma}\subset I_{\mathrm{ho}}^{\Gamma}$.
Related to these nodes are the nodal coordinates $\vek a_{i}$ before
the mapping and $\vek r_{i}$ after the mapping to the cut reference
background element. Note that the nodes in $I_{\mathrm{side}}^{\Omega}$
and $I_{\mathrm{lin}}^{\Gamma}$ superpose in coordinate systems $\vek a$
and $\vek r$. We associate the following set of shape functions with
the above-mentioned nodes: $N_{\mathrm{lin},i}^{\Omega}\left(\vek a\right)$
for $i\in I_{\mathrm{lin}}^{\Omega}$ are the low-order shape functions
of the reference sub-element. $N_{\mathrm{ho},i}^{\Gamma}\left(u\right)$
for $i\in I_{\mathrm{ho}}^{\Gamma}$ are the higher-order shape functions
on the special side of the sub-element. $N_{\mathrm{lin},i}^{\Gamma}\left(u\right)$
for $i\in I_{\mathrm{lin}}^{\Gamma}$ are low-order shape functions
on the special side which coincide with the trace of the low-order
shape functions $N_{\mathrm{lin},i}^{\Omega}\left(\vek a\right)$
for $i\in I_{\mathrm{side}}^{\Omega}$ on edge 2.

The aim is to construct a map for all nodal points from $\vek a$
to $\vek r$ based on the known parametrization of the edges (note
that $\Gamma_{\vek r}^{h}$ is included in some 2D-elements as a curved
edge). The mapping transforms the low-order sub-cells into higher-order
elements with one curved edge. Therefore, we adapt a transfinite interpolation
technique as described e.g.~in \cite{Gordon_1973a,Gordon_1973b,Szabo_2004a,Solin_2003a}.
The scheme yields the blending function method for the quadrilaterals
\cite{Gordon_1973a,Gordon_1973b} and a slightly modified technique
for the triangular elements. We start with the mapping for the elements
having one node in $I_{\mathrm{side}}^{\Omega}$ (or alternatively
the assumption that $\Gamma_{\vek r}^{h}$ is straight). In this case
the mapping from $\vek a$ to $\vek r$ is defined completely by the
vertices and is given by the standard isoparametric mapping for bi-linear
elements as
\begin{equation}
\vek r_{m}\left(\vek a_{n}\right)=\sum_{i\in I_{\mathrm{lin}}^{\Omega}}N_{\mathrm{lin},i}^{\Omega}\left(\vek a_{n}\right)\cdot\vek r_{\mathrm{lin},i}^{\Omega},\quad\vek a_{n},\vek r_{m}\in\mathbb{R}^{2}.\label{eq:IsoParamMapping}
\end{equation}
The coordinates $\vek a_{n}$ and $\vek r_{m}$ are the nodal coordinates
and $\vek r_{\mathrm{lin},i}^{\Omega}$ are the corner nodes of a
sub-element in the cut reference background element in the corresponding
coordinate systems. To include the more elaborate case of a \emph{curved}
higher-order interface element, an additional term to the standard
isoparametric mapping must be included, defined as
\begin{equation}
\vek r_{m}\left(\vek a_{n}\right)=\sum_{i\in I_{\mathrm{lin}}^{\Omega}}N_{\mathrm{lin},i}^{\Omega}\left(\vek a_{n}\right)\cdot\vek r_{\mathrm{lin},i}^{\Omega}+\vek M\left(\vek a_{n}\right),\quad\vek a_{n},\vek r_{m}\in\mathbb{R}^{2}.\label{eq:BldgMethod2}
\end{equation}
Equations (\ref{eq:IsoParamMapping}) and (\ref{eq:BldgMethod2})
are valid for \emph{any} point inside $\Omega_{\vek a}$. Note that
Eq.~(\ref{eq:BldgMethod2}) \emph{must} include Eq.~(\ref{eq:IsoParamMapping})
as a special case of an element with only straight edges. Therefore
the additional vector function $\vek M\left(\vek a\right)\in\mathbb{R}^{d}$
must be defined such that it vanishes for elements with only straight
edges. It is interesting to note that the definition of the geometric
mapping (\ref{eq:BldgMethod2}) formally matches \emph{enriched} approximations
such as those used in the XFEM \cite{Belytschko_1999a,Moes_1999a,Fries_2009b,Babuska_2002a,Babuska_2004a}.

For the definition of $\vek M\left(\vek a\right)$, we impose two
important constraints:
\begin{enumerate}
\item The corner nodes of the higher-order interface element $I_{\mathrm{lin}}^{\Gamma}$
match the corner nodes of the special side of the low-order element
$I_{\mathrm{side}}^{\Omega}$ in the coordinate system $\vek r$. 
\item Those edges of the low-order elements that are not on the special
side shall remain straight after the mapping.
\end{enumerate}
This leads to the requirement that
\begin{equation}
\vek M\left(\vek a_{i}\right)=\vek0\quad\forall i\in I_{\mathrm{lin}}^{\Omega}.\label{eq:BldgMethod3}
\end{equation}
 We find that this can be achieved when splitting $\vek M\left(\vek a\right)$
into the product
\[
\vek M\left(\vek a\right)=\psi\left(\vek a\right)\cdot\vek f\left(u\left(\vek a\right)\right)
\]
 with the definition of $\vek f\left(u(\vek a)\right)$ as 
\begin{equation}
\vek f\left(u(\vek a)\right)=\sum_{i\in I_{\mathrm{ho}}^{\Gamma}}N_{\mathrm{ho},i}^{\Gamma}(u(\vek a))\cdot\vek r_{i}^{\Gamma}-\sum_{i\in I_{\mathrm{lin}}^{\Gamma}}N_{\mathrm{lin},i}^{\Gamma}(u(\vek a))\cdot\vek r_{i}^{\Gamma}\label{eq:BldgMethod4}
\end{equation}
The coordinates $\vek r_{i}^{\Gamma}$ are thereby the element nodes
of the reconstructed interface elements in the coordinate system $\vek r$.
This function corresponds geometrically to subtracting the linear
connection of the nodes in $I_{\mathrm{lin}}^{\Gamma}$ from the curved
higher-order side. The function $\psi$ is differently defined for
triangular and quadrilateral elements and describes the procedure
how the curved edge is blended into an element with only straight
edges. However, note that $\psi$ is not unique and there are many
conceivable alternatives to it. 

As theoretically devised by \cite{Ciarlet_1972b} for 2D simplex elements
and extended by \cite{Lenoir_1986a} and \cite{Bernardi_1989a} for
higher dimensions, certain smoothness properties have to be fulfilled
by the desired coordinate transformation. In general, it can be stated
that the geometrical mapping that is defined by the shape functions
in the reference domain should also be polynomial in the real space.
Further explanations and extensions to transfinite mappings are found
in \cite{Solin_2003a},\cite{Trangenstein_2013a}.

\subsubsection{Mapping for curved triangular elements}

\begin{figure}
\centering
\subfigure[]{\includegraphics[width=5cm]{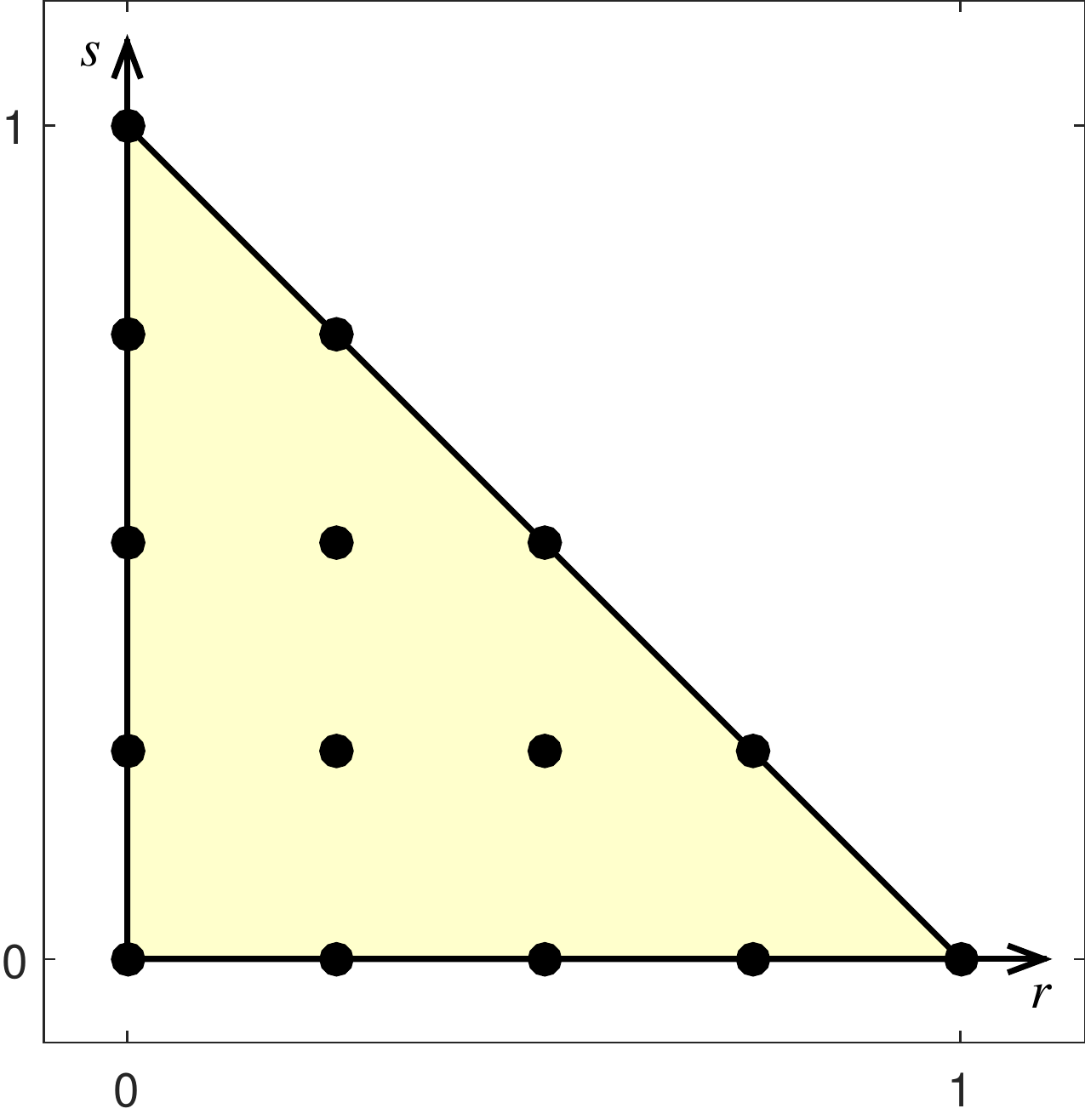}}$\qquad$\subfigure[]{\includegraphics[width=5cm]{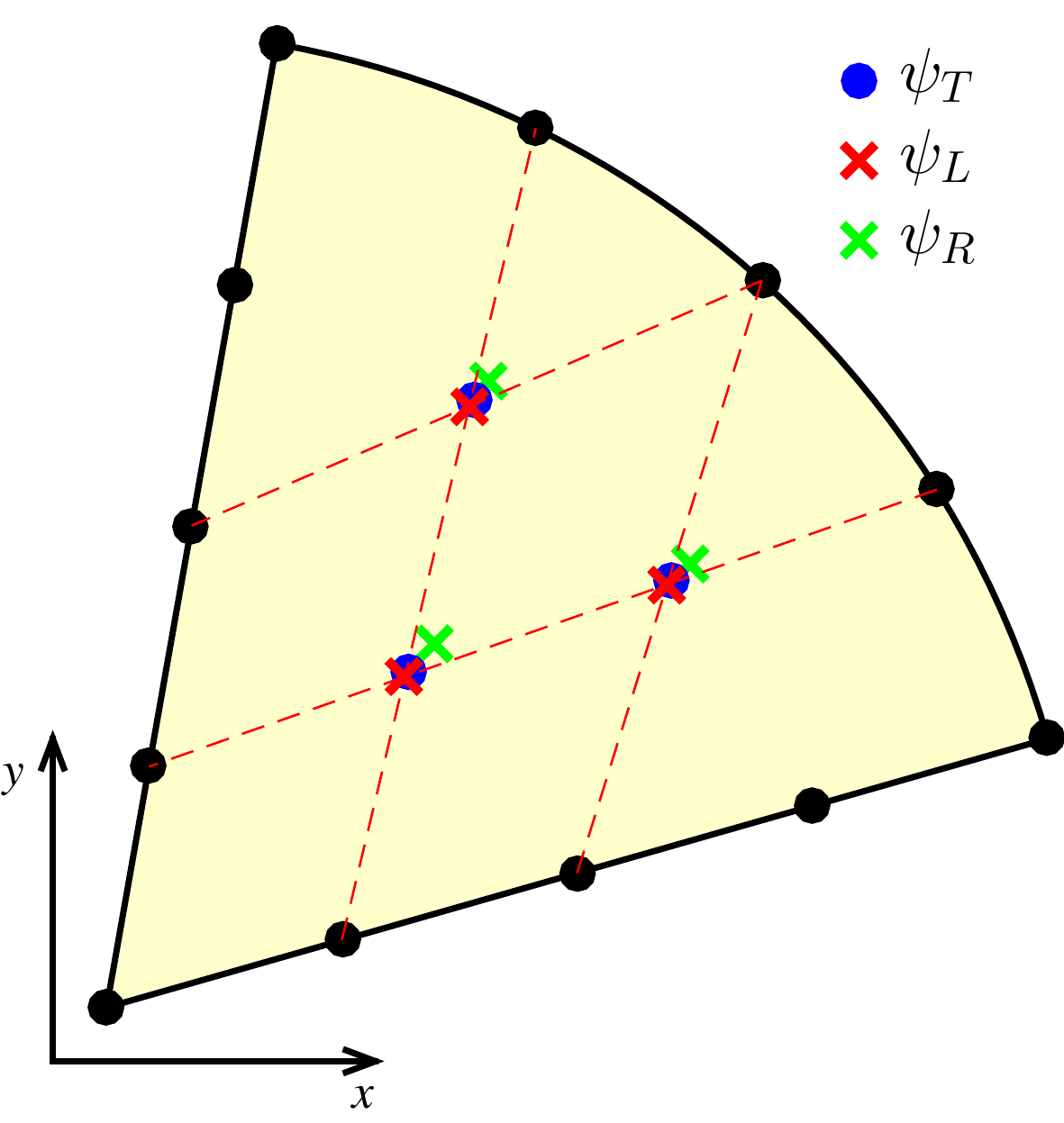}}

\caption{\label{fig:MappingTri}(a) Node distribution inside a reference triangle
for order $m=4$; (b) location of internal nodes for different choices
of $\psi$ in case of a circular arc parametrization of edge 2.}
\end{figure}

It remains to define $u=u_{T}(\vek a)$ and $\psi$. The parametrization
of edge 2 is done by means of vertex functions and gives $u_{T}(\vek a)=N_{\mathrm{lin},3}^{\Omega}\left(\vek a\right)-N_{\mathrm{lin},2}^{\Omega}\left(\vek a\right)$.
The choice of $\psi$ is then given by \cite{Solin_2003a} as
\begin{equation}
\mbox{\ensuremath{\psi}=\ensuremath{\psi_{T}} :}=\frac{\prod\limits _{i\in I_{\mathrm{lin}}^{\Omega}}N_{\mathrm{lin},i}^{\Omega}\left(\vek a\right)}{\prod\limits _{i\in I_{\mathrm{lin}}^{\Gamma}}N_{\mathrm{lin},i}^{\Gamma}\left(u_{T}\left(\vek a\right)\right)}.\label{eq:BldgMethod5-1-1}
\end{equation}
The mapping is valid for all points with exception of nodal points
on $I_{\mathrm{side}}^{\Omega}$ where a division by zero arises.
For these points, however, the mapping of the vertices is defined
by the standard isoparametric mapping. Anyway, also alternative choices
of $\psi$ are conceivable, e.g.~using the ramp function (which is
in fact the mapping of choice for the quadrilaterals), hence $\psi_{R}:=R\left(\vek a\right)=\sum_{I_{\mathrm{side}}^{\Omega}}N_{\mathrm{lin},i}^{\Omega}\left(\vek a\right)$.
Another choice is to use mappings such as in \cite{Ciarlet_1972b,Scott_1973a,Lenoir_1986a}.
In this case, intersections of the boundary nodes are taken as internal
nodes, this is further on designated as $\psi_{L}$ (\emph{Lenoir-mapping}). An issue with
this approach is that the extension to 3D is not trivial. Differences
in the mappings are plotted in Fig.~\ref{fig:MappingTri}. In the
following chapter, we investigate the influence of the choices of
$\psi$ on the convergence rates. Note that in the case of straight
edges, the enrichment function $\vek M_{T}\left(\vek a\right)=\vek0$,
and consequently a standard bi-linear mapping is obtained.

In Fig.~\ref{fig:MappingTri} the inner node distributions are shown
for different versions of $\psi$. The edge is hereby parametrized
by equally spaced nodes on a quarter circle. Note that the mapping
for the boundary nodes coincides for all choices of $\psi$ but differs
inside. As shall be seen below, the node distribution inside the elements
is of large importance for the achieved convergence properties.

\subsubsection{Mapping for curved quadrilateral elements }

For quadrilateral elements, the optimal mapping is found using the
blending function method. Introduced by \cite{Gordon_1973b} for transfinite
(in general non-polynomial) mappings, it can also be used for isoparametric
edge approximations. In this case, a scalar vertex or ramp function
is defined by 
\begin{equation}
\psi:=R\left(\vek a\right)=\sum_{I_{\mathrm{side}}^{\Omega}}N_{\mathrm{lin},i}^{\Omega}\left(\vek a\right).\label{eq:BldgMethod5-1}
\end{equation}
The enrichment function $\vek M\left(\vek a\right)$ for quadrilaterals
is then defined as the linearly blended map
\[
\vek M_{Q}\left(\vek a\right)=R\left(\vek a\right)\cdot\vek f\left(u_{Q}\left(\vek a\right)\right).
\]

For quadrilateral elements the parametrization of edge 2 is given
by $u_{Q}(\vek a)=b$. It is seen that $\vek f\left(u_{Q,i}\right)=\vek0$
for $i\in I_{\mathrm{lin}}^{\Gamma}$ and, hence, also for the nodes
in $I_{\mathrm{side}}^{\Omega}$. Furthermore, $R\left(\vek a_{i}\right)=0$
for $i\in I_{\mathrm{lin}}^{\Omega}\setminus I_{\mathrm{side}}^{\Omega}$,
i.e.~for all remaining nodes in $I_{\mathrm{lin}}^{\Omega}$. Again,
if edge 2 is straight, the enrichment function $\vek M_{Q}\left(\vek a\right)=\vek0$,
and a standard low-order isoparametric mapping is obtained.

\begin{figure}
\centering
\subfigure[]{\includegraphics[width=5cm]{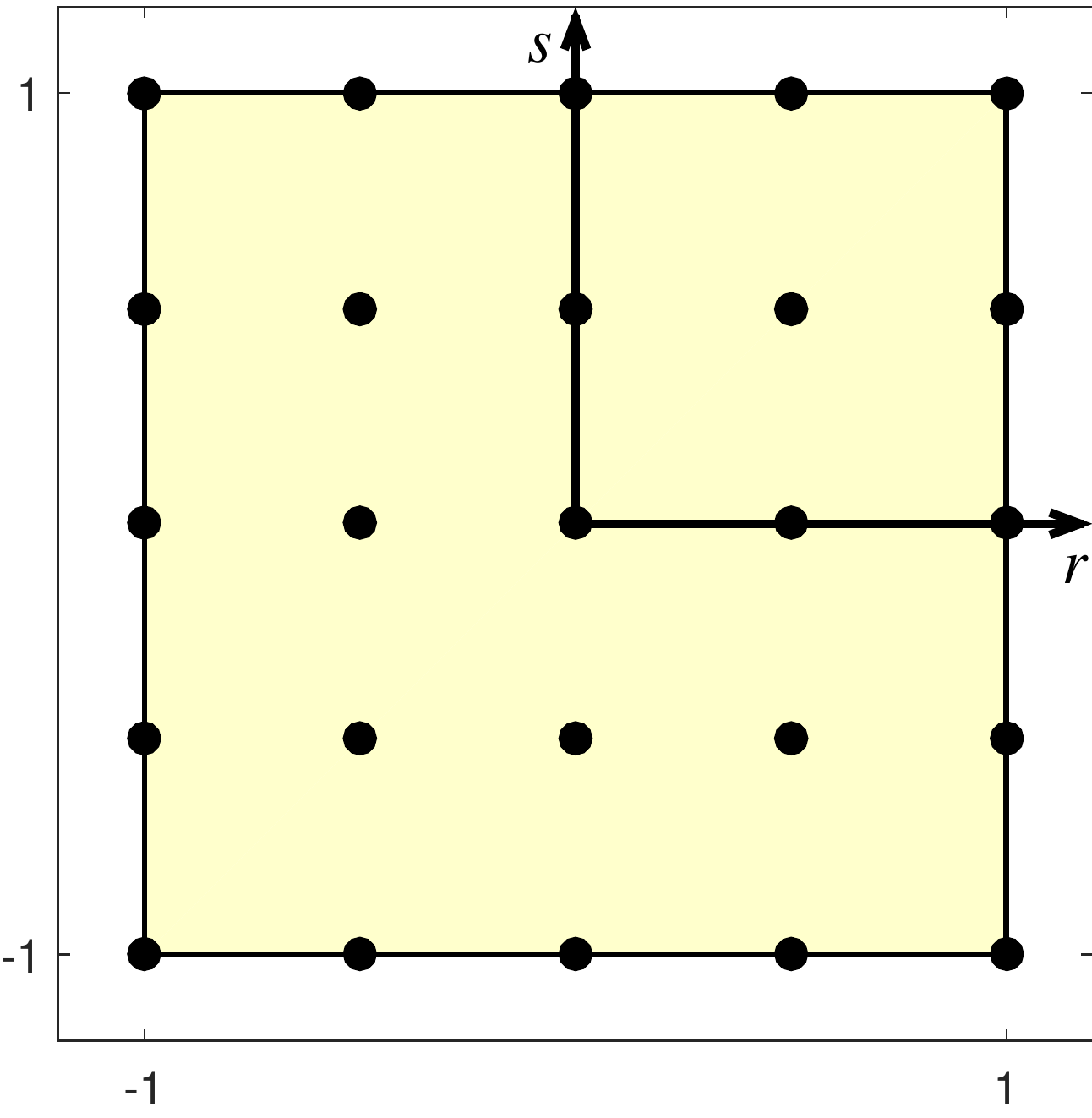}}$\qquad$\subfigure[]{\includegraphics[width=5cm]{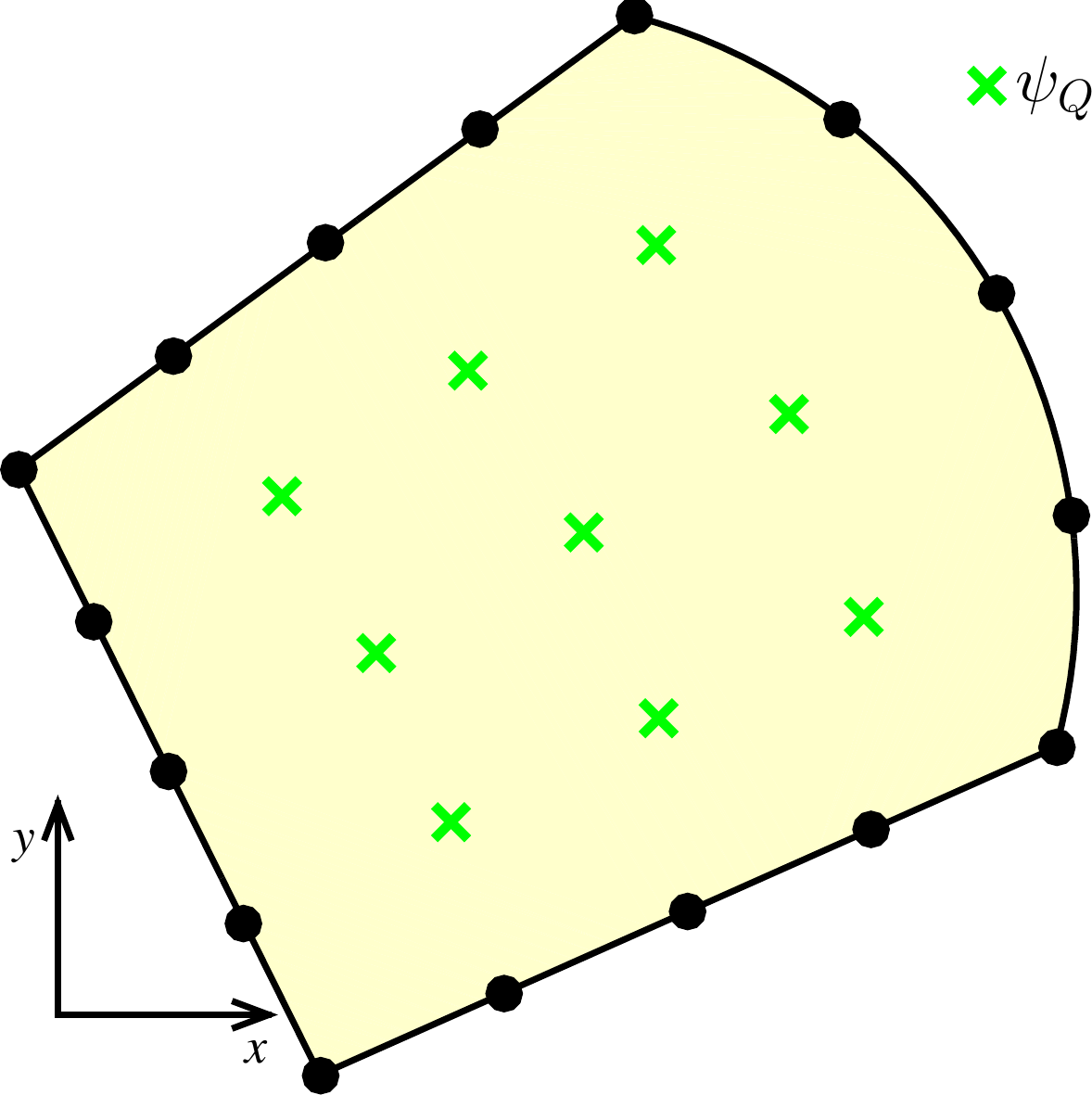}}

\caption{\label{fig:MappingQuad}(a) Node distribution for a quadrilateral
element of order $m=4$, and (b) nodal distribution using the blending
function method for a circular arc parametrization of edge 2 and equidistant
node distribution on the arc.}
\end{figure}

\subsection{Invalid topologies}

\begin{figure}
\centering
\subfigure[]{\includegraphics[width=3.5cm]{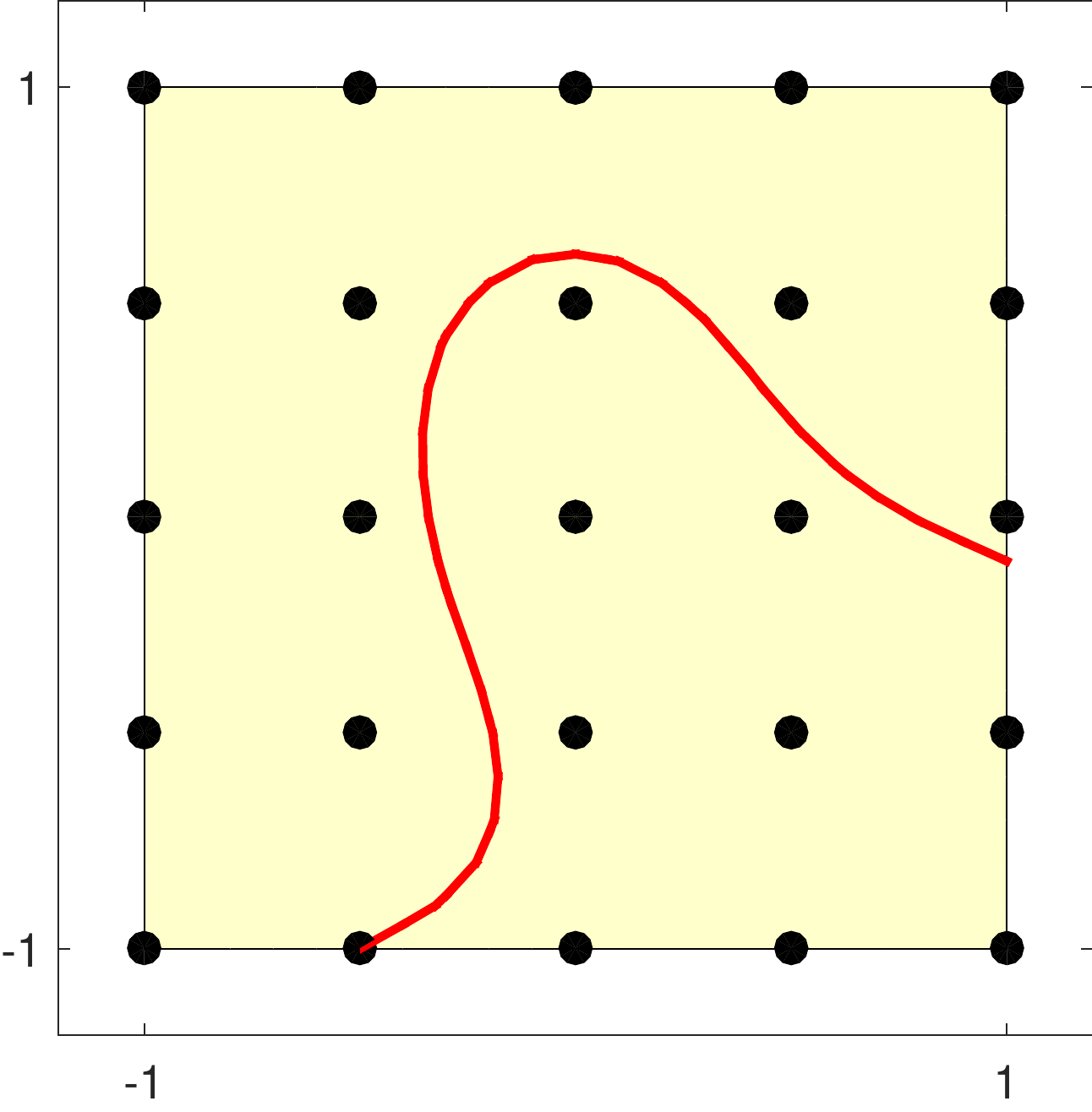}}\subfigure[]{\includegraphics[width=3.5cm]{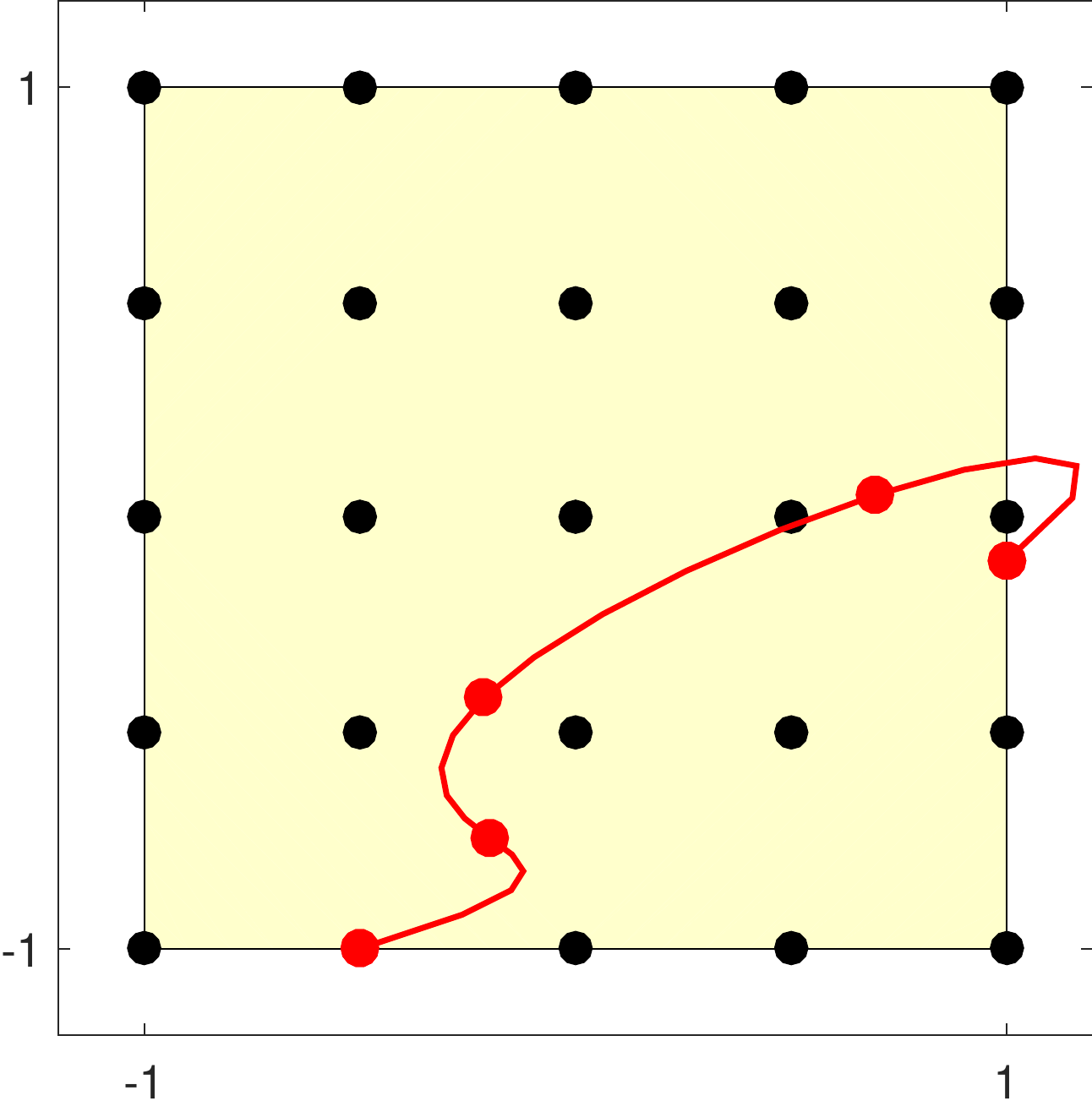}}\subfigure[]{\includegraphics[width=3.5cm]{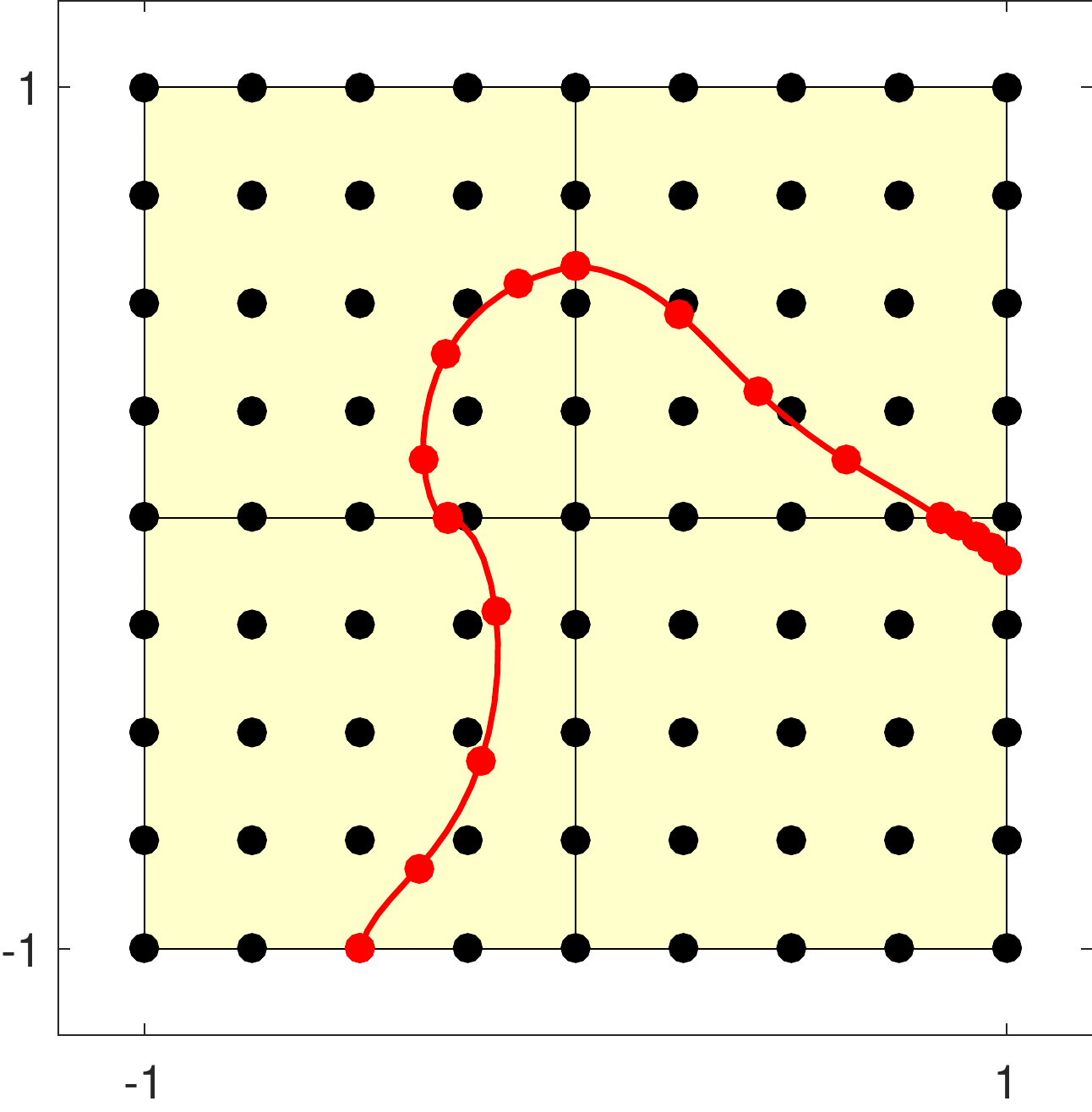}}\subfigure[]{\includegraphics[width=3.5cm]{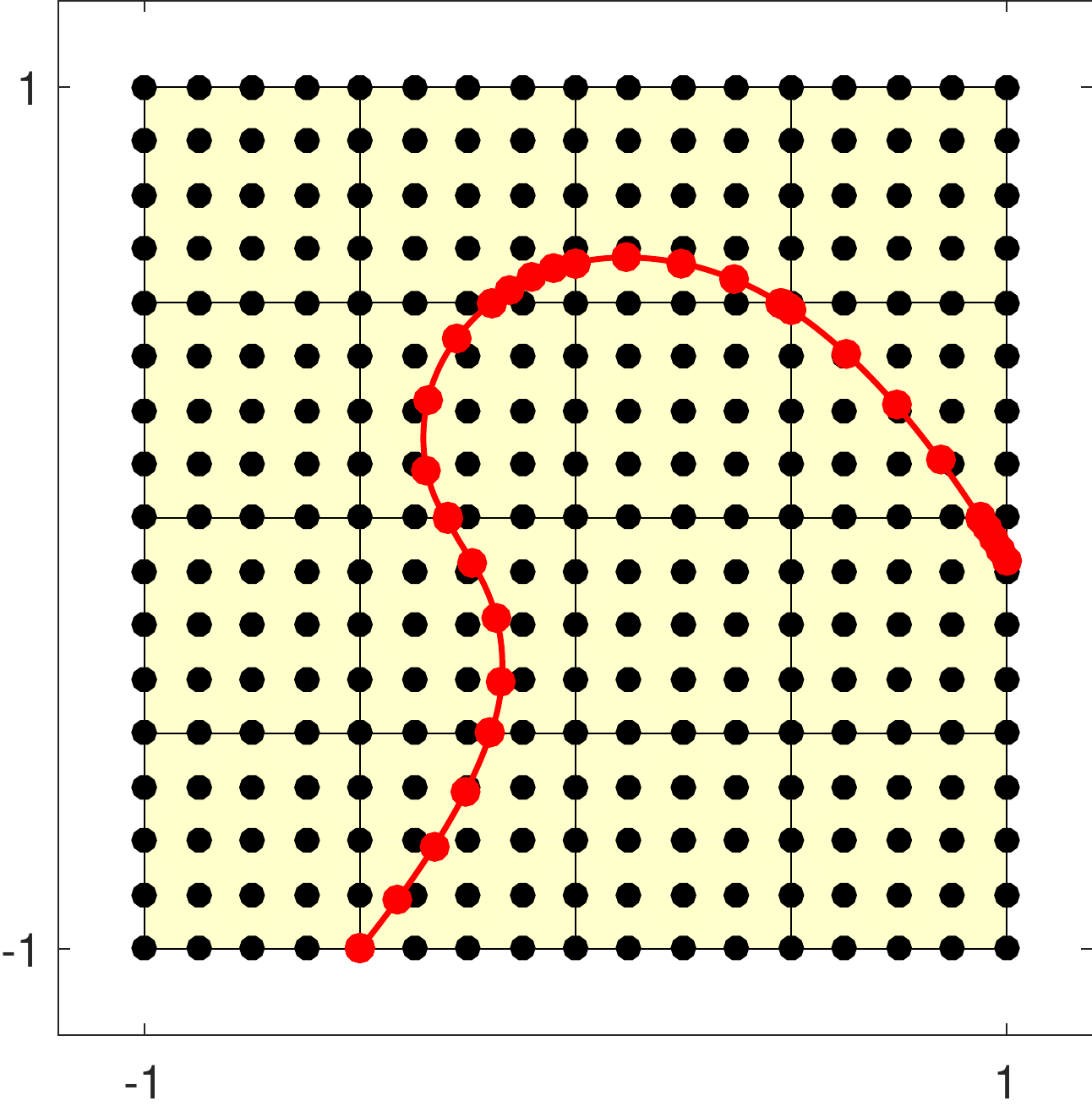}}

\caption{\label{fig:InvalidTopQuad}(a) Exact zero level set inside a higher-order
quadrilateral element. (b) The reconstruction fails due to the high
curvature of the level set. (c) and (d) The problem is resolved by
a uniform refinement of the element until no negative Jacobians are
detected.}
\end{figure}

As mentioned in Section~\ref{sec: Reconstruction}, the assessment
as a valid or invalid topology is based on level set data on the boundary
and does not take into account the shape of the level set inside the
element. However, in some valid configurations the reconstruction
may fail if the curvature of the level set function is too high. A
topological situation leading to such a case is shown in Fig.~\ref{fig:InvalidTopQuad}(a).
The reconstruction, see Fig.~\ref{fig:InvalidTopQuad}(b) fails clearly,
also the Runge phenomenon is apparent. However, as the mesh is successively
refined, see Fig.~\ref{fig:InvalidTopQuad}(c) and Fig.~\ref{fig:InvalidTopQuad}(d),
the situation is finally resolved. For the check, a finite number
of points is chosen and the Jacobian is evaluated. Here the points
are chosen as the integration points of the elements that emerge from
the subdivision scheme but also other choices are possible. The sub-elements
resulting from the decomposition are valid if the Jacobian is strictly
positive at all discrete points.

\subsection{Influence of the subdivision scheme on the condition number }

\begin{figure}
\centering
\subfigure[]{\includegraphics[width=6.5cm]{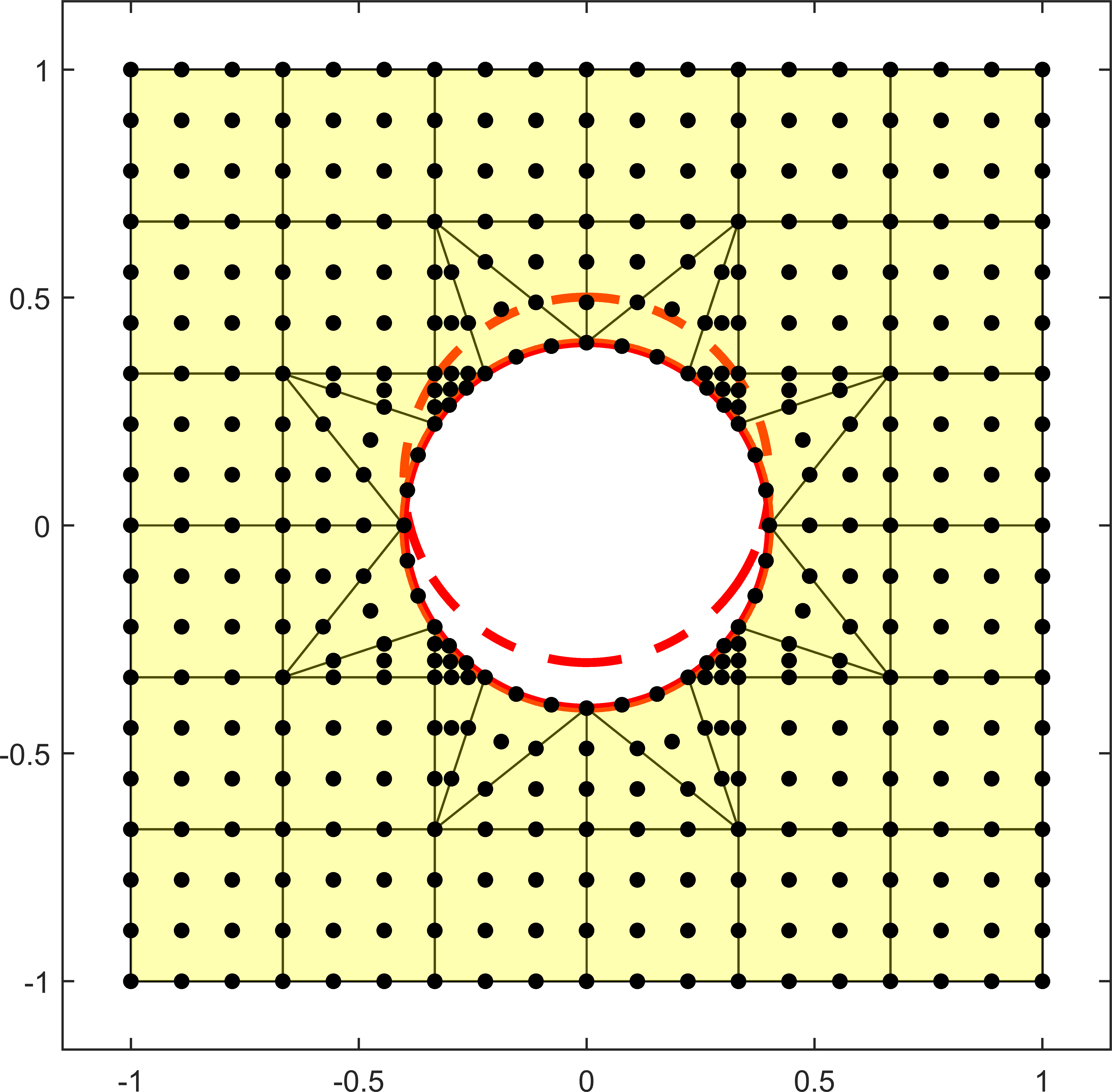}}$\qquad$\subfigure[]{\includegraphics[width=6.5cm]{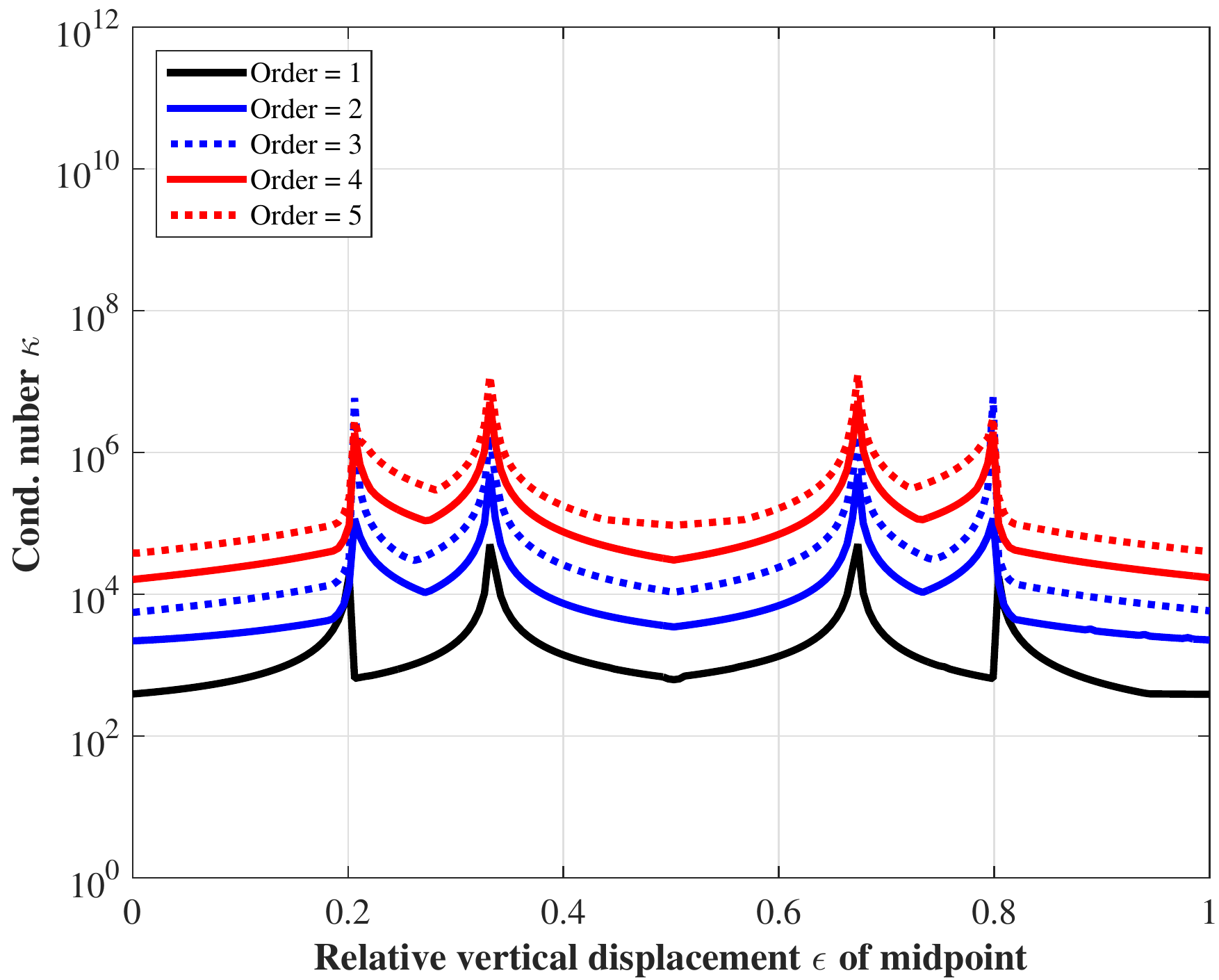}}

\caption{\label{fig:CondNum} (a) Moving discontinuity in a structured mesh
and (b) condition number in dependence of the relative displacement
of the discontinuity.}
\end{figure}

An obvious consequence of the proposed methodology is the presence
of awkward element shapes: Very small inner angles may occur and the
size of neighbouring elements can vary strongly. This may lead to
poorly conditioned system matrices, see e.g.~\cite{Frei_2014a}.
Some remedies for this problem are to use preconditioning schemes
or simple mesh manipulation procedures to prevent very small angles.
The influence on the condition number is shown here using the example
of a plate with a hole, which is later investigated in Section~\ref{sec:Numerical results}.
The mesh consists in total of $36$ elements, the element order is
varied between $m=1$ and $m=5$. The circular inclusion is moved
vertically and for each case the condition number is computed and
shown in Fig.~\ref{fig:CondNum}(b). 

Note that the peaks that represent a very high condition number in
the system of equations (resulting from a FE-analysis of the plate
with hole) emerge only for very few geometries in a small range of
the relative displacement. It is not difficult to identify nodes in
the background mesh which cause these high condition numbers. Several
techniques are available to regulate the situation. Among them are
node manipulation schemes as in \cite{Loehnert_2014a,Rangarajan_2014a,Labelle_2007a}
and stabilization methods as in \cite{Loehnert_2014a}. Herein, direct
solvers are employed for the solution of the system matrices and none
of these techniques are employed.

\section{Numerical results \label{sec:Numerical results}}

In order to illustrate the convergence properties and the stability
of the proposed method we present four different examples with meshes
involving non-fitting interfaces. In the test cases, we focus on (i)
the quality of interface reconstructions, (ii) the approximation properties
in the context of inner-element boundaries and interfaces, and (iii)
the conditioning of the system of equations. The first example is
a geometry approximation using a smooth level set function defining
an interface. The second example investigates an approximation problem
($L_{2}$-projection) of a smooth function on a remeshed, implicitly
defined circular domain. And finally, the third and fourth example
are boundary value problems in linear elasticity. The third example
is a bi-material problem and is representative for methods dealing
with (weak) discontinuities as e.g.~the XFEM. The fourth example
is used to demonstrate the capability of the procedure for a problem
typical for fictitious domain methods. Both BVPs are often used as
benchmarks for the performance of the $p$-version of the FEM.

\subsection{Interface reconstruction of a flower-shaped inclusion}

\begin{figure}
\centering
\subfigure[Cartesian background mesh]{\includegraphics[width=6.5cm]{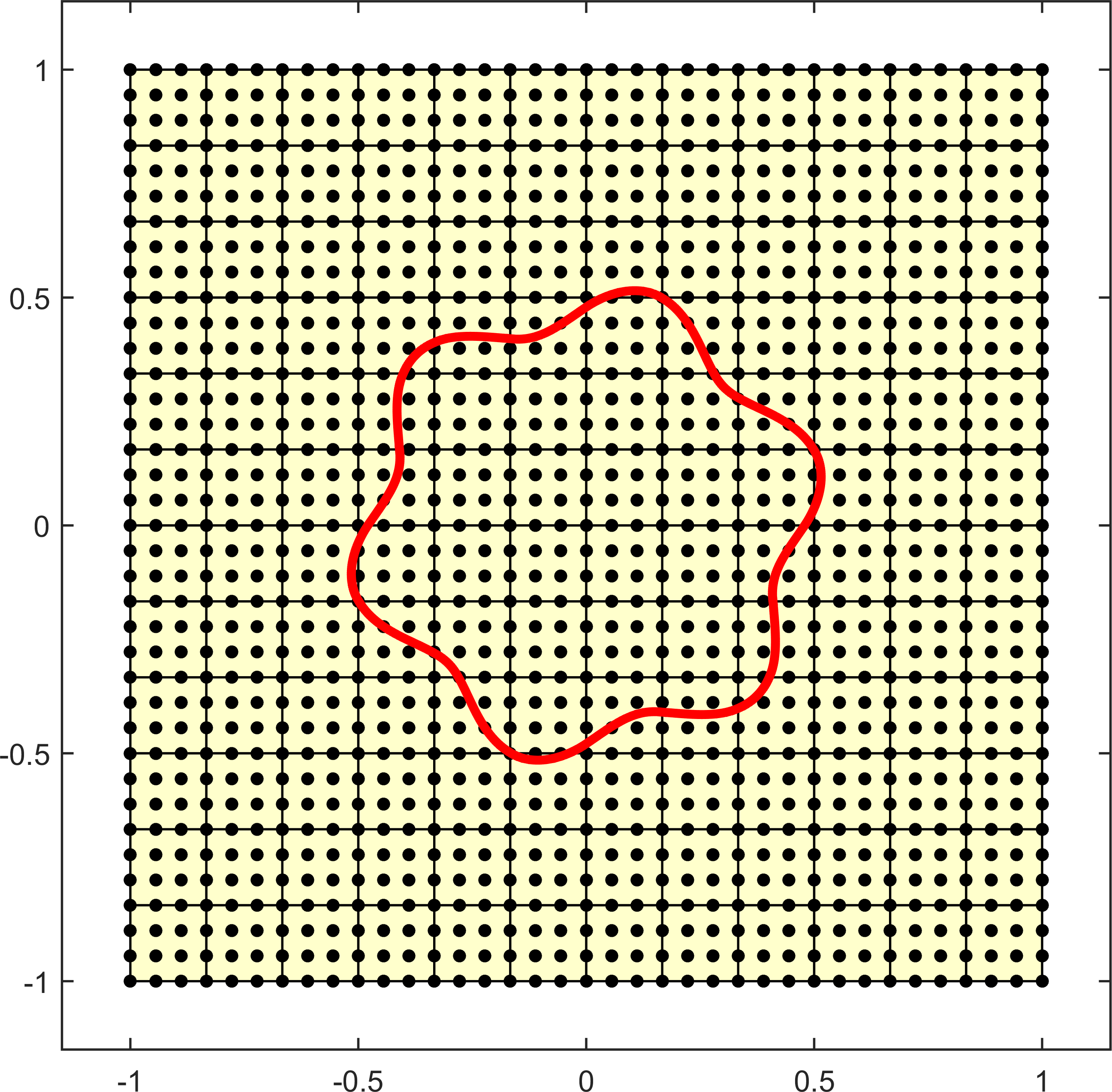}}$\qquad$\subfigure[Deformed background mesh]{\includegraphics[width=6.5cm]{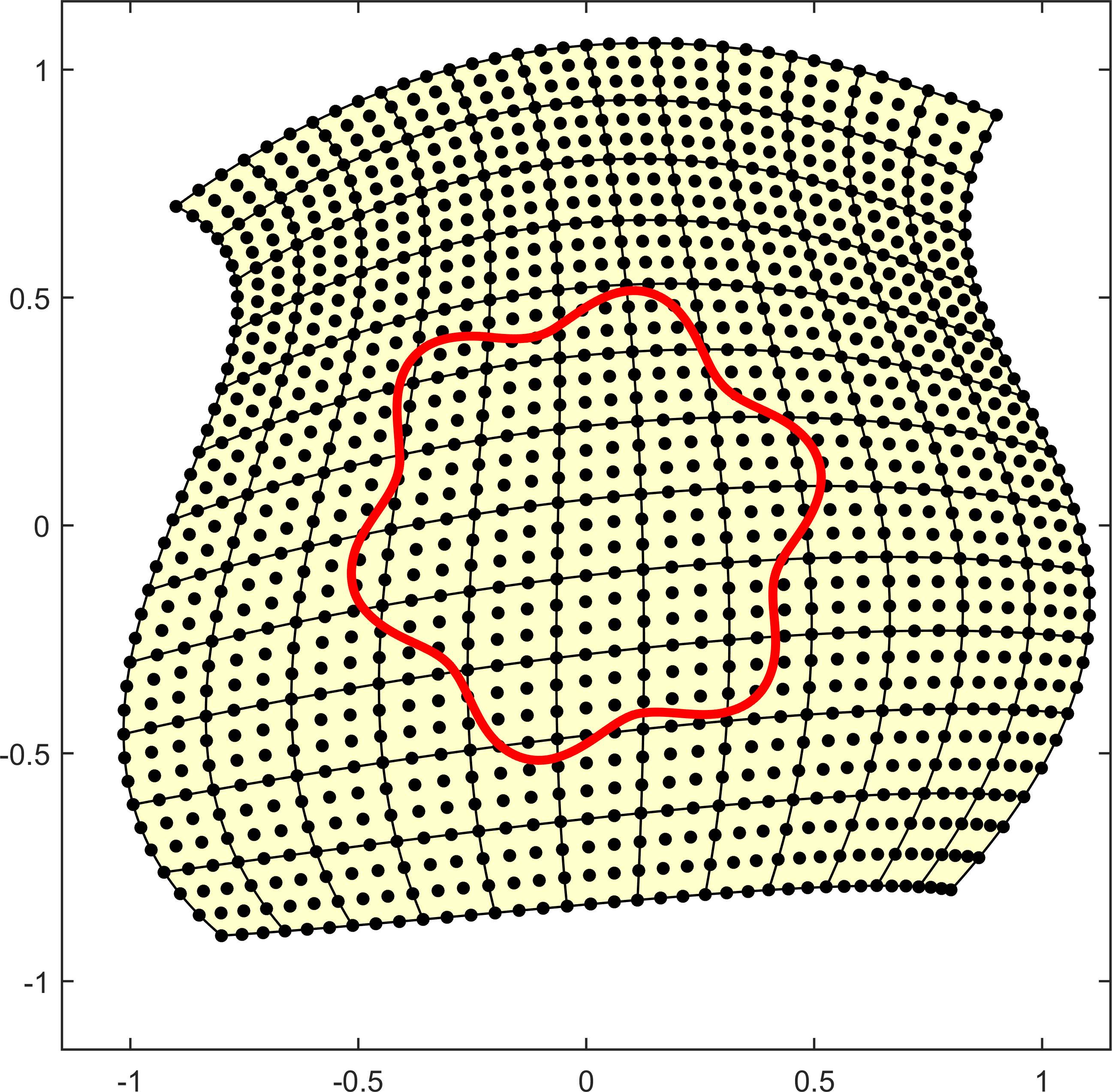}}

\caption{\label{fig:Example1Mesh} Higher-order background meshes with a flower-like
interface. }
\end{figure}

Given a reasonably well graded background mesh, this example shows
that using the proposed reconstruction scheme fairly complicated domains
can be automatically remeshed. The example is inspired by \cite{Li1998_a,Lehrenfeld_2016a}
with the analytical level set function (which could e.g.~describe
the boundary of an inclusion) given by
\begin{equation}
\phi(\vek x,\theta)=||\vek x||-(r+\alpha\sin(\omega\theta))\label{eq:Chap5Ex1LevSet}
\end{equation}
where $||\vek x||$ is the location of a considered point in $\Omega^{h}$,
$\theta\in[0,2\pi]$ the angle connecting the origin and $\vek x$,
and $\{r,\alpha,\omega\}$ a set of variables determining the shape
of the inclusion. For the considered example these variables are chosen
as $\{r=0.48,\alpha=0.05,\omega=6\}$. In Fig.~\ref{fig:Example1Mesh}(a),
a structured Cartesian background mesh is shown and in Fig.~\ref{fig:Example1Mesh}(b)
a deformed background mesh. In the figures, quadrilateral elements
of order $m=3$ are shown. Also the interface is seen, which is the
zero level set of Eq.~(\ref{eq:Chap5Ex1LevSet}).

We study the error at integration points which are distributed inside
the reconstructed interface elements and evaluate the analytical level
set function in these points. Because the interface elements provide
an approximation of the zero level set, a suitable error measure is
defined as 
\begin{equation}
\varepsilon=\sqrt{\int_{\Gamma}\phi(\vek x_{\Gamma})^{2}\,\mathrm{d}\Gamma}.\label{eq:Sec5Ex1Err}
\end{equation}
There, $\vek x_{\Gamma}$ are integration points in the reconstructed
interface and $\phi(\vek x_{\Gamma})$ is the analytical level set
function evaluated at $\vek x_{\Gamma}$. The level set function is
evaluated at the nodes of the background mesh and interpolated inside
the domain using finite element functions. Note that the Newton-Raphson
procedure from Section \ref{sec: Reconstruction} ensures that the
element \emph{nodes} are ``exactly'' on the zero level set of $\phi^{h}(\vek x)$
(with a tolerance of $\approx10^{-12}$). However, for the implied
interface element \emph{in-between}, and consequently also the integration
points, a slight deviation from the zero level set of $\phi^{h}$
occurs. The error norm in Eq.~(\ref{eq:Sec5Ex1Err}) accounts for
the error in $\phi^{h}$ compared to $\phi$ \emph{and} for the error
in the interface elements in approximating the zero level set of $\phi^{h}$.

\begin{figure}
\centering
\subfigure[Cartesian background mesh]{\includegraphics[width=6.5cm]{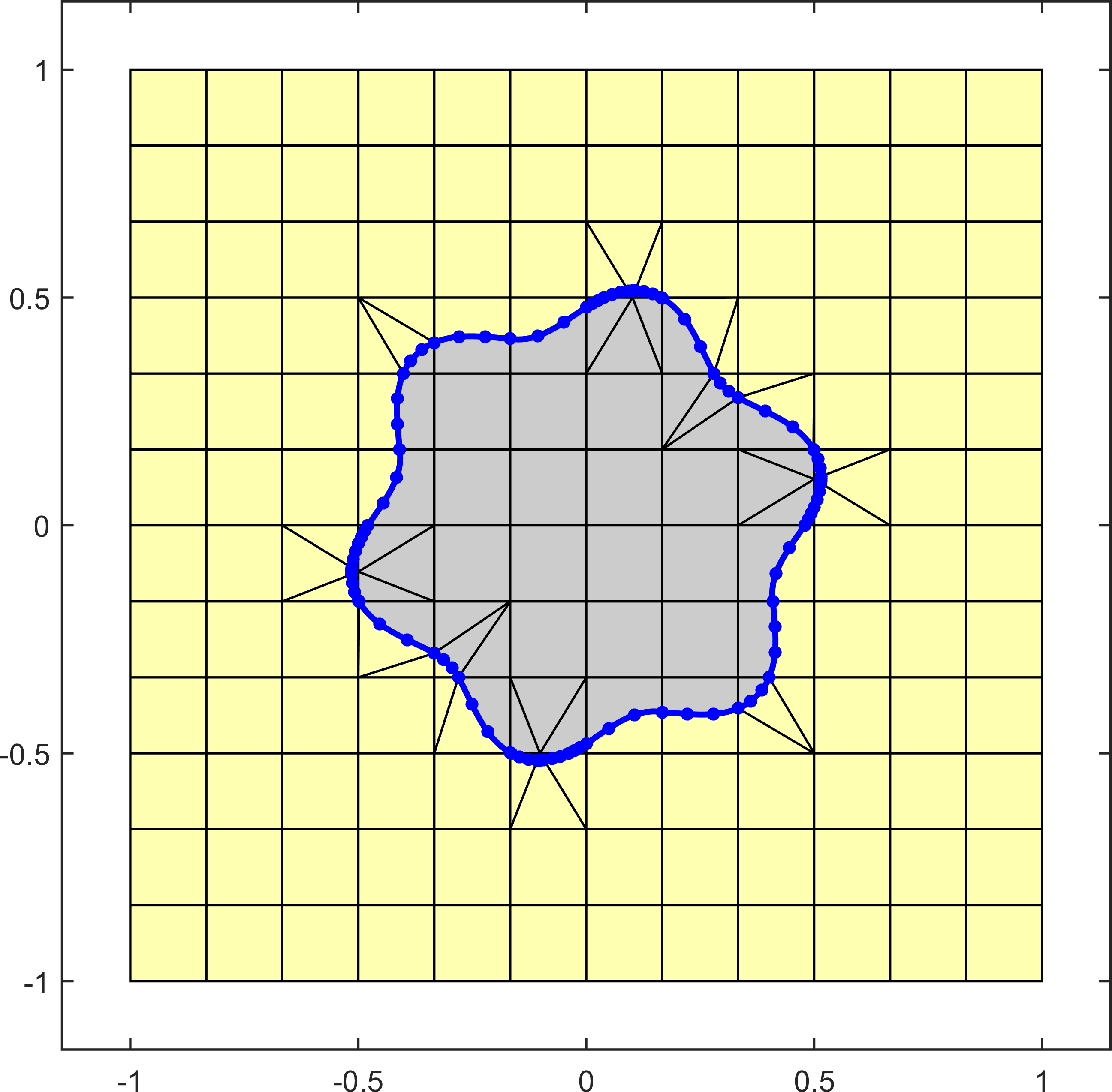}}$\qquad$\subfigure[Deformed background mesh]{\includegraphics[width=6.5cm]{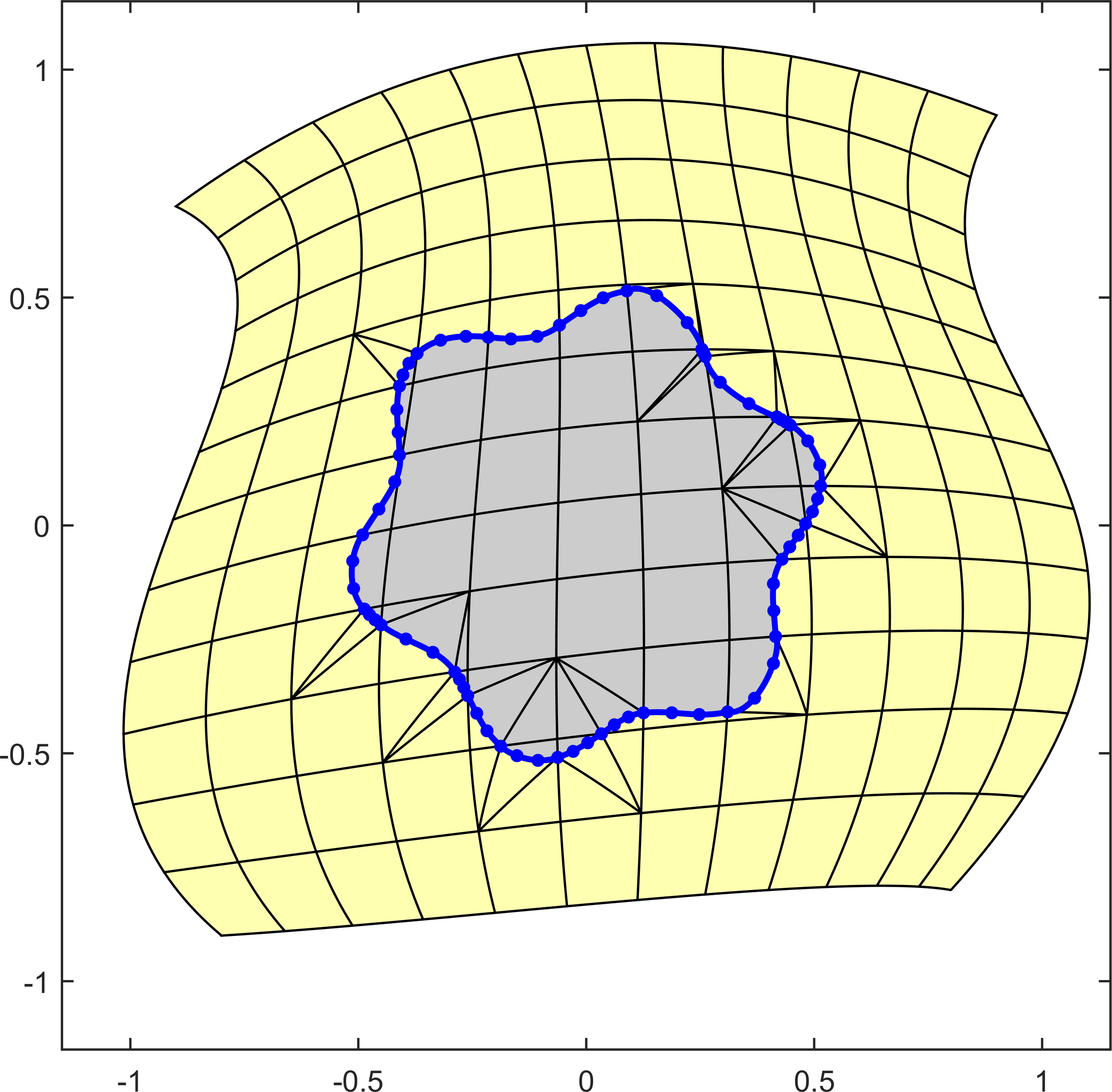}}

\caption{\label{fig:Example1MeshRecon} Updated meshes aligning with the inclusion
for (a) a Cartesian background mesh and (b) a deformed background
mesh.}
\end{figure}

In Fig.~\ref{fig:Example1MeshRecon}, the automatically generated
conformal meshes are visualized. It is emphasized once more, that
in general a mixed mesh consisting of triangular and quadrilateral
elements results. For clarity the element nodes are not plotted in
the meshes.

\begin{figure}
\centering
\subfigure[Cartesian background mesh]{\includegraphics[width=6.5cm]{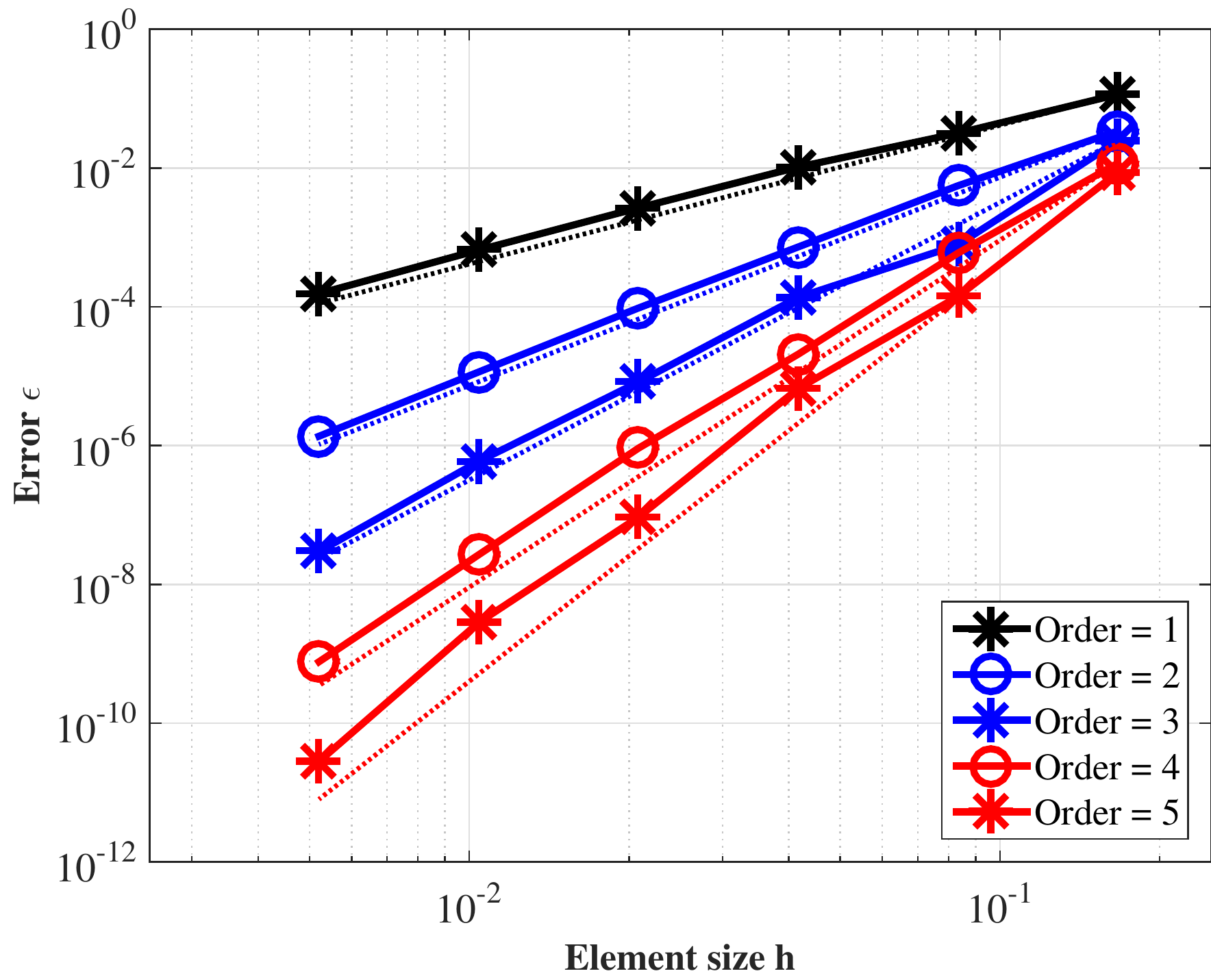}}$\qquad$\subfigure[Deformed background mesh]{\includegraphics[width=6.5cm]{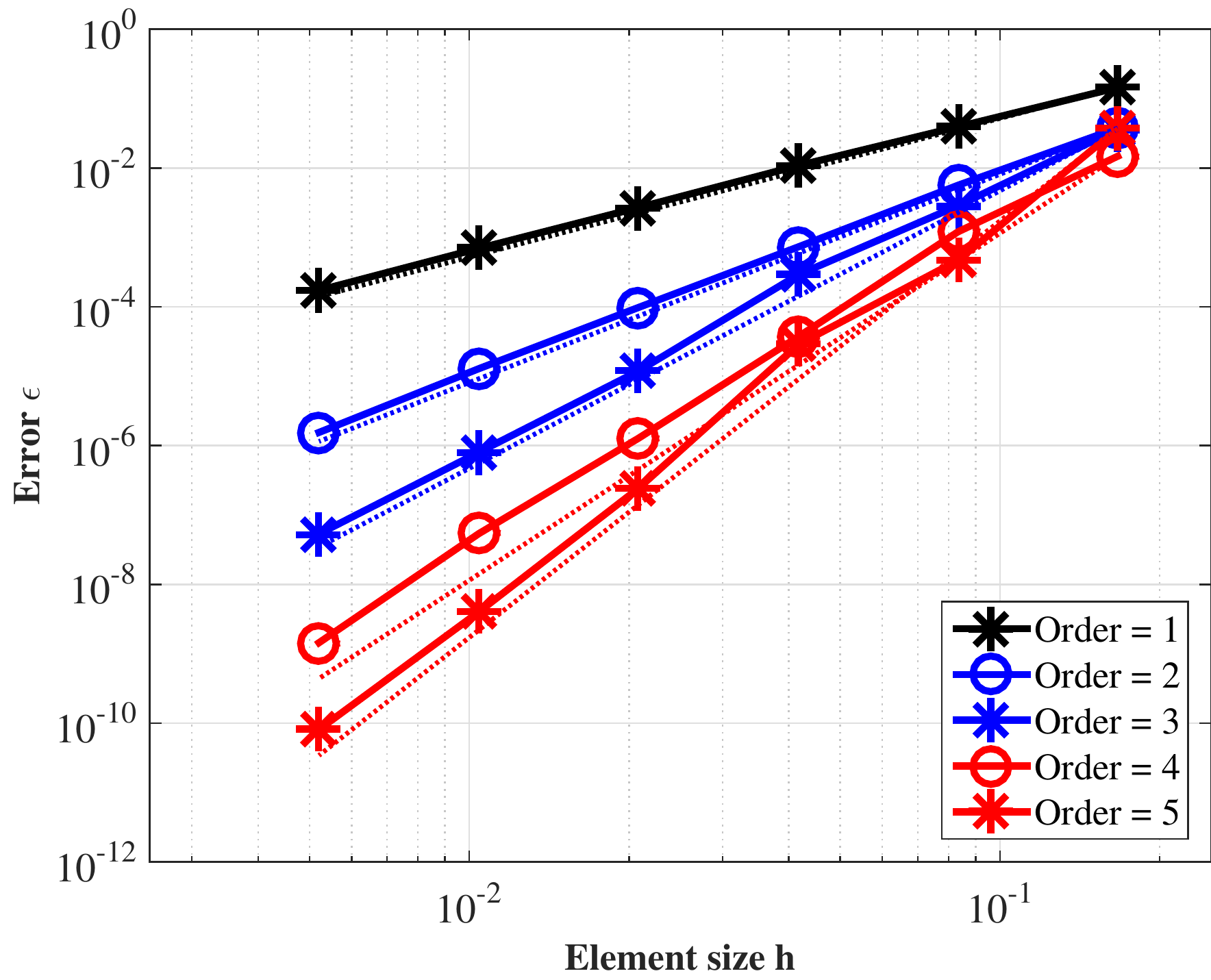}}

\caption{\label{fig:Example1MeshResultConv} $L_{2}$-error, as defined in
Eq.~(\ref{eq:Sec5Ex1Err}), of the interface reconstruction.}
\end{figure}

The analysis is carried out using a series of embedded meshes, indicated
by the characteristic element length $h=diam(K)$, measured in the
background mesh (before the decomposition). As shown in Fig.~\ref{fig:Example1MeshResultConv}(a)
and Fig.~\ref{fig:Example1MeshResultConv}(b), the convergence properties
of the method, based on the Cartesian background mesh as well as on
the deformed background mesh, are optimal after the automatic remeshing.
Note that no node moving was applied and some of the elements look
quite awkward. Nevertheless all elements resulting from the subdivision
procedure are topologically valid. However, for higher curvatures
inside the elements it can be advisable to apply some mesh manipulations
before the subdivision scheme is applied.

\subsection{Approximation problem with circular inclusion}

\begin{figure}
\centering
\subfigure[Cartesian background mesh]{\includegraphics[width=6.5cm]{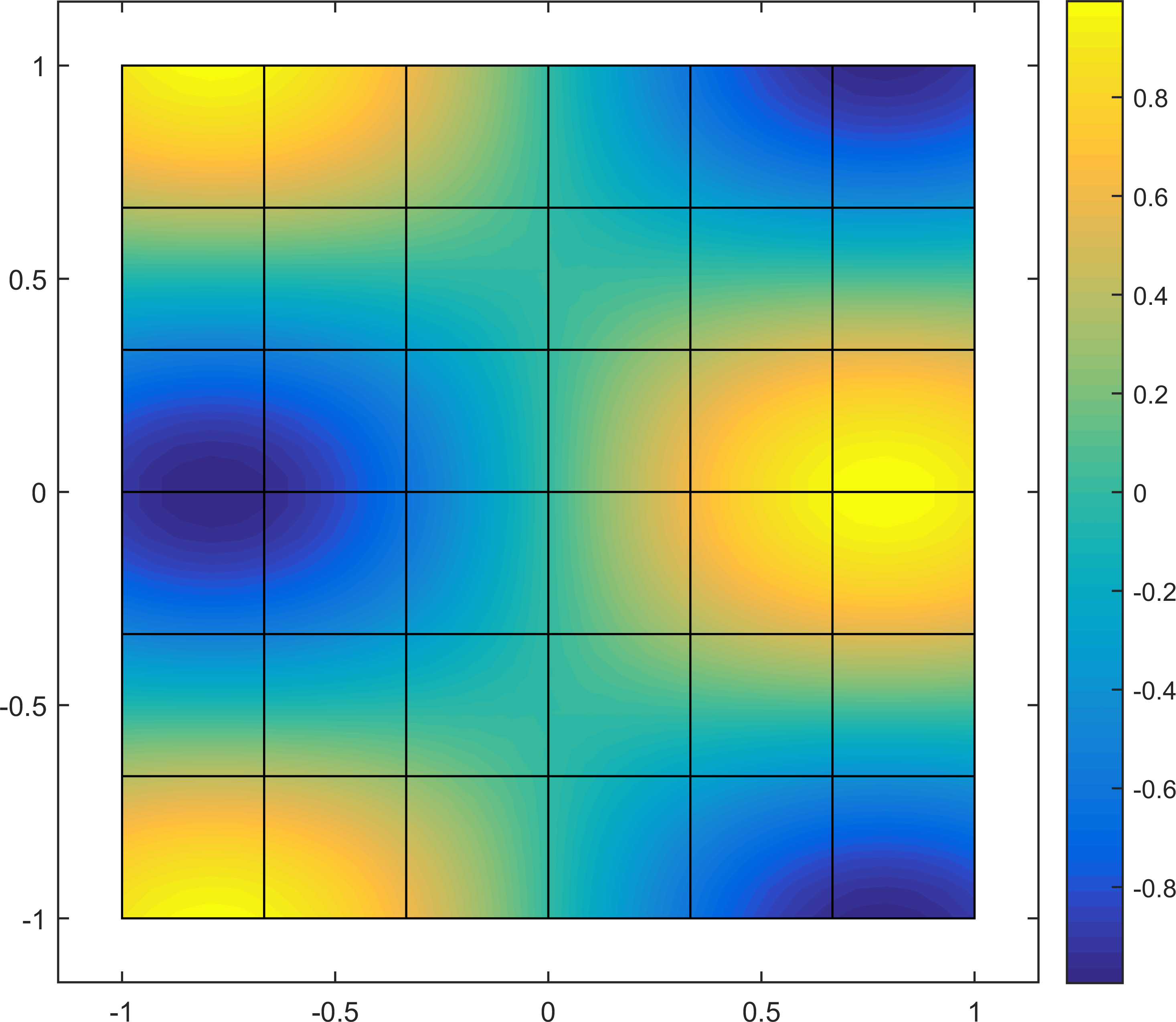}}$\qquad$\subfigure[Deformed background mesh]{\includegraphics[width=6.5cm]{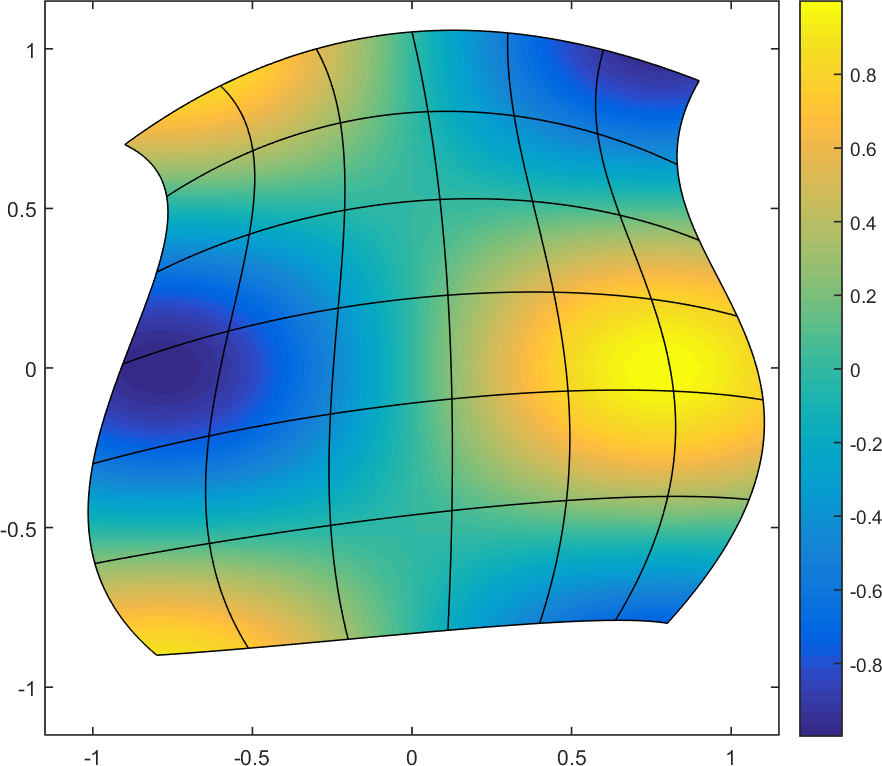}}

\caption{\label{fig:Chap5Ex2Fct} Plot of the function for which an $L_{2}$-projection
is made.}
\end{figure}

\begin{figure}
\centering
\subfigure[Cartesian background mesh]{\includegraphics[width=6.5cm]{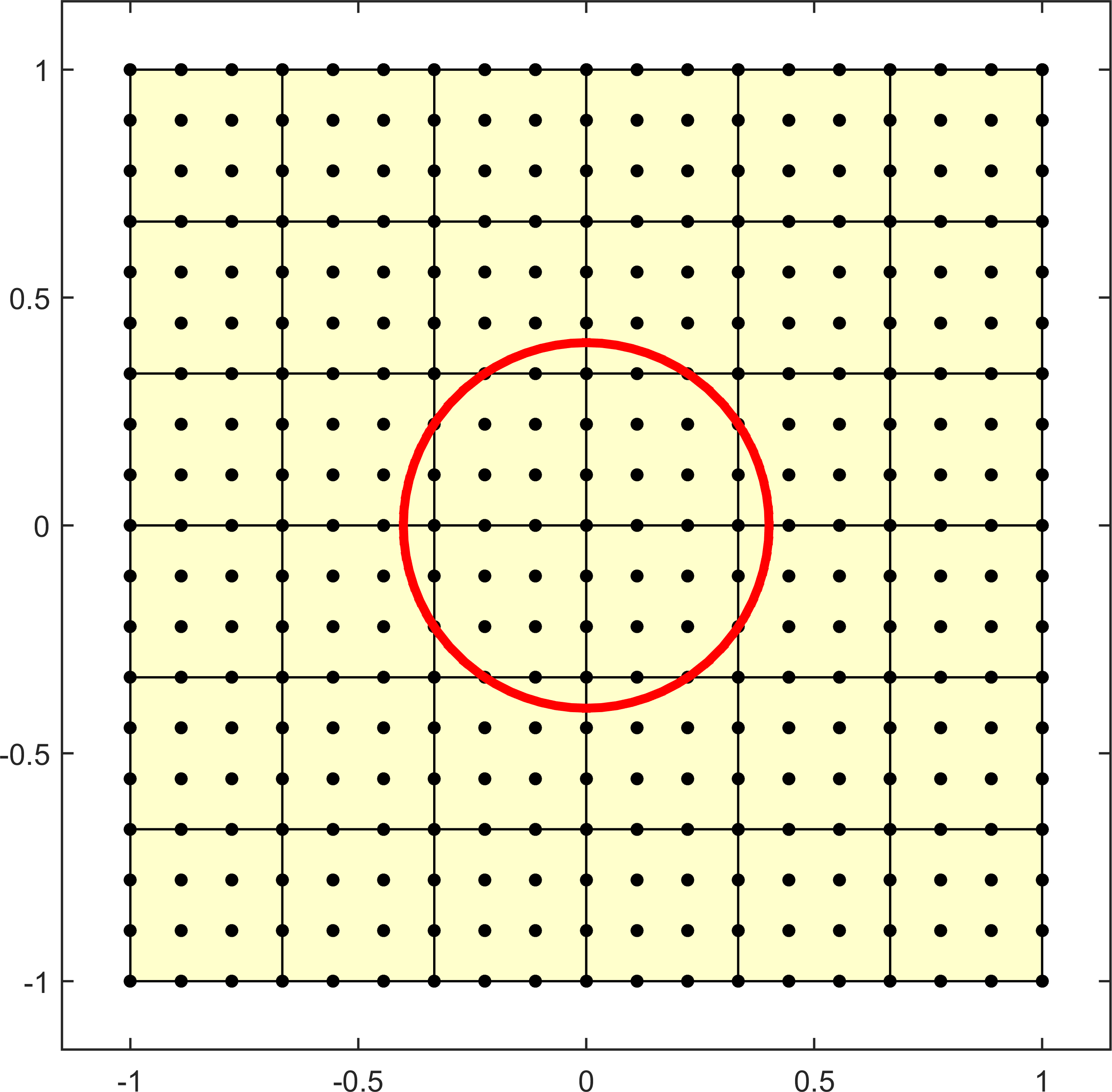}}$\qquad$\subfigure[Deformed background mesh]{\includegraphics[width=6.5cm]{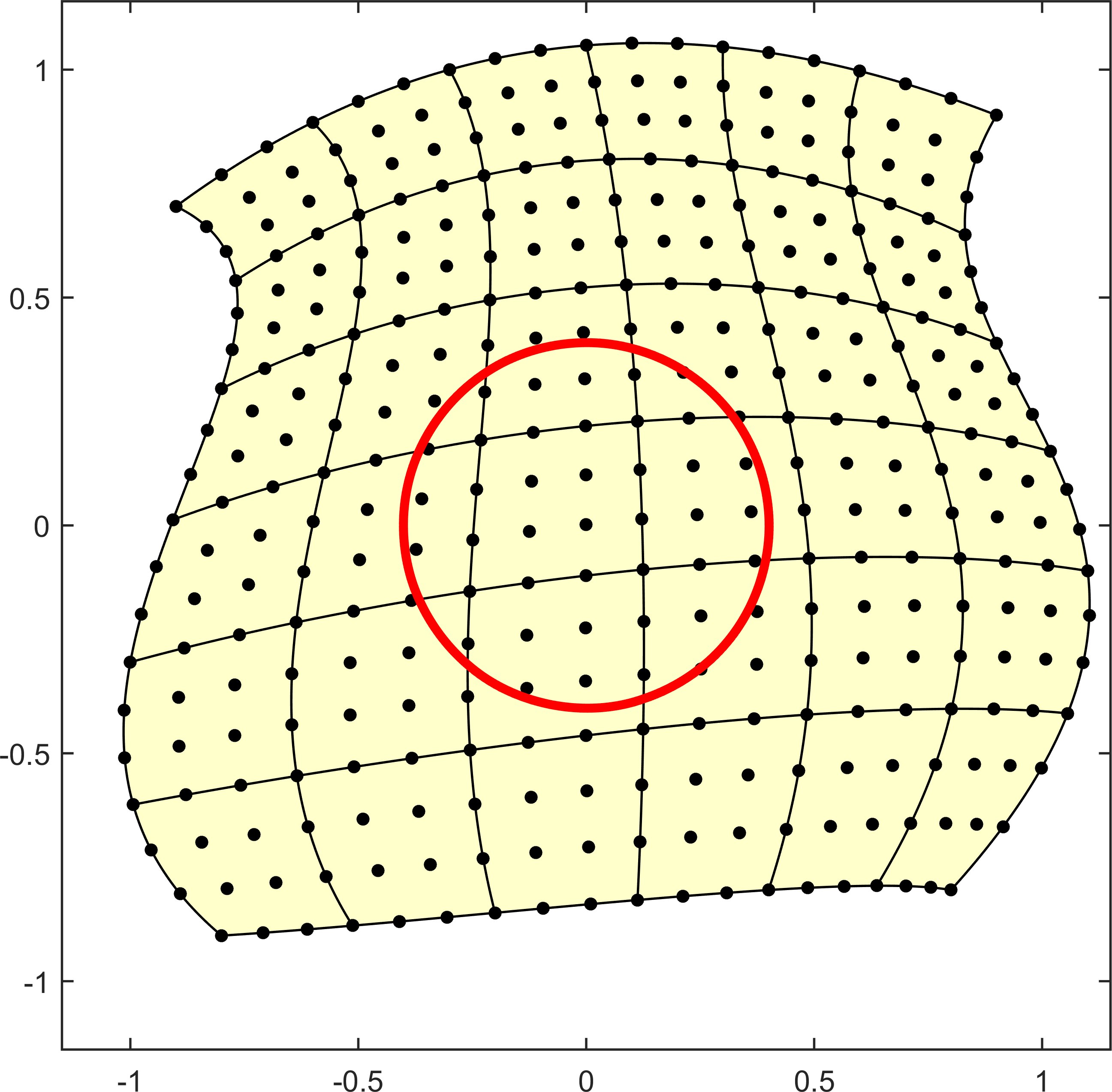}}\\\subfigure[Conformal mesh based on (a)]{\includegraphics[width=6.5cm]{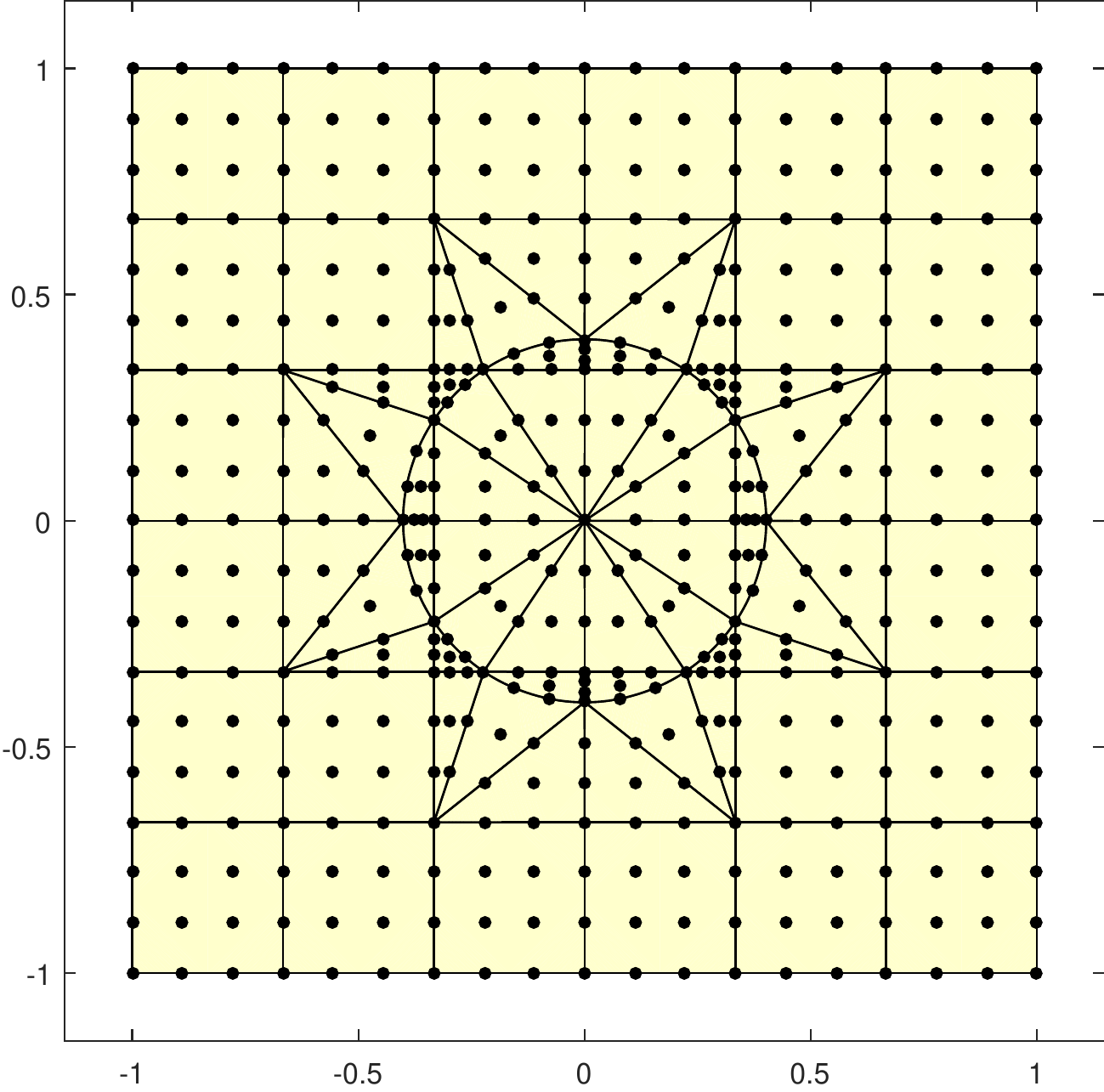}}$\qquad$\subfigure[Conformal mesh based on (b)]{\includegraphics[width=6.5cm]{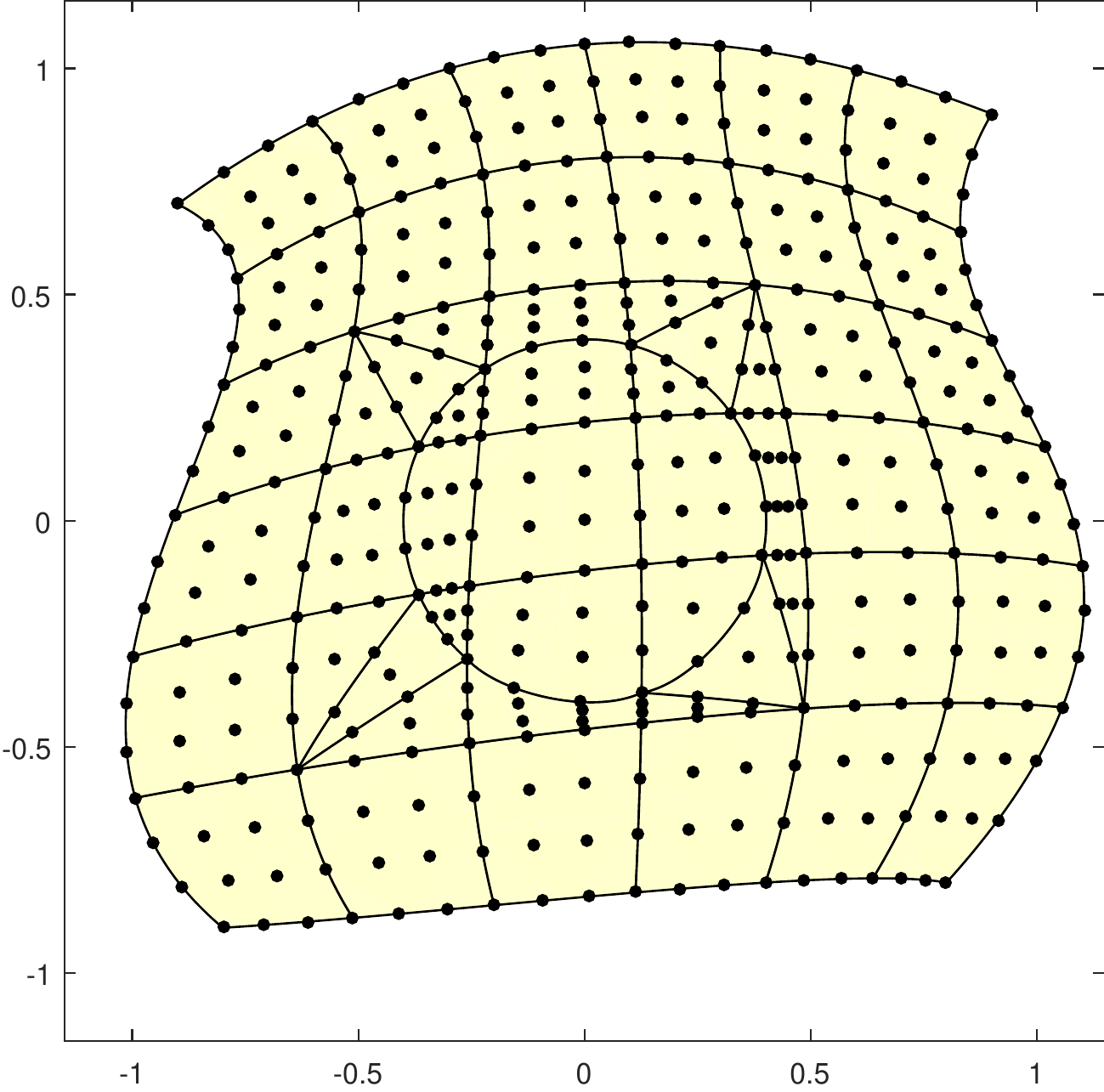}}

\caption{\label{fig:Chap5Ex2Meshes} Plot of the involved meshes: (a) and (b)
are background meshes, (c) and (d) show the automatically generated
conforming meshes.}
\end{figure}

The second example shows that the automatic decomposition of the background
mesh retains optimal approximation properties of the standard $L_{2}$-projection.
Also the convergence properties of the different mappings for triangular
elements, as introduced in Section \ref{sec: Remeshing}, are investigated.
Again, a Cartesian and a deformed background mesh, respectively, are
automatically decomposed to align with a circular inclusion, see Fig.~\ref{fig:Chap5Ex2Meshes}(a)
and Fig.~\ref{fig:Chap5Ex2Meshes}(b). The radius of the inclusion
is given as $r=0.4$ and centred in the origin. The function to be
interpolated is given as $f(x,y)=\sin(2x)\cos(3y)$ and therefore
smooth over the whole domain. In Fig.~\ref{fig:Chap5Ex2Fct}(a) and
Fig.~\ref{fig:Chap5Ex2Fct}(b) plots of the function over the background
meshes are shown. 

The error is measured in the standard $L_{2}$-norm as well as in
the $H_{E}$-norm, to show that in both error measures optimal convergence
rates are obtained. 

To begin with, we consider the approximation error of $f$ in the
$L_{2}$-norm defined by 
\[
\varepsilon_{L_{2}}=\sqrt{\int_{\Omega}\frac{(f(\vek x)-f^{h}(\vek x))^{2}}{f(\vek x)^{2}}\,\mathrm{d}\Omega}\qquad\forall\vek x\in\Omega
\]

\begin{figure}
\centering
\subfigure[Cartesian background mesh]{\includegraphics[width=6.5cm]{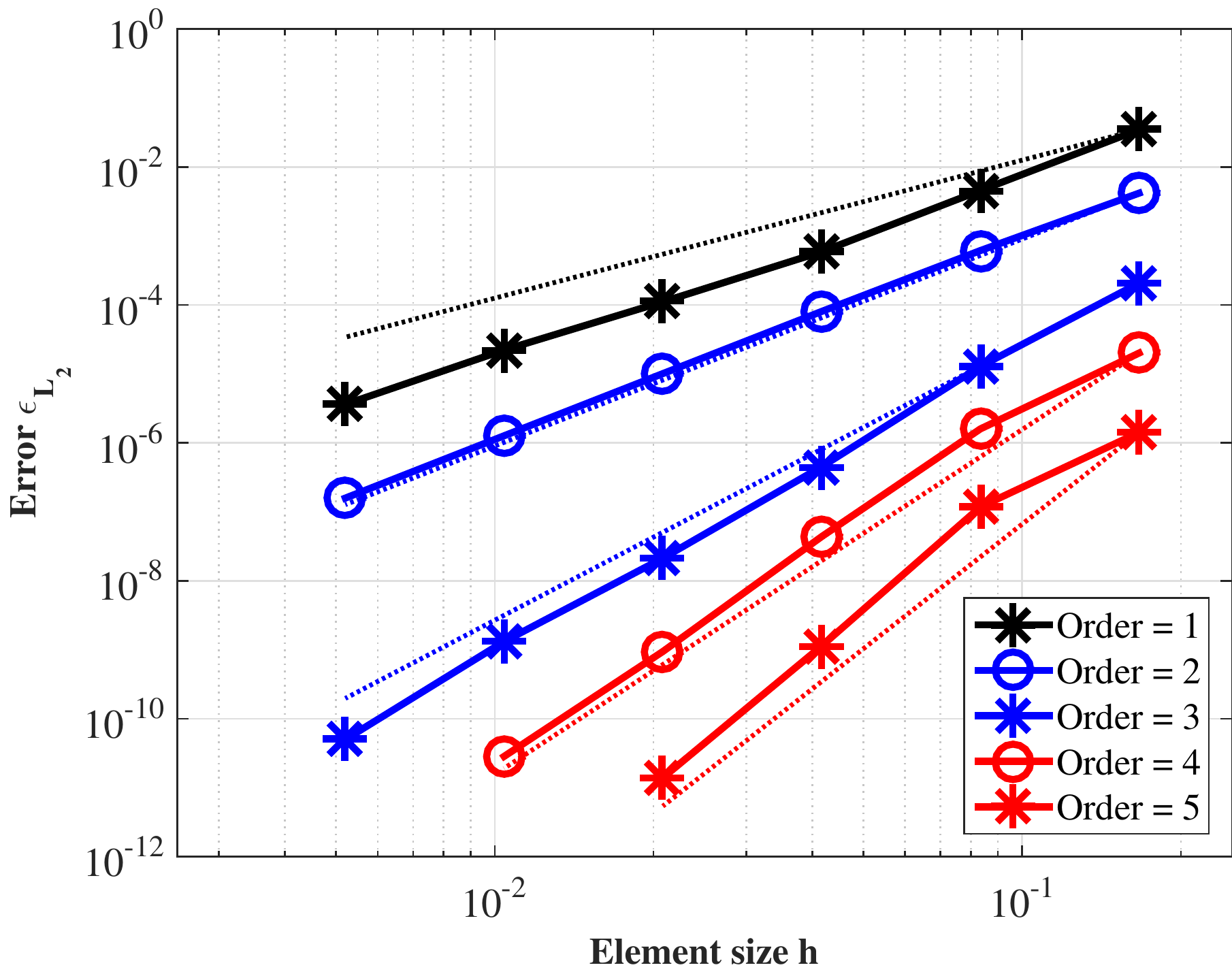}}$\qquad$\subfigure[Deformed background mesh]{\includegraphics[width=6.5cm]{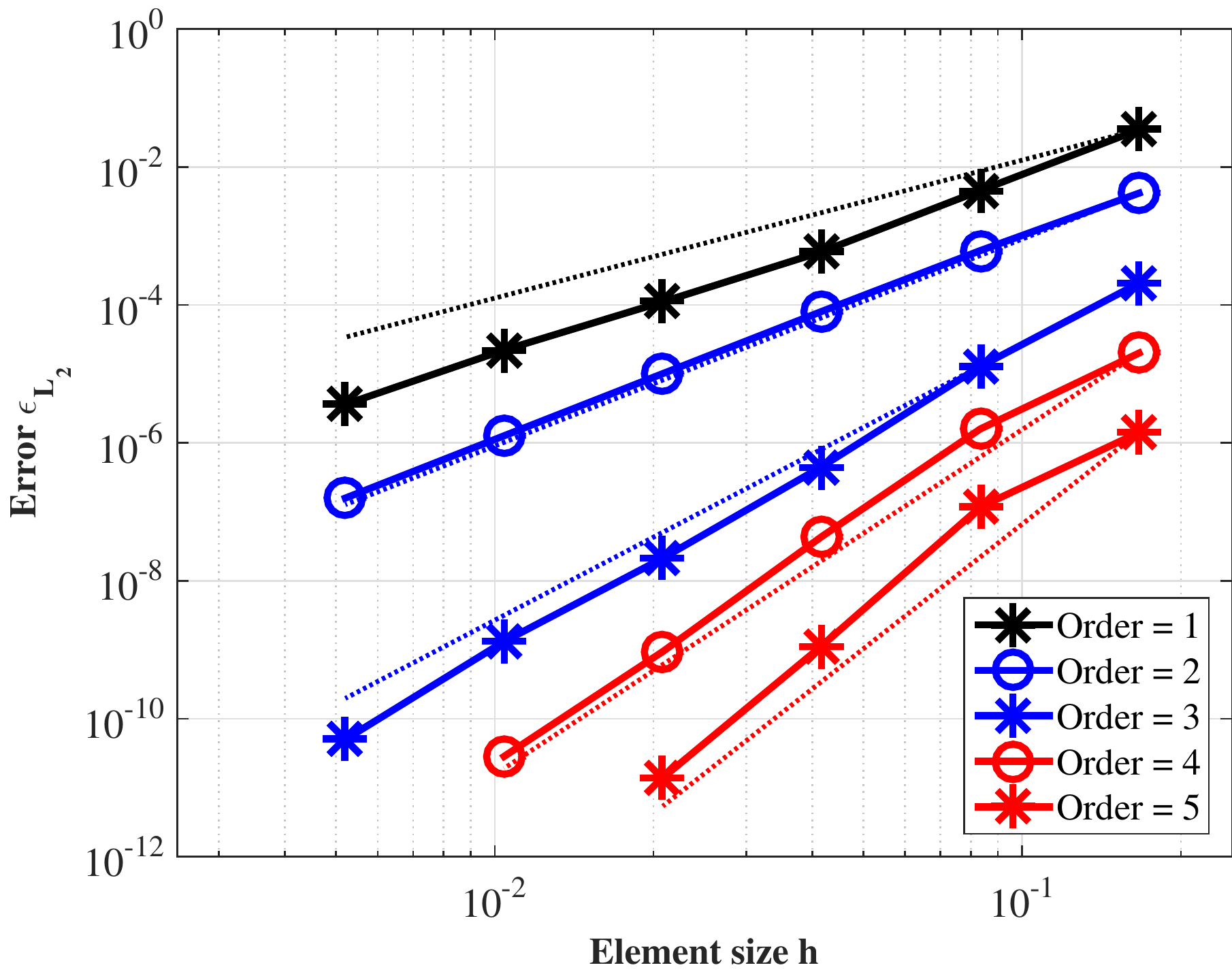}}

\caption{\label{fig:Chap5Ex2ResultL2} Approximation error measured in the
$L_{2}$-norm. }
\end{figure}

As depicted in Fig.~\ref{fig:Chap5Ex2ResultL2}, optimal convergence
rates are achieved. Also the \emph{gradient} approximation error is
investigated in the same way, the error is measured in the energy
norm $H_{E}$

\[
\varepsilon_{H_{E}}=\sqrt{\int_{\Omega}\frac{(\nabla f(\vek x)-\nabla f^{h}(\vek x))^{2}}{\nabla f^{h}(\vek x)^{2}}\,\textrm{\ensuremath{\mathrm{d}\Omega}}}\qquad\forall\vek x\in\Omega
\]

\begin{figure}[h]
\centering
\subfigure[Cartesian background mesh]{\includegraphics[width=6.5cm]{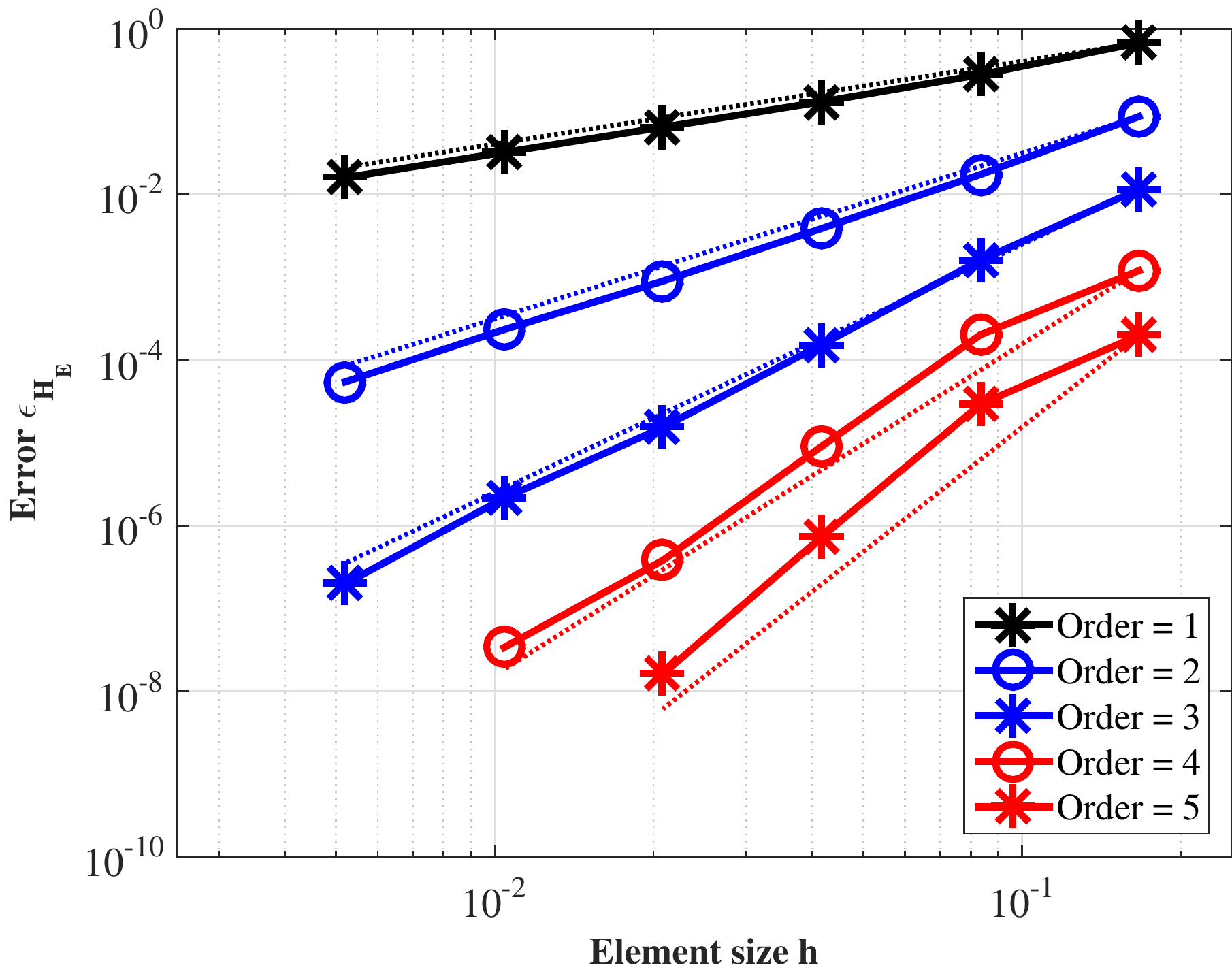}}$\qquad$\subfigure[Deformed background mesh]{\includegraphics[width=6.5cm]{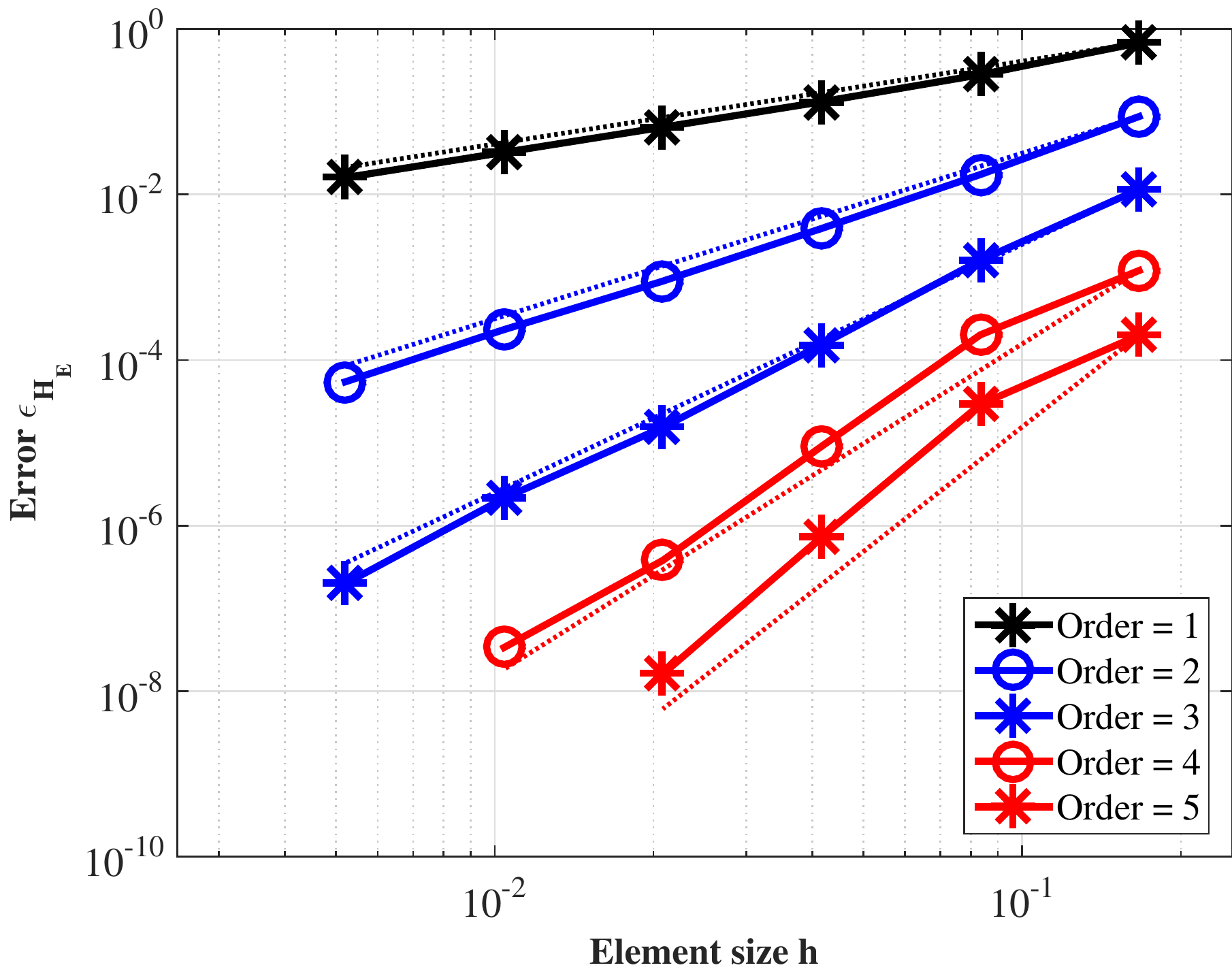}}

\caption{\label{fig:Chap5Example2ConvResultHE} Approximation error measured
in the $H_{E}$-Norm.}
\end{figure}

The results in Fig.~\ref{fig:Chap5Ex2ResultL2} and Fig.~\ref{fig:Chap5Example2ConvResultHE}
show that the remeshing strategy, indeed, achieves optimal convergence
properties asymptotically, in both the $L_{2}$-norm and the $H_{E}$-norm. 

For the sake of completeness, the $L_{2}$-error for the other mappings
of the element nodes to the sub-cells, discussed in Section \ref{sec: Remeshing},
namely the \emph{Lenoir-mapping} and the \emph{blending function mapping},
is also shown for the Cartesian background mesh. The Lenoir-mapping
gives optimal approximation properties (but is more difficult to implement
than the actually applied mapping). The blending function mapping,
on the contrary, clearly gives sub-optimal results and degrades the
convergence rates significantly. This is due to the properties of
the Jacobian and well-known in the meshing community.

\begin{figure}
\centering
\subfigure[Blending function mapping]{\includegraphics[width=6.5cm]{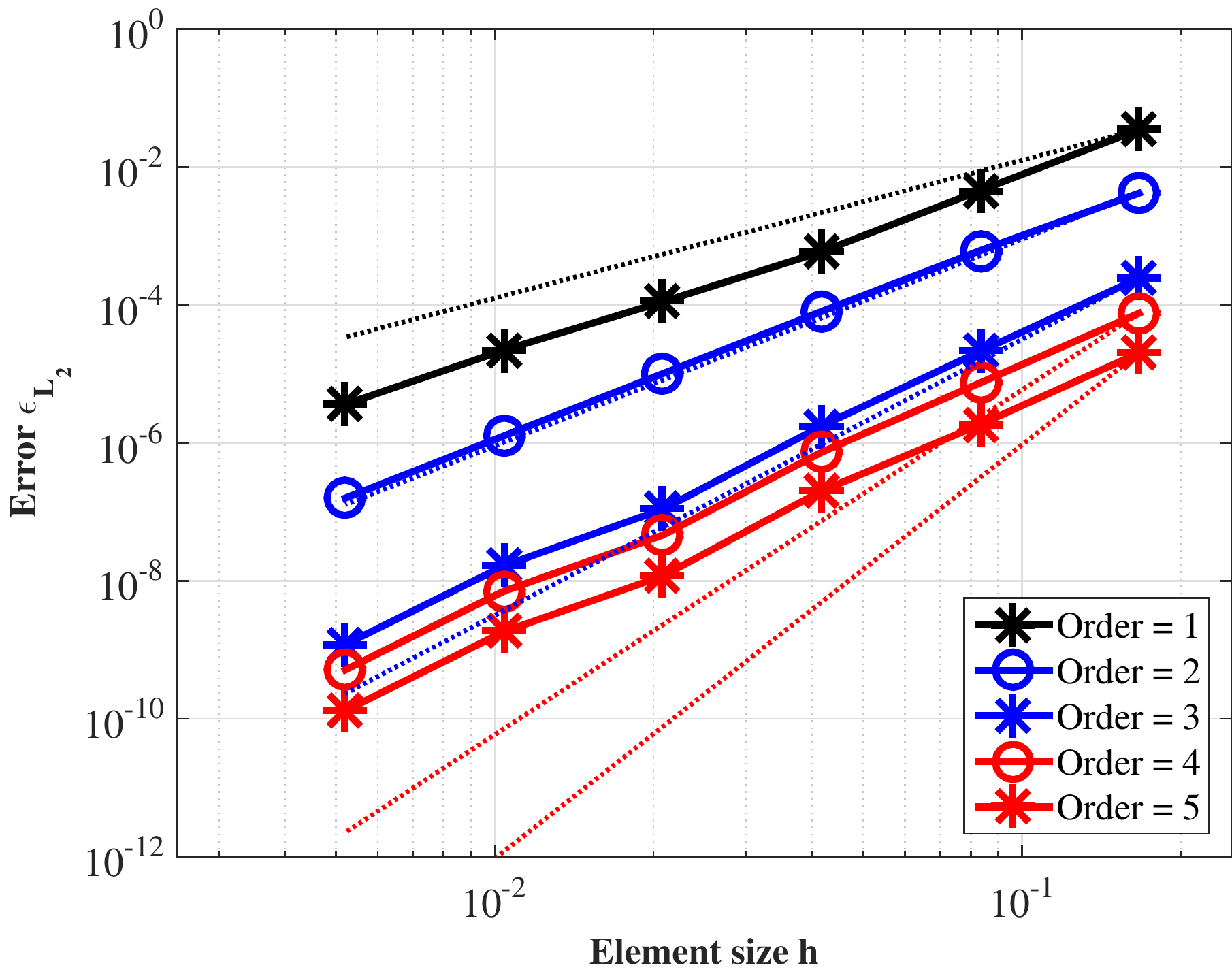}}\qquad\subfigure[Lenoir mapping]{\includegraphics[width=6.5cm]{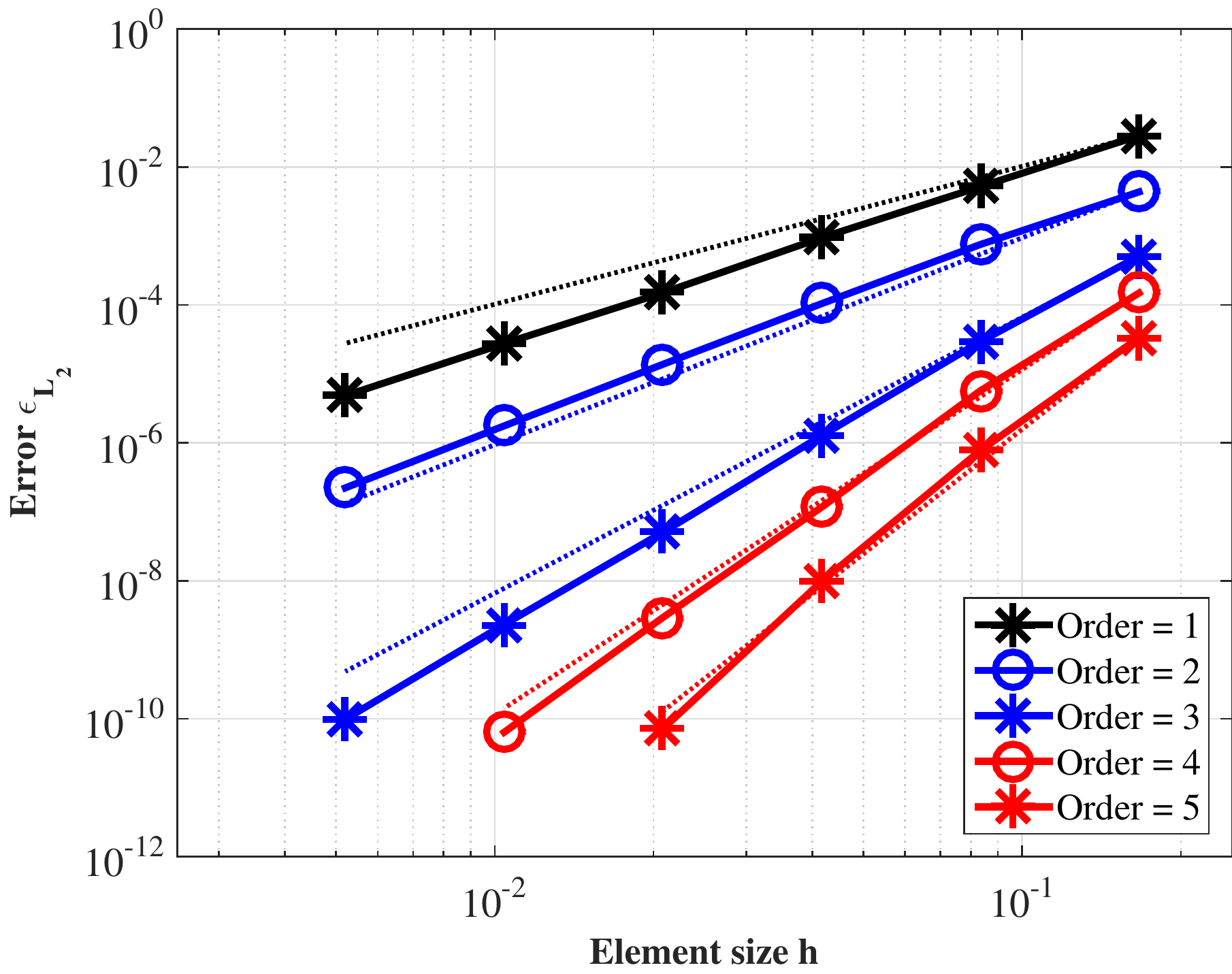}}

\caption{Approximation error measured in $L_{2}$-norm for the mappings introduced
in Section \ref{sec: Remeshing}. \label{fig:Example1MeshResultConv-1-1-1}.}
\end{figure}

\subsection{Bi-material boundary value problem}

As a third example, a plate with a circular inclusion is considered
to illustrate the performance for the approximation of general boundary-value
problems with a curved (weak) discontinuity. The example is a well
known benchmark problem for higher-order extended finite element procedures,
see e.g.~\cite{Cheng_2009a,Sukumar_2001c}. The error in the displacements
(measured in the $L_{2}$-norm), the error in stresses (measured in
the $H_{E}$-norm) and also the condition number of the resulting
system of equations are investigated. Regarding the employed error
measures for the convergence studies for this and also for the next
example, using the analytical solution $\vek u$ and the finite element
solution $\vek u^{h}$, the $L_{2}$-norm on $\Omega$ is defined
as
\begin{equation}
\varepsilon_{L_{2}}=||\vek u-\vek u^{h}||_{L_{2}(\Omega)}=\sqrt{\int(\vek u-\vek u^{h}){}^{\intercal}(\vek u-\vek u^{h})\,\mathrm{d}\Omega}.\label{eq:L2Norm}
\end{equation}
To investigate the convergence behaviour of the function derivatives,
we use the strain energy $H_{E}$, defined as
\begin{equation}
\varepsilon_{H_{E}}=||\vek u-\vek u^{h}||_{H_{E}(\Omega)}=\sqrt{\int(\vek\epsilon-\vek\epsilon^{h}){}^{\intercal}\mathbf{C}(\vek\epsilon-\vek\epsilon^{h})\,\mathrm{d}\Omega}\label{eq:HENorm}
\end{equation}
with $\mathbf{C}$ being the elastic stiffness tensor and $\vek\epsilon(\vek u)$
and $\vek\epsilon^{h}(\vek u^{h})$ the strain fields of the exact
solution and the finite element approximation, respectively. We use
relative versions of these two norms, therefore $\varepsilon_{L_{2}(\Omega)}=||\vek u-\vek u^{h}||_{L_{2}(\Omega)}/||\vek u||_{L_{2}(\Omega)}$
and $\varepsilon_{H_{E}(\Omega)}=||\vek u-\vek u^{h}||_{H_{E}(\Omega)}/||\vek u||_{H_{E}(\Omega)}$.
Assuming that a $C_{\infty}$-smooth function is approximated on a
regular domain using a polynomial of degree $k$ and with a discretization
aligning with the interface, optimal convergence rates are expected,
see \cite{Szabo_1991a,Strang_2008a}. Subject to these limitations,
the convergence rates in the error norms defined above are
\[
||\vek u-\vek u^{h}||_{L_{2}(\Omega)}=\mathcal{O}(h^{k+1})\quad\textrm{{and}}\quad||\vek u-\vek u^{h}||_{H_{E}(\Omega)}=\mathcal{O}(h^{k}).
\]
 In case of unfitted interfaces, e.g.~using a piecewise straight
approximation of curved boundaries or material interfaces crossing
the elements, the approximation quality is significantly degraded.
The convergence rates are then independent of the ansatz order of
the finite element functions and are limited by the smoothness of
the approximated function or the accuracy of the geometry description. 

\begin{figure}
\centering
\subfigure[Cartesian background mesh]{\includegraphics[width=6.5cm]{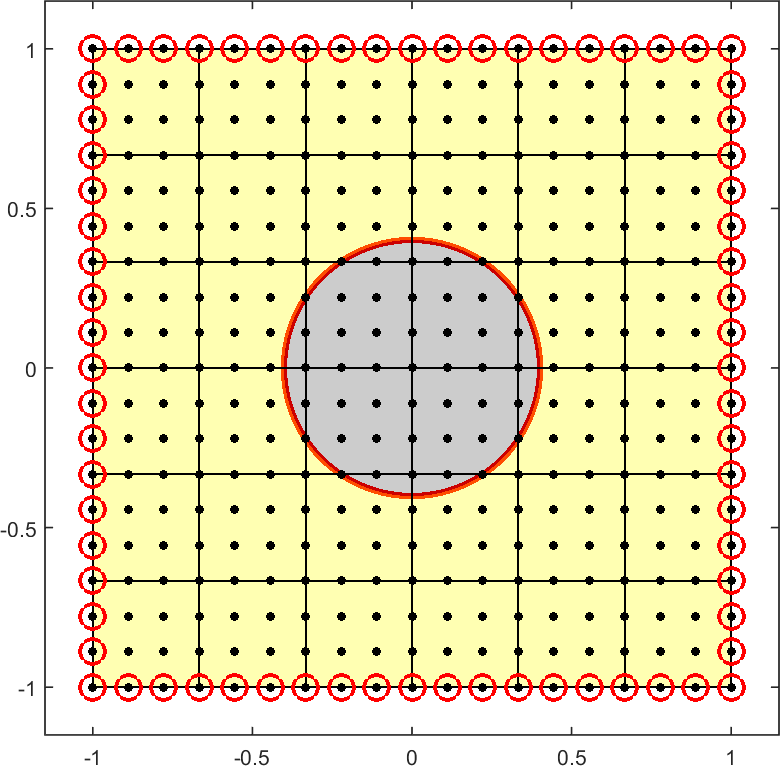}}$\qquad$\subfigure[Deformed background mesh]{\includegraphics[width=6.5cm]{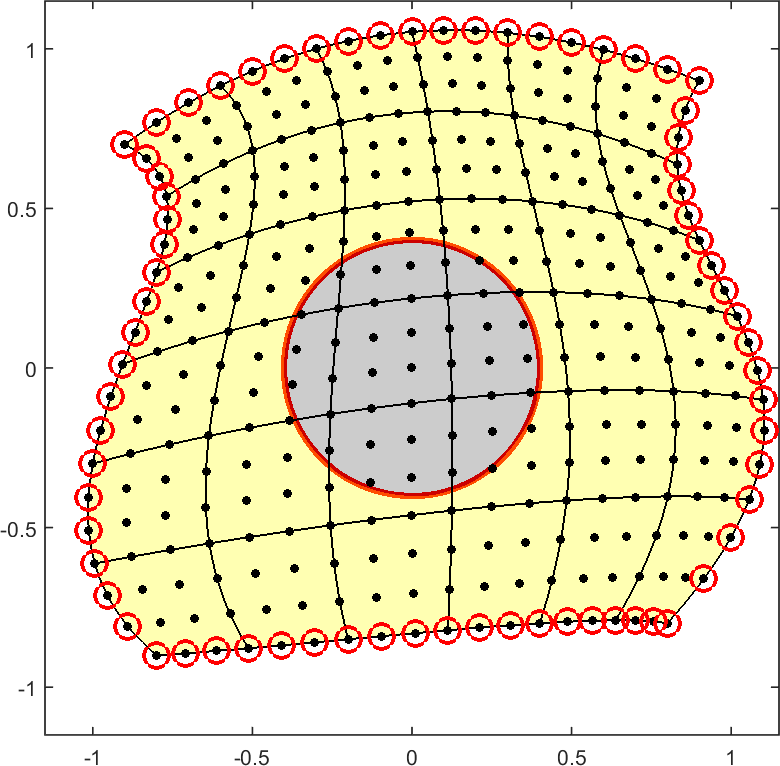}}

\caption{\label{fig:Example3Mesh-2} Background meshes for the bi-material
test case. The outer element edges conform with the boundary. The
interface between the materials however, is within the elements. }
\end{figure}

The problem statement is an infinite domain consisting of material
$A$ with a circular inclusion of material $B$. As before, two different
background meshes are chosen, both are depicted in Fig~\ref{fig:Example3Mesh-2}(a)
and Fig~\ref{fig:Example3Mesh-2}(b). The circle with the radius
$a=0.4$ is centered in the origin. On the entire outer domain boundary,
represented by red circles in in Fig~\ref{fig:Example3Mesh-2}, Dirichlet
boundary conditions are imposed using the analytical solution given
in Eqs.~(\ref{eq:DispBiMatur}) and (\ref{eq:DispBiMatut}).

The solution for the displacement components $u_{x}$ and $u_{y}$
in polar coordinates is given by \cite{Sukumar_2001c,Fries_2006b}
as
\begin{eqnarray}
u_{r}(r,\theta) & = & \left\{ \begin{alignedat}{1}\left(1-\frac{b^{2}}{a^{2}}\right)\alpha+\frac{b^{2}}{a^{2}}r, & \qquad0\leq r\leq a\\
\left(r-\frac{b^{2}}{r}\right)\alpha+\frac{b^{2}}{r}, & \qquad a\leq r\leq b,
\end{alignedat}
\right.\label{eq:DispBiMatur}\\
u_{\theta}(r,\theta) & =0\label{eq:DispBiMatut}
\end{eqnarray}
with the parameter $\alpha$ defined as 
\[
\alpha=\frac{(\lambda_{1}+\mu_{1}+\mu_{2})b^{2}}{(\lambda_{2}+\mu_{2)}a^{2}+(\lambda_{1}+\mu_{1})(b^{2}-a^{2})+\mu_{2}b^{2}}.
\]
$(\theta,r)$ are the standard polar coordinates, $t_{x}$ the traction
in $x$-direction, and $\mu$ the shear modulus. For our example,
plane strain conditions are considered, therefore the Kolosov constant
is defined as $\kappa=3-4\nu$, with $\nu$ being the Poisson ratio.
For the numerical computations, the material constants for material
$A$ and material $B$ are chosen as $\{E_{1}=10,\nu_{1}=0.3\}$ and
$\{E_{2}=1,\nu_{2}=0.25\}$, respectively. For all following examples,
the mesh is sequentially subdivided into $\ell_{h}=\{8,16,32,64,128\}$
quadrilateral elements per side. The characteristic element length
is given as $h=1/\ell_{h}$. A standard displacement-based finite
element formulation is used to solve the equations. 

\begin{figure}
\centering
\subfigure[Cartesian background mesh]{\includegraphics[width=6.5cm]{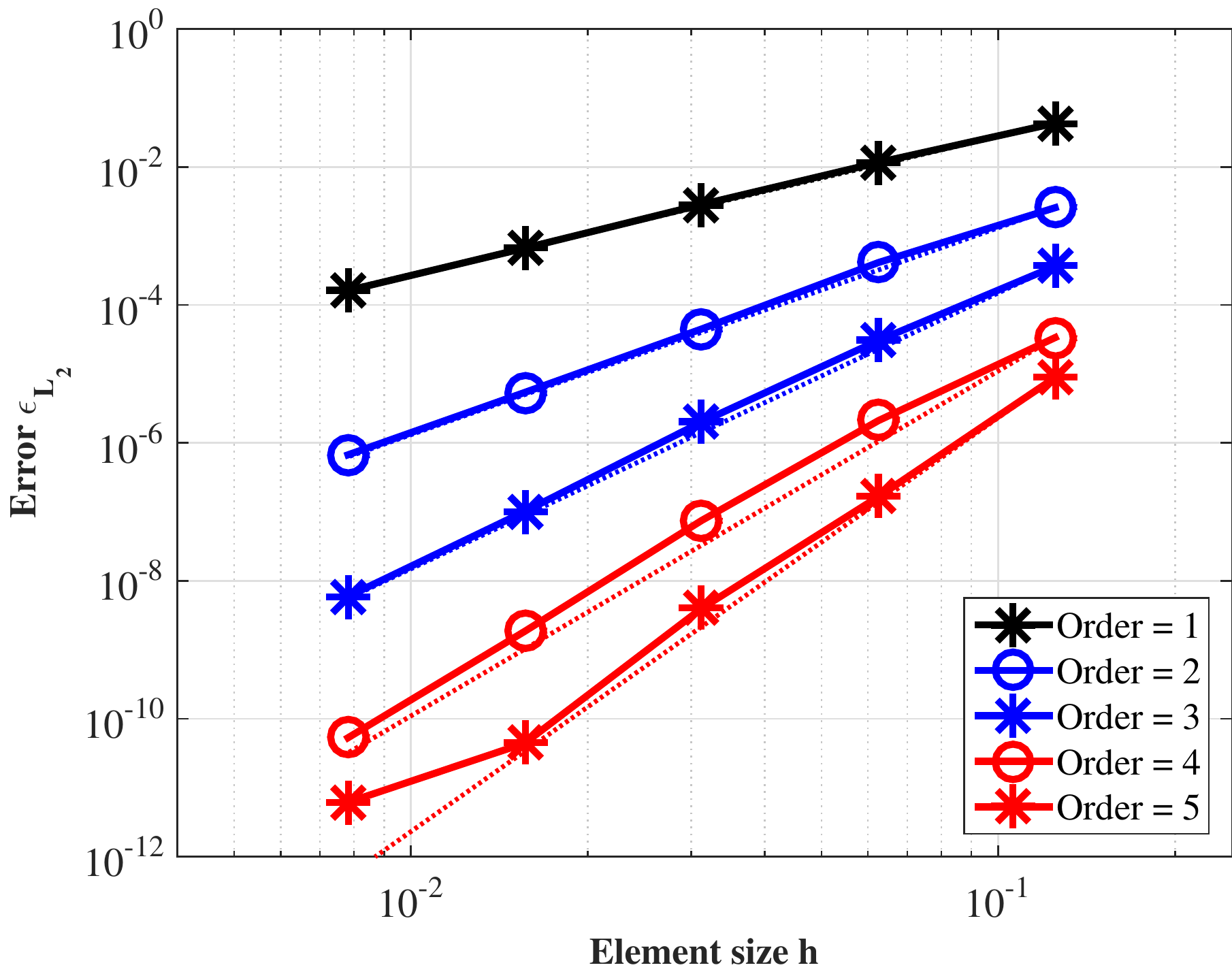}}$\qquad$\subfigure[Deformed background mesh]{\includegraphics[width=6.5cm]{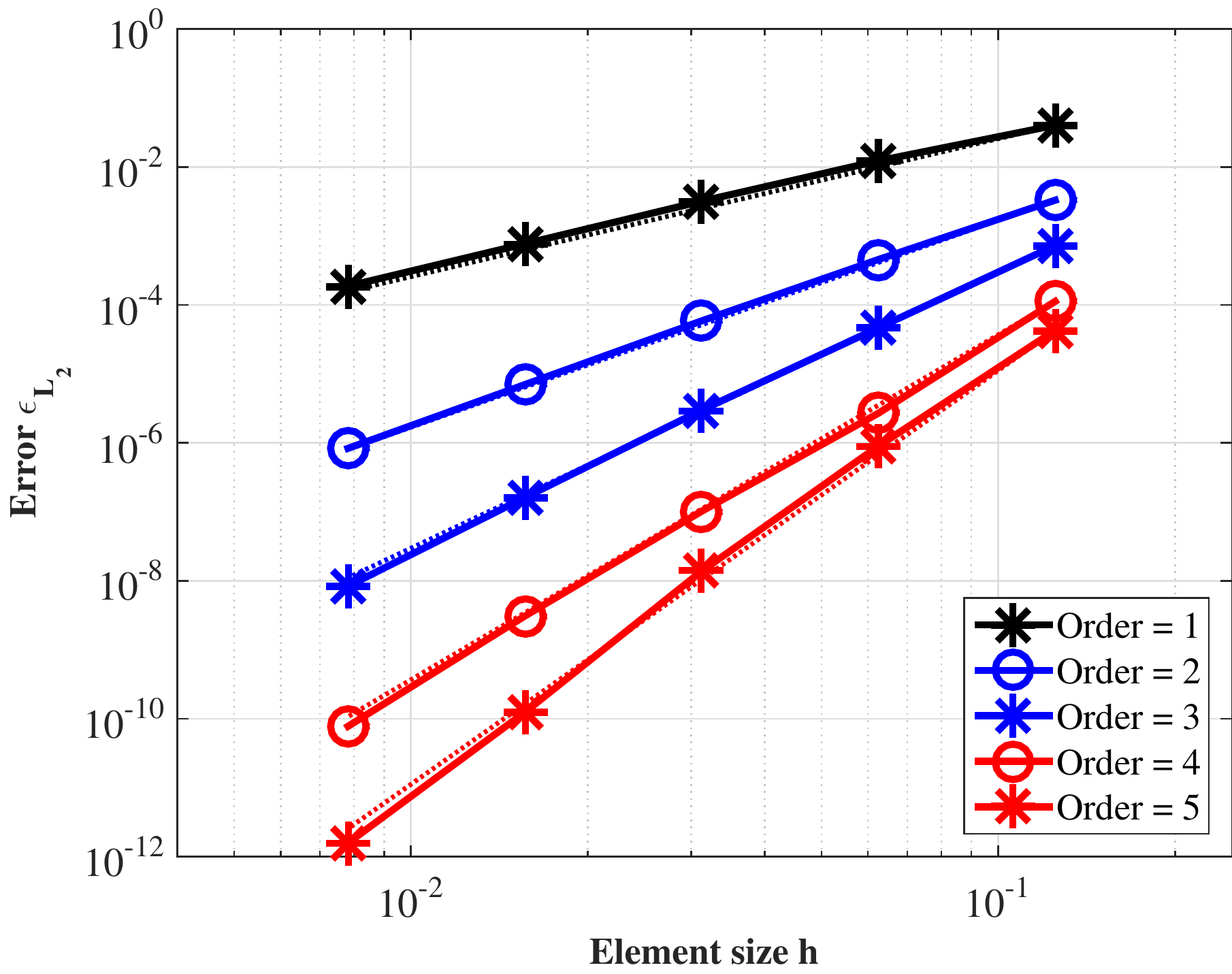}}

\caption{\label{fig:Example3MeshResultConvL2-1} Convergence results for the
bi-material problem with the error measured in the $L_{2}$-norm.}
\end{figure}

The first convergence study measures the displacement error in the
$L_{2}$-norm, $\varepsilon_{L_{2}}$, as introduced in Eq.~(\ref{eq:L2Norm}).
The dotted lines show the theoretically optimal convergence rates.
Fig.~\ref{fig:Example3MeshResultConvL2-1} displays the $L_{2}$-error
as a function of the characteristic element size $h$. The results
validate the proposed approach and show excellent agreement between
the optimal and the calculated convergence rates.

\begin{figure}
\centering
\subfigure[Cartesian background mesh]{\includegraphics[width=6.5cm]{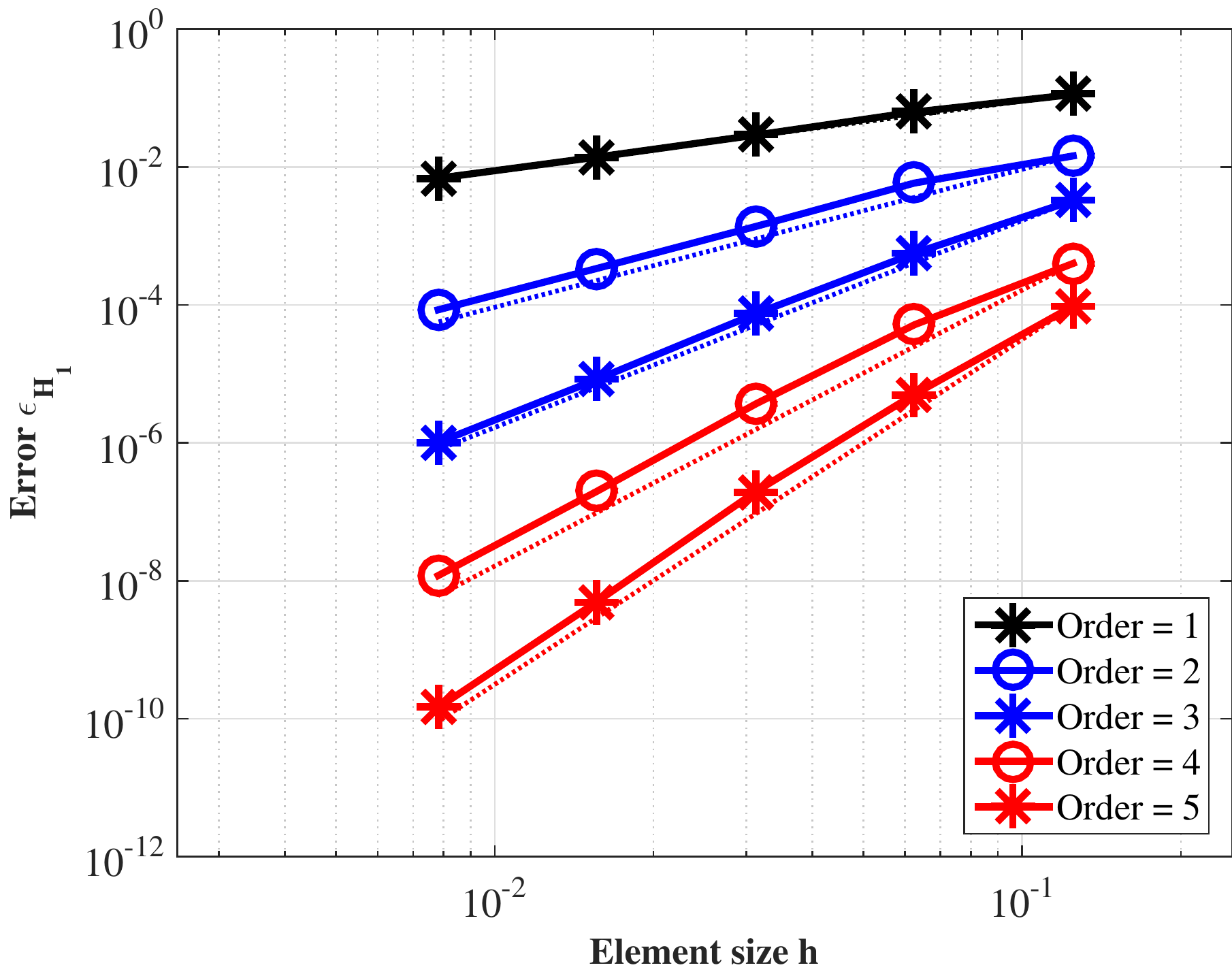}}$\qquad$\subfigure[Deformed background mesh]{\includegraphics[width=6.5cm]{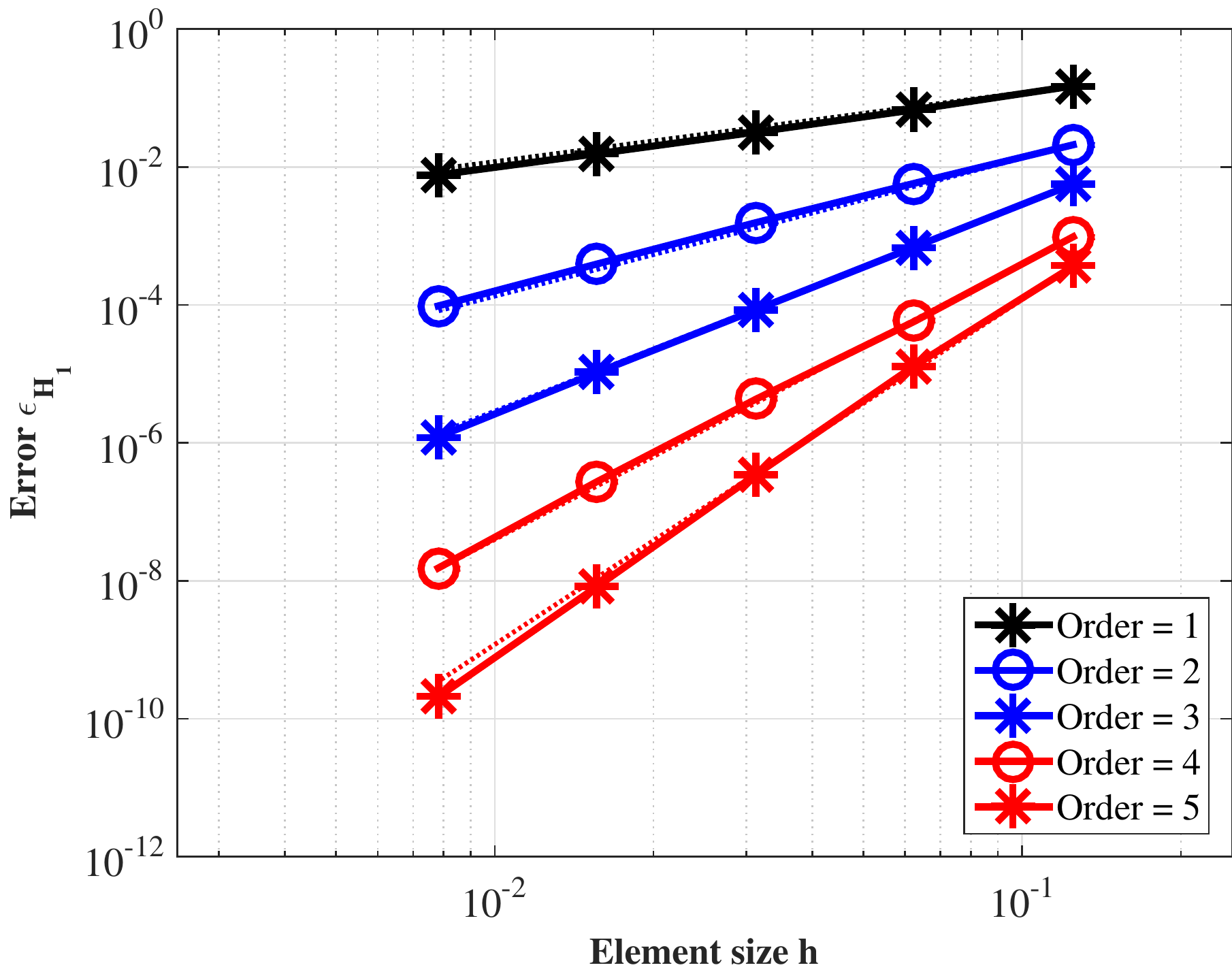}}

\caption{\label{fig:Example1MeshResultConv-1-2-1} Convergence results for
the bi-material problem with the error measured in the $H_{E}$-norm.}
\end{figure}

Also the error in strain energy is calculated using the analytical
solution as presented in \cite{Sukumar_2001c,Fries_2006b}. The results
in Fig.~\ref{fig:Example1MeshResultConv-1-2-1} clearly show optimal
convergence rates also for the strain energy norm.

\begin{figure}
\centering
\subfigure[Cartesian background mesh]{\includegraphics[width=6.5cm]{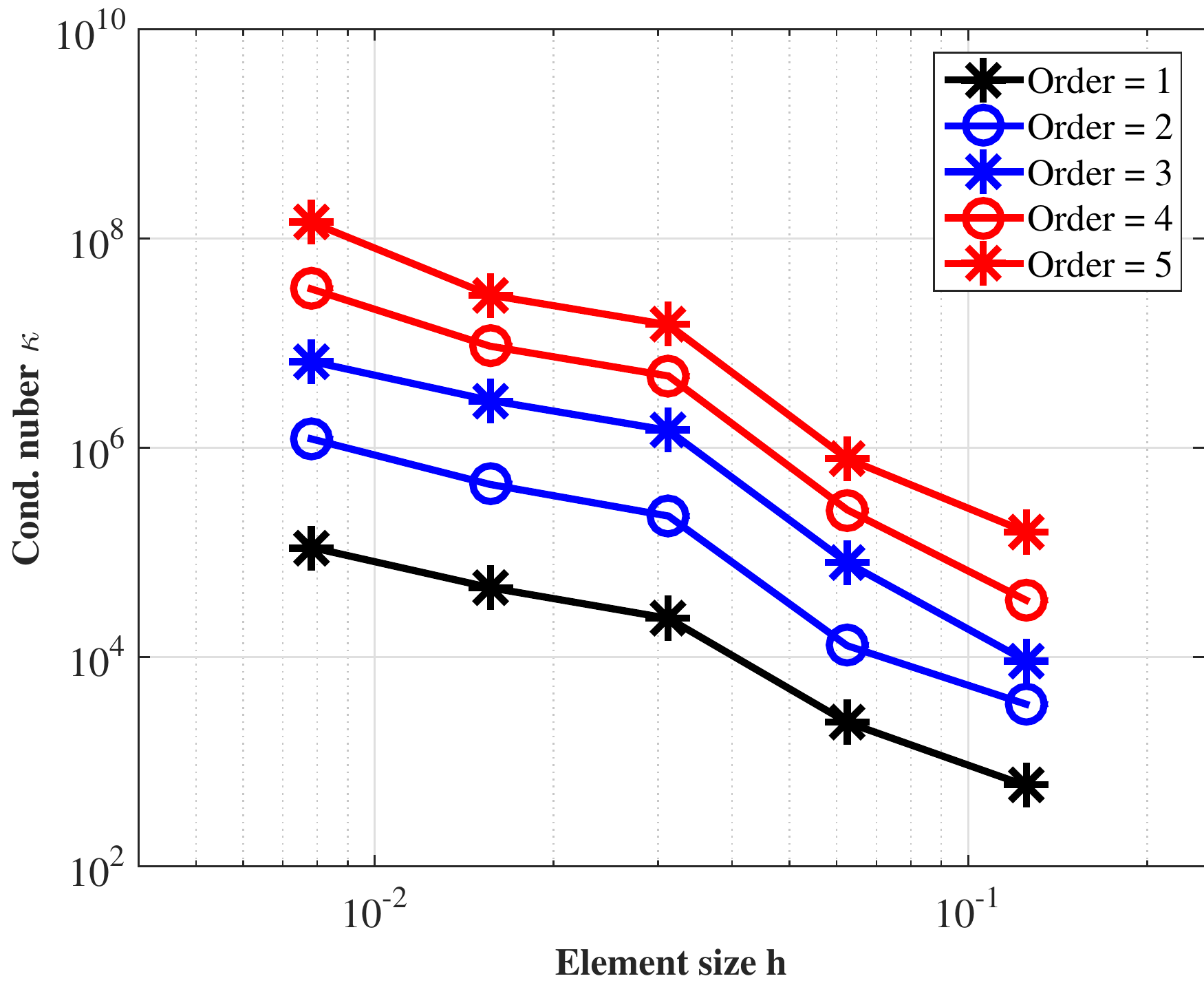}}$\qquad$\subfigure[Deformed background mesh]{\includegraphics[width=6.5cm]{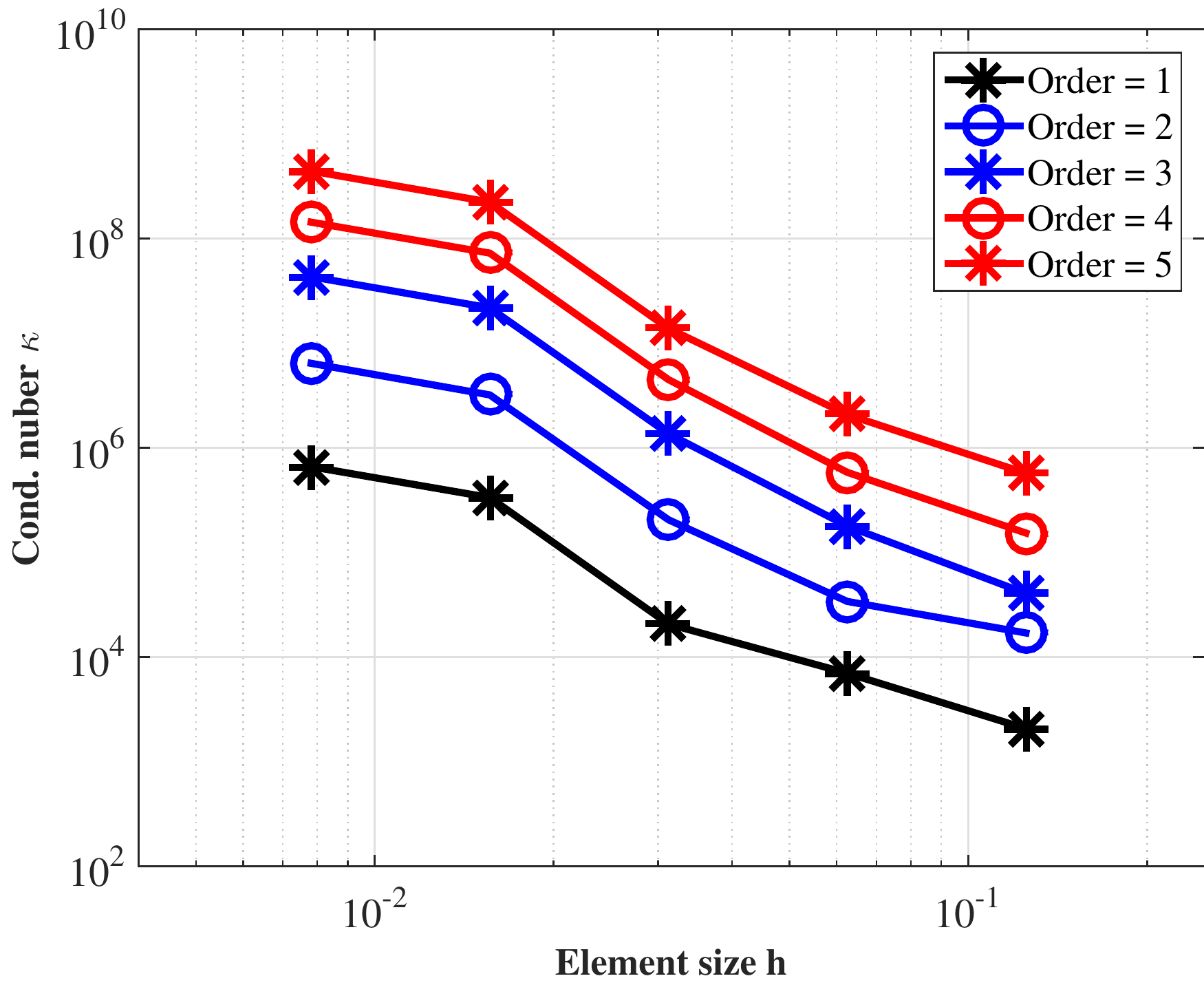}}

\caption{Dependency of the condition number on the element size for different
element orders.\label{fig:Chap5Ex3ResultCondNum}}
\end{figure}

Investigating the condition number $\kappa$ of the system stiffness
matrix and plotting $\kappa$ over the characteristic element length
$h$, Fig.~\ref{fig:Chap5Ex3ResultCondNum}(a) is obtained for the
Cartesian background mesh and Fig.~\ref{fig:Chap5Ex3ResultCondNum}(b)
for the deformed background mesh. It is seen that in this example,
the condition number $\kappa$ is not excessively high, therefore,
no effort was made to manipulate the node positions in the background
mesh. 

\subsection{Infinite plate with a circular hole}

The last example is chosen to illustrate the performance of the proposed
higher-order CDFEM in the context of fictitious domain methods where
the (curved) domain boundary is described implicitly using the level
set function. The example is, as before, a benchmark problem for higher
order finite element procedures, see e.g.~\cite{Cheng_2009a}, \cite{Daux_2000b}.
The error in the displacements (measured in a $L_{2}$-norm), the
error in displacement derivatives (measured in $H_{E}$-norm) and
also the condition number of the equation system are investigated.

\begin{figure}
\centering
\subfigure[Cartesian background mesh]{\includegraphics[width=6.5cm]{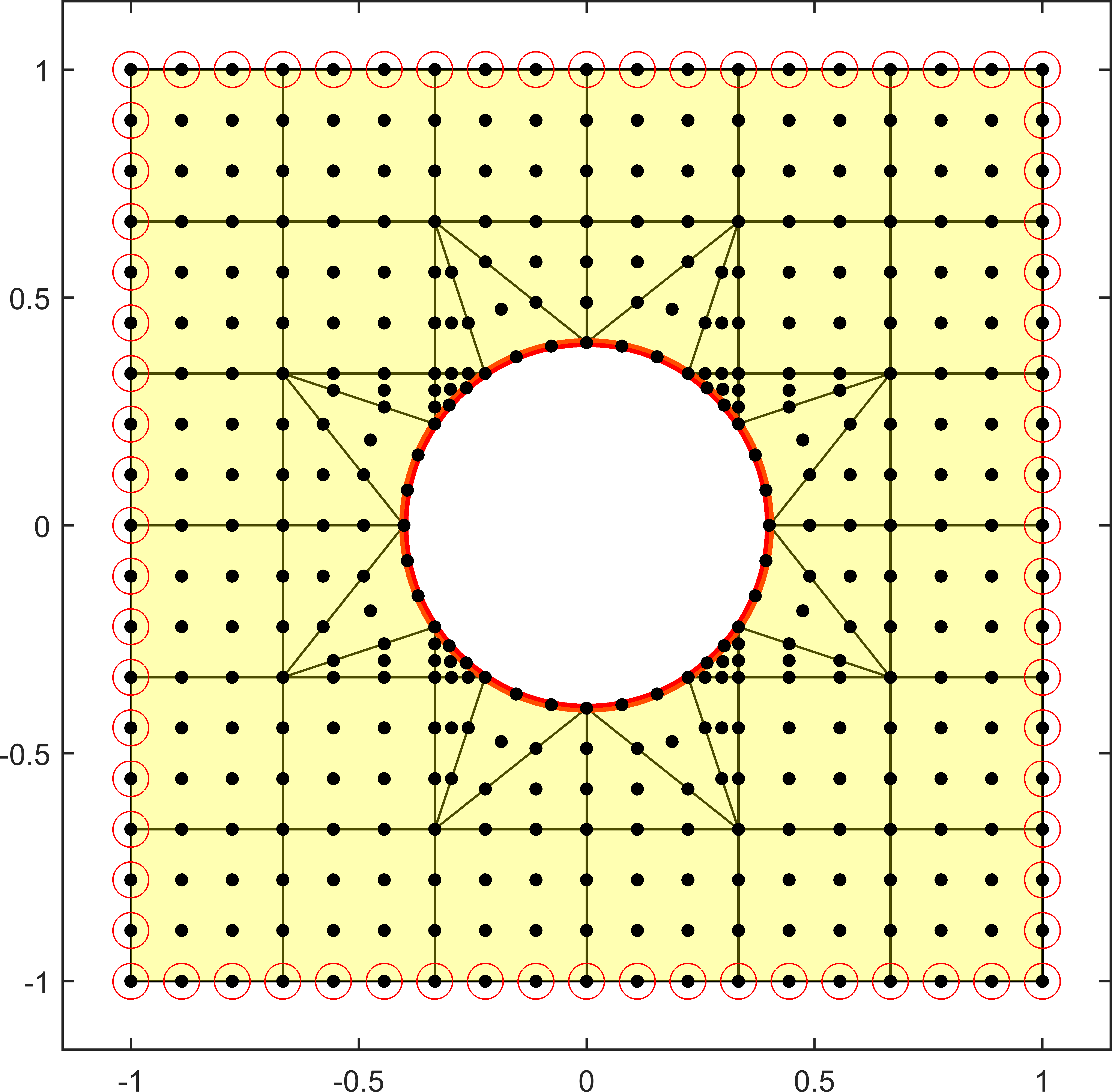}}$\qquad$\subfigure[Deformed mesh]{\includegraphics[width=6.5cm]{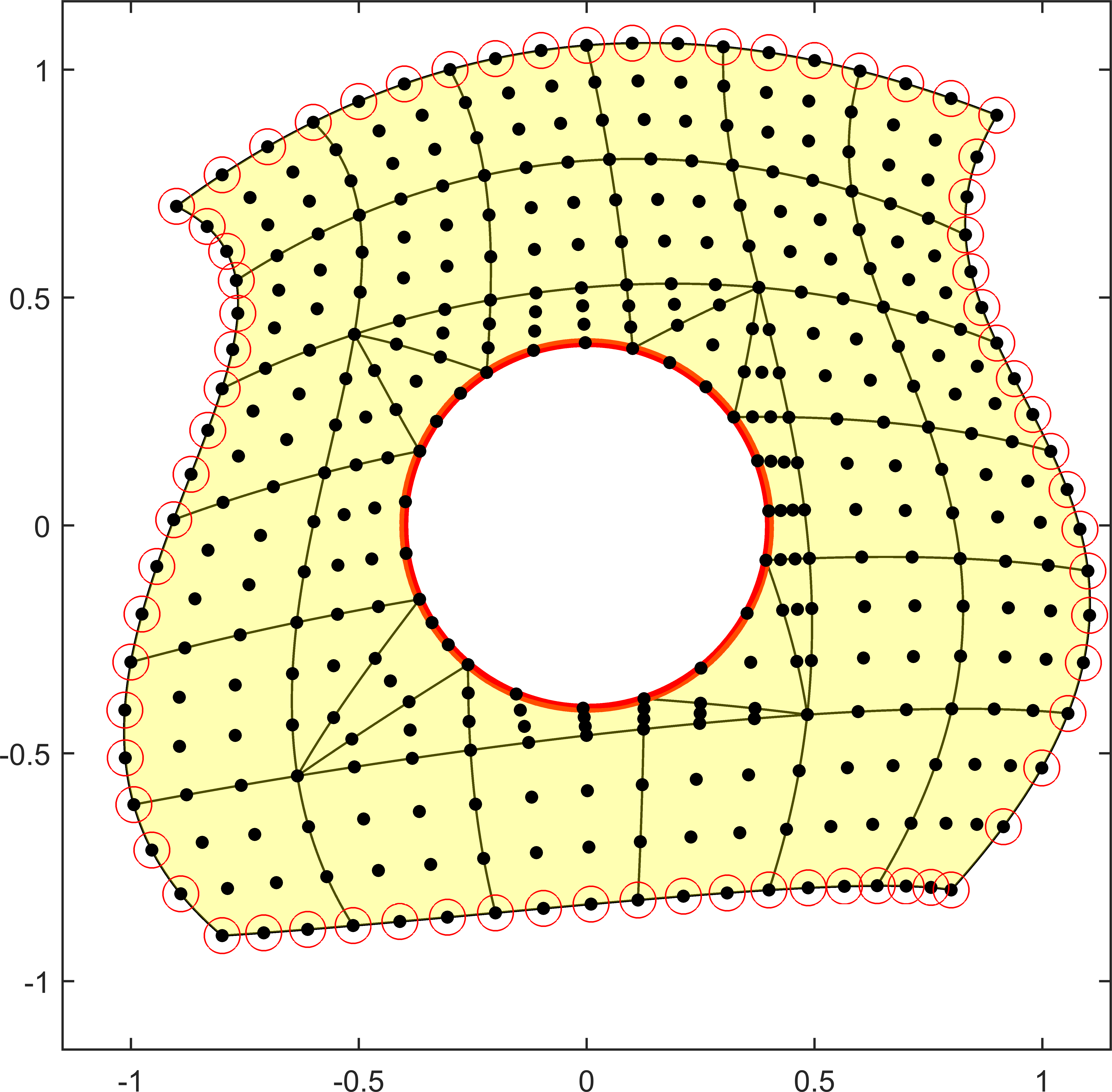}}

\caption{Structured (a) and unstructured mesh (b) with Dirichlet boundary conditions.
\label{fig:Example3Mesh}}
\end{figure}

The problem statement is an infinite domain with a circular inclusion
and unidirectional tension in $x$-direction. As before, two different
background meshes are chosen, a Cartesian background mesh and a deformed
background mesh, both are depicted in Fig.~\ref{fig:Example3Mesh}(a)
and (b). The circular hole with the radius $a=0.4$ is centered in
the origin. Instead of traction boundary conditions, on the entire
outer domain boundary, represented by red circles in Fig~\ref{fig:Example3Mesh},
Dirichlet boundary conditions are imposed using the analytical solution
given in Eq.~(\ref{eq:PlateHoleux}) and Eq.~(\ref{eq:PlateHoleuy}).

The solution for the displacement components $u_{x}$ and $u_{y}$
in Cartesian coordinates is given in \cite{Szabo_1991a} as
\begin{alignat}{1}
u_{x}(r,\theta) & =\frac{t_{x}a}{8\mu}\left[\frac{r}{a}(\kappa+1)\cos\theta+\frac{2a}{r}((1+\kappa)\cos\theta+\cos3\theta)-\frac{2a^{3}}{r^{3}}\cos3\theta\right]\label{eq:PlateHoleux}\\
u_{y}(r,\theta) & =\frac{t_{x}a}{8\mu}\left[\frac{r}{a}(\kappa-3)\sin\theta+\frac{2a}{r}((1-\kappa)\sin\theta+\sin3\theta)-\frac{2a^{3}}{r^{3}}\sin3\theta\right]\label{eq:PlateHoleuy}
\end{alignat}
with $(\theta,r)$ being the standard polar coordinates, $t_{x}$
the traction in $x$-direction, and $\mu$ the shear modulus. For
our example, plane strain conditions are considered, therefore the
Kolosov constant is given as $\kappa=3-4\nu$, with $\nu$ as Poisson
ratio. For the following numerical computations, the material constants
are chosen as $E=10^{4}$ and $\nu=0.3$. 

\begin{figure}
\centering
\subfigure[Cartesian background mesh]{\includegraphics[width=6.5cm]{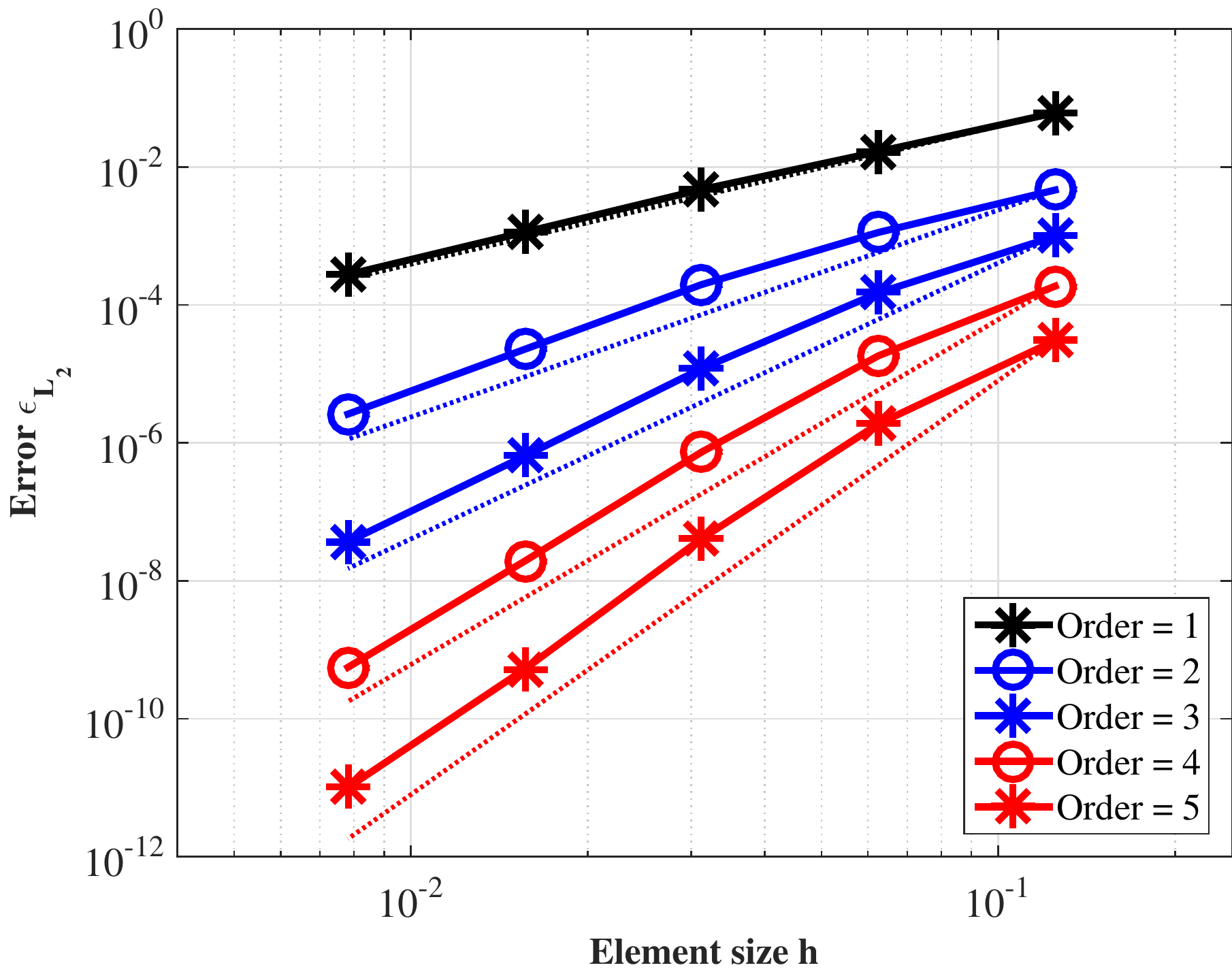}}$\qquad$\subfigure[Deformed background mesh]{\includegraphics[width=6.5cm]{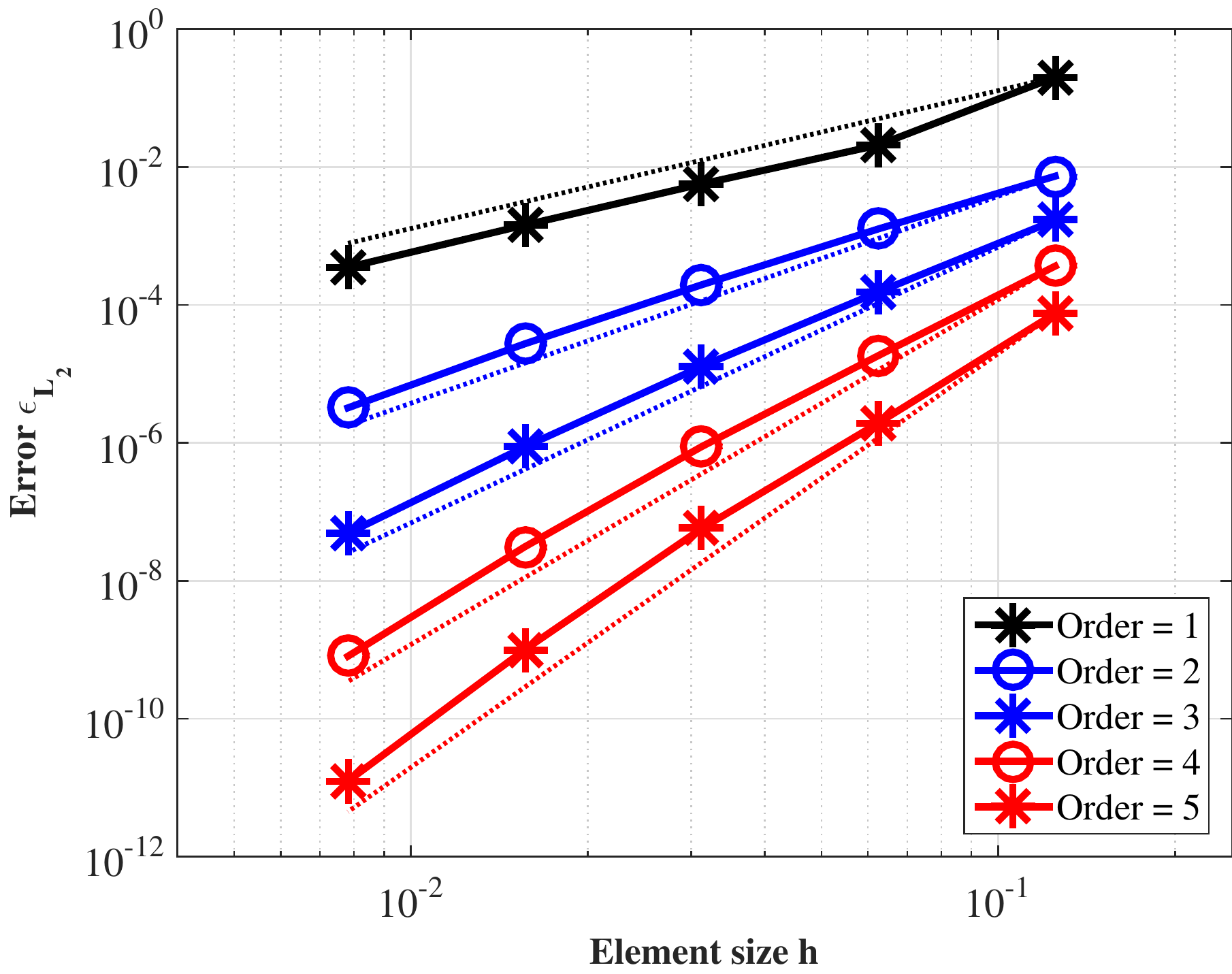}}

\caption{Convergence study for the displacements in the $L_{2}$-norm, Eq.~(\ref{eq:L2Norm}).\label{fig:Example3MeshResultConvL2}}
\end{figure}

As before, the mesh is sequentially subdivided in $\ell_{h}=\{8,16,32,64,128\}$
quadrilateral elements per side. The automatically generated elements
inside the inclusion are simply neglected when integrating the weak
form of the partial differential equation. The first convergence study
measures the displacement error in $L_{2}$-norm, $\varepsilon_{L_{2}}$
as introduced in Eq.~(\ref{eq:L2Norm}). The dotted lines show the
theoretically optimal convergence rates. Fig.~\ref{fig:Example3MeshResultConvL2}
displays the error as a function of the characteristic element length
$h$. The results show again excellent agreement between the optimal
and the calculated convergence rates. 

\begin{figure}
\centering
\subfigure[Cartesian background mesh]{\includegraphics[width=6.5cm]{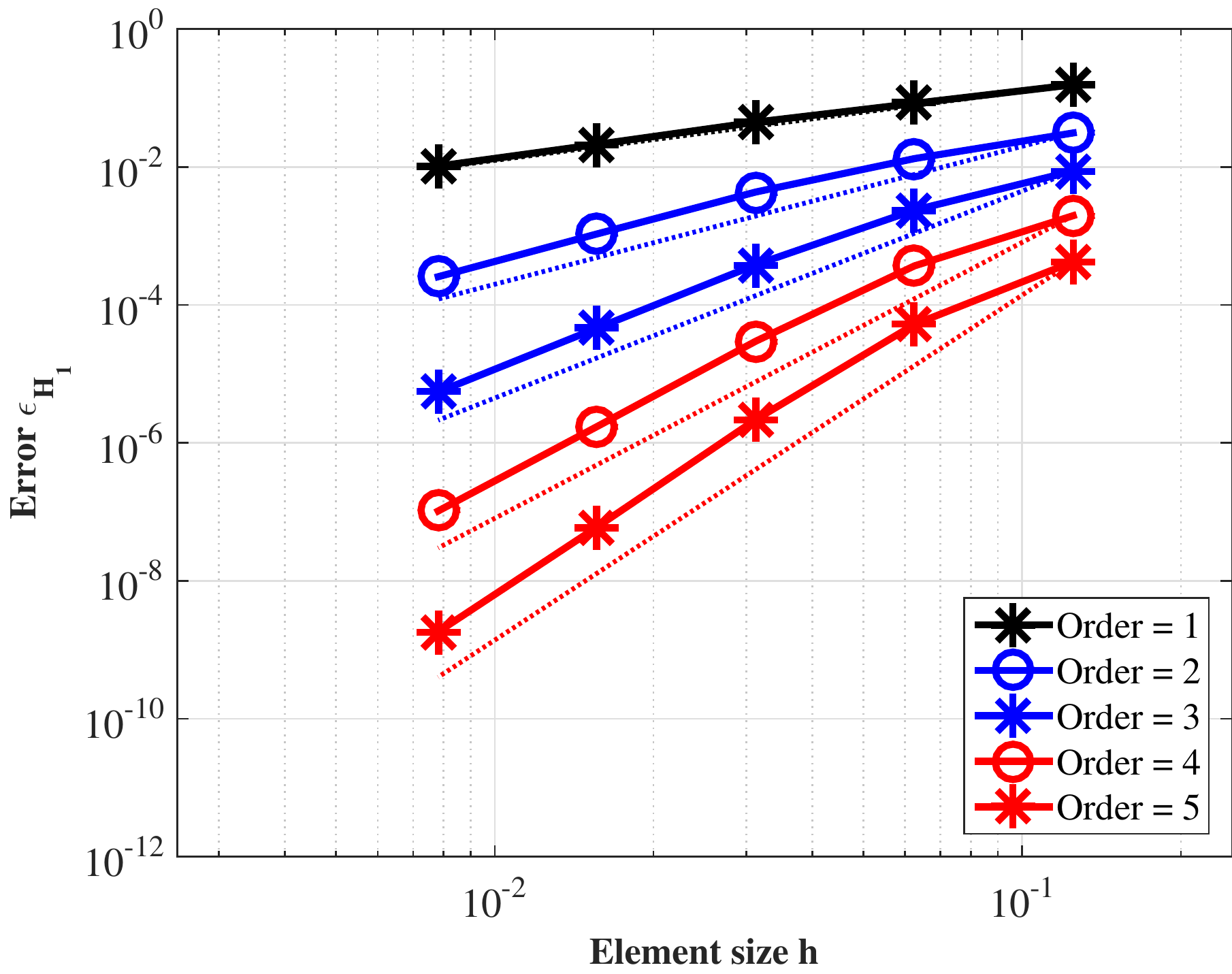}}$\qquad$\subfigure[Deformed background mesh]{\includegraphics[width=6.5cm]{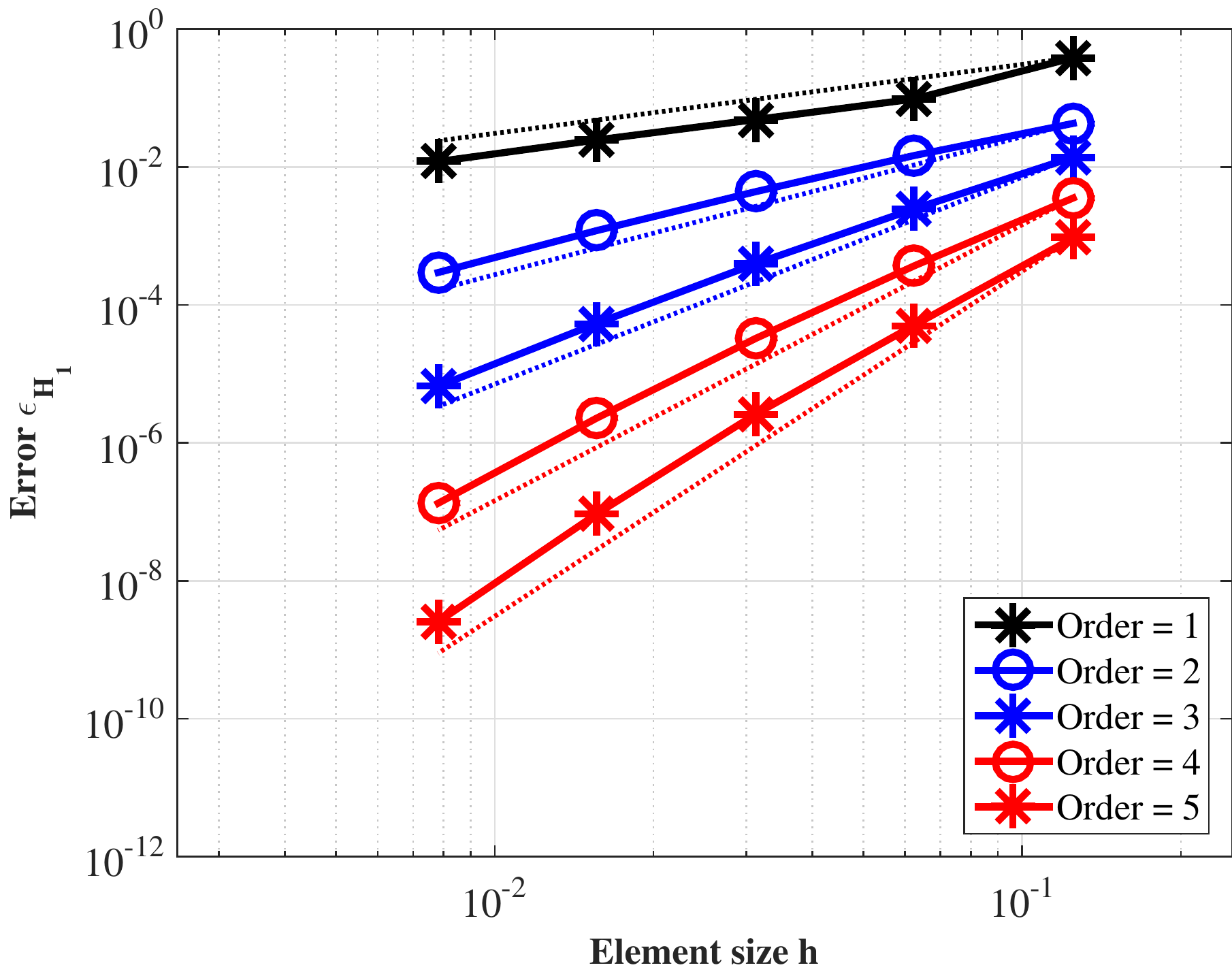}}

\caption{Convergence results for the plate with a circular hole problem in
the $H_{E}$-norm.\label{fig:Chap5Ex4HeNorm}}
\end{figure}

\begin{figure}
\centering
\subfigure[Cartesian background mesh]{\includegraphics[width=6.5cm]{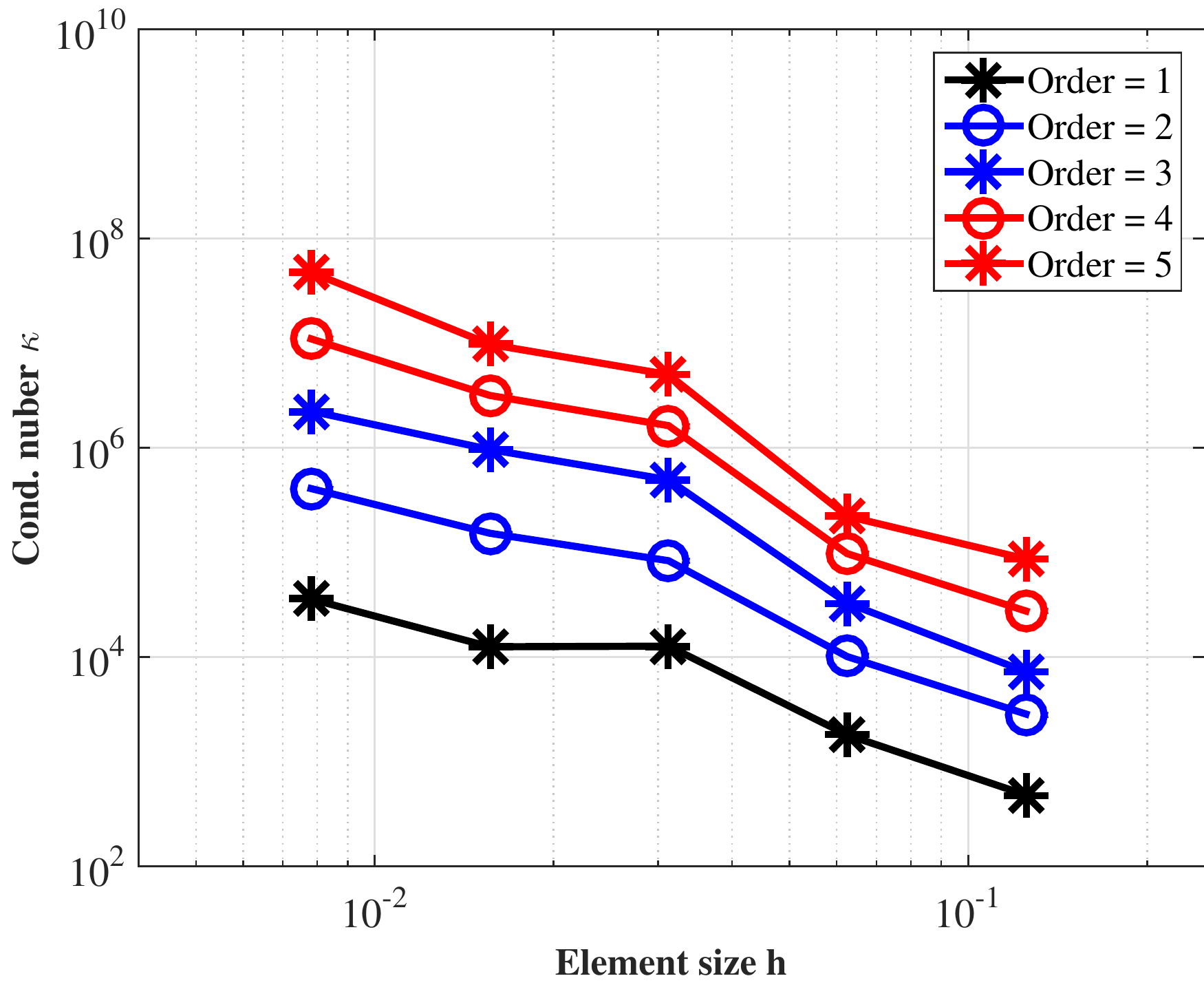}}$\qquad$\subfigure[Deformed background mesh]{\includegraphics[width=6.5cm]{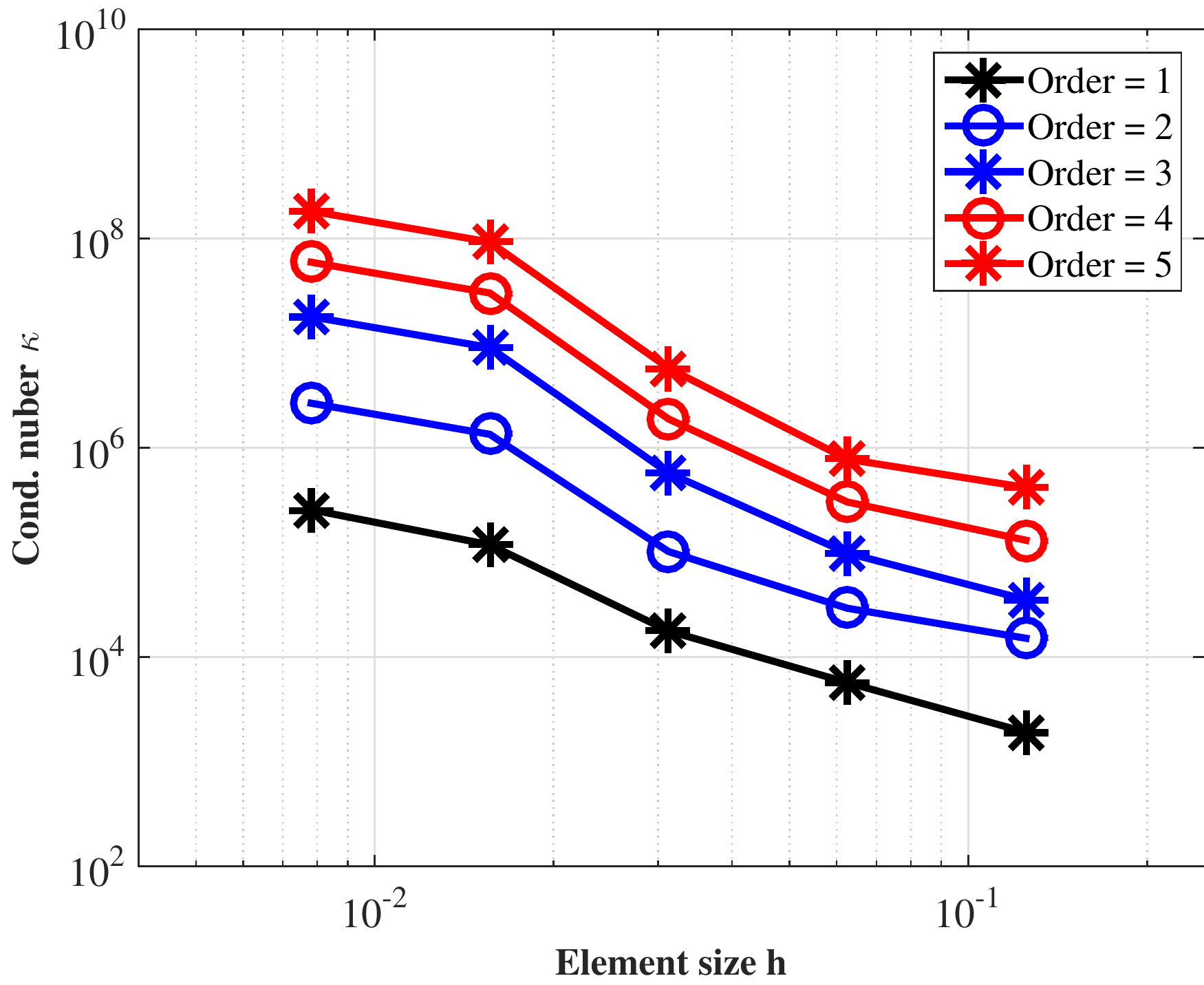}}

\caption{Dependency of the condition number on the element size for different
element orders.\label{fig:Chap5Ex4ResultCondNum}}
\end{figure}

The exact stresses and strains are given in e.g.~\cite{Cheng_2009a},
\cite{Daux_2000b}. The error in strain energy is calculated using
the error measure in Eq.~(\ref{eq:HENorm}). The results are shown
in Fig.~\ref{fig:Chap5Ex4HeNorm} and are again optimal. The condition
number for the Cartesian and the deformed background mesh are plotted
in Fig.~\ref{fig:Chap5Ex4ResultCondNum}(a) and Fig.~\ref{fig:Chap5Ex4ResultCondNum}(b).

\section{Conclusions, Limitations and Perspectives\label{sec:Conclusions}}

In this work, a higher-order conformal decomposition finite element
method (CDFEM) has been proposed. It combines the advantages of the
FEM on body fitted meshes (application of boundary and interface conditions,
numerical integration) with those of embedded domain methods (no explicit
manual mesh generation). The key-aspect is the automatic decomposition
of cut background elements into higher-order conforming sub-elements.

Therefore, a two-step procedure is suggested: In the first part, the
implicit interfaces described by the level set method are detected
and meshed. Then, based on the topological situation in the cut element,
higher-order elements are mapped on the two sides of the interface.
It is shown that the location of the inner element nodes, depending
on the concrete mapping, have important implications for the achieved
convergence rates. 

The reconstruction and remeshing scheme is executed
in the reference element, therefore also non-structured background
meshes can be used. Under restrictions on the curvature of the level
set function, all topological situations can be covered. Used in a
\emph{FEM} simulation, the automatically determined sub-elements may
lead to an ill-conditioning of the system matrix. Several solution
strategies will be investigated in future works. Herein, none of these
techniques have been employed and, nevertheless, excellent have been
obtained with direct solvers.

The proposed approach shows excellent agreement between optimal and
calculated convergence rates. It is very flexible with respect to
the employed element types and not limited to structured background
meshes. It shows a stable behaviour even for geometrically complex
interfaces. The method is able to automatically remesh a broad class
of problems with an implicitly sharp interface representation. A number
of extensions for the method can be conceived, for example, evolving
interfaces that occur in plasticity or two-phase flows. Also an extension
to three-dimensional problems is currently in progress. It is also
planned to combine the approach with adaptive refinements such that
even more complex geometrical situations can be covered automatically.

\bibliographystyle{plain}
\bibliography{ArxivPrePrint}

\end{document}